\documentclass{amsart}

\usepackage{amsmath, amsthm, amssymb}
\usepackage{amsfonts}
\usepackage{graphicx}
\usepackage{epic}

\newtheorem{Theorem}{Theorem}

\newtheorem{Definition}{Definition}
\newtheorem{Example}{Example}

\begin{document}

\title{$k$-Ribbon Fibonacci Tableaux}         

\author{Naiomi Cameron}
\address{Department of Mathematical Sciences, Lewis \& Clark College, Portland, Oregon 97219}
\email{ncameron@lclark.edu} 

\author{Kendra Killpatrick}
\address{Department of Mathematics, Pepperdine University, Malibu, California 90265} \email{Kendra.Killpatrick@pepperdine.edu}
      


\date{}          


\begin{abstract}  We extend the notion of $k$-ribbon tableaux to the Fibonacci lattice, a differential poset defined by R. Stanley in 1975.  Using this notion, we describe an insertion algorithm that takes $k$-colored permutations to pairs of $k$-ribbon Fibonacci tableaux of the same shape, and we demonstrate a color-to-spin property, similar to that described by Shimozono and White for ribbon tableaux.  We give an evacuation algorithm which relates the pair of $k$-ribbon Fibonacci tableaux obtained through the insertion algorithm to the pair of $k$-ribbon Fibonacci tableaux obtained using Fomin's growth diagrams.  In addition, we present an analogue of Knuth relations for $k$-colored permutations and $k$-ribbon Fibonacci tableaux.
\end{abstract}

\maketitle

\section{Introduction}

Young's lattice of partitions ordered by inclusion and its generalization $Y^r$ has been well-studied.  It is one of the prime examples of an $r$-differential poset, as defined by R. Stanley in 1988 \cite{St1}.  Stanley showed that for any $r$-differential poset $P$, 

\begin{equation} 
\sum_{\lambda \in P_n} e(\lambda)^2 = r^n n!
\label{r^n}
\end{equation}
where $P_n$ is the set of elements of rank $n$ and $e(\lambda)$ is the number of chains in $P$ from $\hat{0}$ to $\lambda$.  

For the case when $r=1$, there are several combinatorial proofs for (\ref{r^n}), including the Schensted insertion algorithm which takes permutations in $S_n$ to pairs of standard Young tableaux $(P, Q)$ of the same shape.   For the case when $r=2$, Shimozono and White \cite{ShW} describe the domino (or 2-ribbon) insertion algorithm given by Barbasch and Vogan \cite{BaV} and Garfinkle \cite{Gar} which takes 2-colored permutations to pairs of domino tableaux and proved that the total color of the permutation is equal to the spin of the domino tableaux.  In 2002, Shimozono and White extended these ideas for arbitrary values of $r$ by describing color-to-spin insertion algorithms for ribbon tableaux \cite{SW2}.

The other classic example of an $r$-differential poset is the Fibonacci lattice $Z(r)$, described by R. Stanley \cite{St2}.  When $r=1$, the identity in (\ref{r^n}) can be shown by means of a Schensted-like insertion which takes permutations to pairs of Fibonacci tableaux of the same shape \cite{Rob}.  In 2006, the authors developed a notion of domino tableaux for the Fibonacci lattice $Z(2)$ and described an insertion algorithm on 2-colored permutations with a color-to-spin property similar to that described by Shimozono and White \cite{CamKill}.  In the present paper, we extend the notion of domino Fibonacci tableaux and its consequences to the general Fibonacci lattice, in the form of  $k$-ribbon Fibonacci tableaux.

In the next section, we provide the necessary background for the rest of the paper, while in Section \ref{kribbon}, we formally define $k$-ribbon Fibonacci tableaux.  The $k$-ribbon Fibonacci insertion and evacuation algorithms are fully described in Sections \ref{insertion} and \ref{evacuation}.  A geometric interpretation of the insertion algorithm is presented in Section \ref{geometric}, which helps to illuminate the major results of the insertion in Section \ref{theorems} and makes the ``color-to-spin" property of the insertion a relatively straightforward matter in Section \ref{colortospin}.  Finally, Section \ref{knuth} is devoted to describing the analogue of Knuth relations for $k$-colored permutations and $k$-ribbon Fibonacci tableaux.

\section{Background and Definitions}      
\label{background}

In this section we give the necessary background and definitions for the theorems in this paper.  The interested reader is encouraged to read Chapter 5 of {\it{The Symmetric Group, 2nd Edition}} by Bruce Sagan \cite{Sag} for general reference. 

Stanley \cite{St1} gave the following definition of an {\it{$r$-differential poset}}:

\begin{Definition}
An {\it {$r$-differential poset}} is a poset which satisfies the following three conditions:
\end{Definition}

\begin{enumerate}
\item  $P$ has a $\hat{0}$ element, is graded and is locally finite.
\item  If $x \neq y$ and there are exactly $m$ elements in $P$ which are covered by $x$ and by $y$, then there are exactly $m$ elements in $P$ which cover both $x$ and $y$.  
\item  For $x \in P$, if $x$ covers exactly $m$ elements of $P$, then $x$ is covered by exactly $m+r$ elements of $P$.
\end{enumerate}

The classic example of a $1$-differential poset is Young's lattice, which is the poset of the set of partitions together with the binary relation $\lambda \leq \mu$ if and only if $\lambda_i \leq \mu_i$ for all $i$.

A second kind of $1$-differential poset is the Fibonacci differential poset.  The general definition of a Fibonacci $r$-differential poset was given by Richard Stanley in \cite{St2} (Definition 5.2).  

Let $A = \{ 1_1, 1_2, \dots, 1_r, 2 \}$ and let $A^*$ be the set of all finite words $a_1 a_2 \cdots a_l$ of elements of $A$ (including the empty word).

\begin{Definition}
The {\it {Fibonacci differential poset}} $Z(r)$ has as its set of elements the set of words in $A^*$.  If $w \in Z(r)$, then define $z$ to be covered by $w$ (i.e. $z \lessdot w$) in $Z(r)$ if either:
\end{Definition}

\begin{enumerate}
\item  $z$ is obtained from $w$ by changing a $2$ to a $1_i$ for some $i$ if the only letters to the left of this $2$ are also $2$'s, or
\item $z$ is obtained from $w$ by deleting the leftmost $1$ of any type.  
\end{enumerate}

The first four rows of the Fibonacci lattice $Z(2)$ are shown below.  
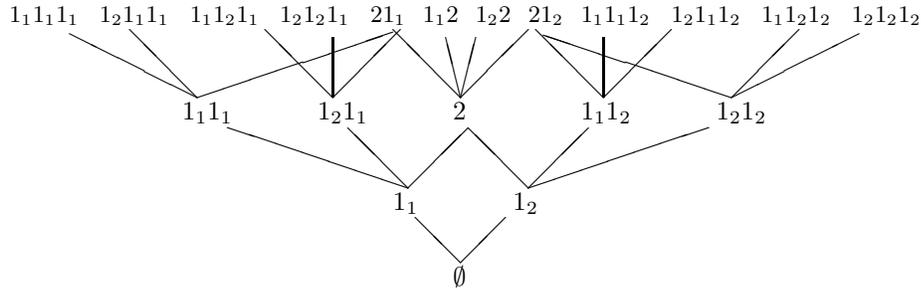
\begin{center}
\setlength{\unitlength}{1cm}
\begin{center}
\begin{picture}(9,4)(1,1)
\put(5.9,1){$\emptyset$}
\put(5.1,2){$1_1$}
\put(6.7,2){$1_2$}
\put(6, 1.3){\line(-1,1){.6}}
\put(6, 1.3){\line(1,1){.6}}
\put(2.3, 3.2){$1_1 1_1$}
\put(4.1, 3.2){$1_2 1_1$}
\put(5.9,3.2){2}
\put(9.4, 3.2){$1_2 1_2$}
\put(7.6, 3.2){$1_1 1_2$}

\put(5.3, 2.3){\line(-1,1){.8}}
\put(5.3, 2.3){\line(1,1){.8}}
\put(5.3, 2.3){\line(-3,1){2.4}}
\put(6.9, 2.3){\line(-1,1){.8}}
\put(6.9, 2.3){\line(3,1){2.4}}
\put(6.9, 2.3){\line(1,1){.8}}

\put(0,4.5){\small{$1_1 1_1 1_1$}}
\put(1.2, 4.5){\small{$1_2 1_1 1_1$}}
\put(2.4, 4.5){\small{$1_1 1_2 1_1$}}
\put(3.6, 4.5){\small{$1_2 1_2 1_1$}}
\put(4.8, 4.5){\small{$21_1$}}
\put(5.5, 4.5){\small{$1_1 2$}}
\put(6.2, 4.5){\small{$1_2 2$}}
\put(6.9, 4.5){\small{$2 1_2$}}
\put(7.6, 4.5){\small{$1_1 1_1 1_2$}}
\put(8.8, 4.5){\small{$1_2 1_1 1_2$}}
\put(10.0, 4.5){\small{$1_1 1_2 1_2$}}
\put(11.2, 4.5){\small{$1_2 1_2 1_2$}}

\put(2.5, 3.5){\line(-2,1){1.7}}
\put(2.5, 3.5){\line(-1,1){.9}}
\put(2.5, 3.5){\line(3,1){2.6}}

\put(4.3, 3.5){\line(-1,1){.9}}
\put(4.3, 3.5){\line(0,1){.8}}
\put(4.3, 3.5){\line(1,1){.9}}

\put(6, 3.5){\line(-1,1){.9}}
\put(6, 3.5){\line(-1,4){.2}}
\put(6, 3.5){\line(1,4){.2}}
\put(6, 3.5){\line(1,1){.9}}

\put(7.9, 3.5){\line(-1,1){.9}}
\put(7.9, 3.5){\line(0,1){.8}}
\put(7.9, 3.5){\line(1,1){.9}}

\put(9.6, 3.5){\line(-3,1){2.5}}
\put(9.6, 3.5){\line(1,1){.9}}
\put(9.6, 3.5){\line(2,1){1.7}}

\end{picture}
\end{center}
\end{center}
\vspace{.5in}

Fomin \cite{Fom} gave a general method for representing a permutation with a square diagram and then using a growth function to create a pair of saturated chains in a differential poset.  In particular, Fomin's method can be applied to the square diagram of a $k$-colored permutation to create a pair of saturated chains in $Z(k)$, giving a further proof of Stanley's result (\ref{r^n}).

Given a permutation in $S_n$, we can create a $k$-colored permutation by assigning to each element one of $k$ colors.  If element $x_i$ in the permutation is colored by color $j$ we write $x_i^j$.  Note that the set of all $k$-colored permutations certainly includes permutations which do not use all $k$-colors.

  For a $k$-colored permutation written in two line notation:
\[
\begin{array}{cccccc}
\pi&=&1&2&\cdots&n\\
 & &x_1^{j_1}&x_2^{j_2}&\cdots&x_n^{j_n}
\end{array}
\]
with each $x_i$ colored by color $j_i$ for $1 \leq j_i \leq k$, we create a square diagram by placing an $X^{j_1}$ in column $i$ and row $x_i$ (indexed from left to right, bottom to top) if  \shortstack{$i$\\$x_i^{j_i}$} is a column in the permutation $\pi$.  For example, consider the permutation
\[
\pi = \begin{array}{ccccccc}
1&2&3&4&5&6&7\\
2^3&7^1&1^1&5^4&6^3&4^2&3^4
\end{array}
\]
as a $5$-colored permutation.  Then we obtain the following square diagram:

\setlength{\unitlength}{1cm}
\begin{center}
\begin{picture}(7,6.5)(1,1)
\thinlines
\put(1,0){\line(0,1){7}}
\put(2,0){\line(0,1){7}}
\put(3,0){\line(0,1){7}}
\put(4,0){\line(0,1){7}}
\put(5,0){\line(0,1){7}}
\put(6,0){\line(0,1){7}}
\put(7,0){\line(0,1){7}}
\put(8,0){\line(0,1){7}}
\put(1,0){\line(1,0){7}}
\put(1,1){\line(1,0){7}}
\put(1,2){\line(1,0){7}}
\put(1,3){\line(1,0){7}}
\put(1,4){\line(1,0){7}}
\put(1,5){\line(1,0){7}}
\put(1,6){\line(1,0){7}}
\put(1,7){\line(1,0){7}}
\put(1.3, 1.3){$X^3$}
\put(2.3, 6.3){$X^1$}
\put(3.3, .3){$X^1$}
\put(4.3, 4.3){$X^4$}
\put(5.3, 5.3){$X^3$}
\put(6.3, 3.3){$X^2$}
\put(7.3, 2.3){$X^4$}
\end{picture}
\end{center}
\vspace{.5in}

Fomin's method gives a way to translate this square diagram into a pair of saturated chains in $Z(k)$ in the following manner.  Begin by placing $\emptyset$'s along the lower edge and the left edge at each corner.  Label the remaining corners in the diagram by following the rules given below (called a {\it{growth function}}).  If we have

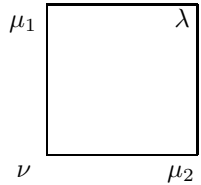
\setlength{\unitlength}{1cm}
\begin{center}
\begin{picture}(2,2)(0,0)
\thinlines
\put(0,0){\line(0,1){2}}
\put(2,0){\line(0,1){2}}
\put(0,0){\line(1,0){2}}
\put(0,2){\line(1,0){2}}
\put(-.4,-.3){$\nu$}
\put(-.5, 1.7){$\mu_1$}
\put(1.6, -.3){$\mu_2$}
\put(1.7, 1.7){$\lambda$}
\end{picture}
\end{center}

\noindent with each side of the square representing a cover relation in $Z(k)$ or an equality, then:

\begin{enumerate}
\item If $\mu_1 \gtrdot \nu$ and $\mu_2 = \nu$ then $\lambda=\mu_1$ (and similarly for $\mu_1$ and $\mu_2$ interchanged).
\item If $\mu_1 \gtrdot \nu$, $\mu_2 \gtrdot \nu$ then $\lambda$ is obtained from $\nu$ by prepending a $2$.
\item If $\mu_1 = \nu = \mu_2$ and the box contains an $X^i$, then obtain $\lambda$ from $\nu$ by prepending a $1_i$.
\item If $\mu_1 = \nu = \mu_2$ and the box does not contain an $X^i$, then $\lambda = \nu$.
\end{enumerate}

By following this procedure on our previous example, we obtain the complete growth diagram:

\setlength{\unitlength}{1.5cm}
\begin{center}
\begin{picture}(7,6.5)(1,1)
\thinlines
\put(1,0){\line(0,1){7}}
\put(2,0){\line(0,1){7}}
\put(3,0){\line(0,1){7}}
\put(4,0){\line(0,1){7}}
\put(5,0){\line(0,1){7}}
\put(6,0){\line(0,1){7}}
\put(7,0){\line(0,1){7}}
\put(8,0){\line(0,1){7}}
\put(1,0){\line(1,0){7}}
\put(1,1){\line(1,0){7}}
\put(1,2){\line(1,0){7}}
\put(1,3){\line(1,0){7}}
\put(1,4){\line(1,0){7}}
\put(1,5){\line(1,0){7}}
\put(1,6){\line(1,0){7}}
\put(1,7){\line(1,0){7}}
\put(.8,-.3){$\emptyset$}
\put(.8,.7){$\emptyset$}
\put(.8,1.7){$\emptyset$}
\put(.8,2.7){$\emptyset$}
\put(.8,3.7){$\emptyset$}
\put(.8,4.7){$\emptyset$}
\put(.8,5.7){$\emptyset$}
\put(.8,6.7){$\emptyset$}
\put(1.8,-.3){$\emptyset$}
\put(2.8,-.3){$\emptyset$}
\put(3.8,-.3){$\emptyset$}
\put(4.8,-.3){$\emptyset$}
\put(5.8,-.3){$\emptyset$}
\put(6.8,-.3){$\emptyset$}
\put(7.8,-.3){$\emptyset$}
\put(1.4, 1.4){$X^3$}
\put(2.4, 6.3){$X^1$}
\put(3.4, .4){$X^1$}
\put(4.4, 4.4){$X^4$}
\put(5.3, 5.3){$X^3$}
\put(6.3, 3.3){$X^2$}
\put(7.4, 2.4){$X^4$}
\put(1.8, .7){$\emptyset$}
\put(2.8, .7){$\emptyset$}
\put(3.7, .7){$1_1$}
\put(4.7, .7){$1_1$}
\put(5.7, .7){$1_1$}
\put(6.7, .7){$1_1$}
\put(7.7, .7){$1_1$}
\put(1.7, 1.7){$1_3$}
\put(2.7, 1.7){$1_3$}
\put(3.7, 1.7){2}
\put(4.7, 1.7){2}
\put(5.7, 1.7){2}
\put(6.7, 1.7){2}
\put(7.7, 1.7){2}
\put(1.7, 2.7){$1_3$}
\put(2.7, 2.7){$1_3$}
\put(3.7, 2.7){2}
\put(4.7, 2.7){2}
\put(5.7, 2.7){2}
\put(6.7, 2.7){2}
\put(7.6, 2.7){$1_42$}
\put(1.7, 3.7){$1_3$}
\put(2.7, 3.7){$1_3$}
\put(3.7, 3.7){2}
\put(4.7, 3.7){2}
\put(5.7, 3.7){2}
\put(6.5, 3.7){$1_2 2$}
\put(7.6, 3.7){22}
\put(1.7, 4.7){$1_3$}
\put(2.7, 4.7){$1_3$}
\put(3.7, 4.7){2}
\put(4.6, 4.7){$1_4 2$}
\put(5.5, 4.7){$1_4 2$}
\put(6.6, 4.7){22}
\put(7.3, 4.7){$2 1_2 2$}
\put(1.7, 5.7){$1_3$}
\put(2.7, 5.7){$1_3$}
\put(3.7, 5.7){2}
\put(4.5, 5.7){$1_4 2$}
\put(5.2, 5.7){$1_3 1_4 2$}
\put(6.3, 5.7){$2 1_4 2$}
\put(7.4, 5.7){222}
\put(1.7, 6.7){$1_3$}
\put(2.4, 6.7){$1_1 1_3$}
\put(3.5, 6.7){$2 1_3$}
\put(4.6, 6.7){22}
\put(5.3, 6.7){$2 1_4 2$}
\put(6.2, 6.7){$2 1_3 1_4 2$}
\put(7.2, 6.7){$2 2 1_4 2$}
\end{picture}
\end{center}
\vspace{1in}

Notice that the elements along the rightmost column of the growth diagram represent a chain in $Z(k)$ as do the elements along the topmost row of the diagram.  In the following section, we give a method for representing a chain in $Z(k)$ as a $k$-ribbon Fibonacci path tableau and we define standard $k$-ribbon Fibonacci tableaux.

\section{$k$-Ribbon Fibonacci Tableaux}
\label{kribbon}

An element of $Z(k)$ can be represented by a {\it $k$-ribbon Fibonacci shape} by letting $1_j$ correspond to a group of $k$ squares consisting of $j$ adjacent squares in a single column followed by $k-j$ adjacent squares in the first row, and letting a $2$ correspond to a group of $2k$ squares consisting of a column of $k+1$ squares followed by $k-1$ adjacent squares in the first row.  For example, the element $1_4 1_1 2 2 1_5 2 1_2$ in $Z(5)$ is represented by

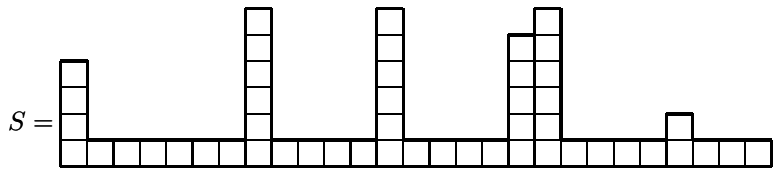
\setlength{\unitlength}{.7cm}
\begin{center}
\begin{picture}(14,3)(0,0)
\thicklines

\put(-.5,.7){$S=$}

\put(.5,0){\line(1,0){13.5}}
\put(.5,0){\line(0,1){2}}
\put(1,.5){\line(1,0){3}}
\put(.5,2){\line(1,0){.5}}

\put(1,.5){\line(0,1){1.5}}
\put(1.5,0){\line(0,1){.5}}
\put(4,0){\line(0,1){3}}

\put(4,3){\line(1,0){.5}}
\put(4.5,.5){\line(0,1){2.5}}
\put(4.5,.5){\line(1,0){2}}

\put(6.5,0){\line(0,1){3}}
\put(6.5,3){\line(1,0){.5}}
\put(7,.5){\line(0,1){2.5}}
\put(7,.5){\line(1,0){2}}

\put(9,0){\line(0,1){2.5}}
\put(9,2.5){\line(1,0){.5}}
\put(9.5,3){\line(1,0){.5}}
\put(9.5,0){\line(0,1){3}}
\put(10,.5){\line(0,1){2.5}}
\put(10,.5){\line(1,0){2}}

\put(12,0){\line(0,1){1}}
\put(12,1){\line(1,0){.5}}
\put(12.5,.5){\line(0,1){.5}}
\put(12.5,.5){\line(1,0){1.5}}
\put(14,0){\line(0,1){.5}}

\thinlines

\put(-.5,.7){$S=$}

\put(.5,0){\line(1,0){13.5}}
\put(.5,0){\line(0,1){2}}
\put(.5,.5){\line(1,0){13.5}}
\put(.5,1){\line(1,0){.5}}
\put(.5,1.5){\line(1,0){.5}}
\put(.5,2){\line(1,0){.5}}

\put(1,0){\line(0,1){2}}
\put(1.5,0){\line(0,1){.5}}
\put(2,0){\line(0,1){.5}}
\put(2.5,0){\line(0,1){.5}}
\put(3,0){\line(0,1){.5}}
\put(3.5,0){\line(0,1){.5}}
\put(4,0){\line(0,1){3}}

\put(4,1){\line(1,0){.5}}
\put(4,1.5){\line(1,0){.5}}
\put(4,2){\line(1,0){.5}}
\put(4,2.5){\line(1,0){.5}}
\put(4,3){\line(1,0){.5}}
\put(4.5,0){\line(0,1){3}}

\put(5,0){\line(0,1){.5}}
\put(5.5,0){\line(0,1){.5}}
\put(6,0){\line(0,1){.5}}
\put(6.5,0){\line(0,1){3}}

\put(6.5,1){\line(1,0){.5}}
\put(6.5,1.5){\line(1,0){.5}}
\put(6.5,2){\line(1,0){.5}}
\put(6.5,2.5){\line(1,0){.5}}
\put(6.5,3){\line(1,0){.5}}
\put(7,0){\line(0,1){3}}

\put(7.5,0){\line(0,1){.5}}
\put(8,0){\line(0,1){.5}}
\put(8.5,0){\line(0,1){.5}}
\put(9,0){\line(0,1){2.5}}

\put(9,1){\line(1,0){1}}
\put(9,1.5){\line(1,0){1}}
\put(9,2){\line(1,0){1}}
\put(9,2.5){\line(1,0){1}}
\put(9.5,3){\line(1,0){.5}}
\put(9.5,0){\line(0,1){3}}
\put(10,0){\line(0,1){3}}

\put(10.5,0){\line(0,1){.5}}
\put(11,0){\line(0,1){.5}}
\put(11.5,0){\line(0,1){.5}}
\put(12,0){\line(0,1){1}}
\put(12,1){\line(1,0){.5}}
\put(12.5,0){\line(0,1){1}}
\put(13,0){\line(0,1){.5}}
\put(13.5,0){\line(0,1){.5}}
\put(14,0){\line(0,1){.5}}

\end{picture}
\end{center}

We now describe two types of $k$-ribbons which can be used to cover any $k$-ribbon Fibonacci shape.  First, a {\em$k$-ribbon of height $j$} can be a grouping of $k$ squares so that there are $j$ adjacent squares stacked vertically in a single column starting in the first row followed directly by $k-j$ squares in the first row.  Notice any shape corresponding to a $1_j$ can be covered by a $k$-ribbon of height $j$.  Second, a {\em$k$-ribbon of height $j$} can be a grouping of $k$ squares which is made up of $j$ adjacent squares stacked vertically atop a column of height $k+1$ and $k-j$ squares in the rightmost positions of the adjacent $k-1$ squares in the first row.


For example, a $7$-ribbon of height 3 occupies the squares with dots in the example below while the squares without dots form a $7$-ribbon of height 5:

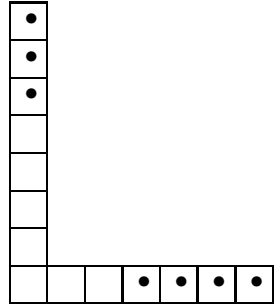
\setlength{\unitlength}{1cm}
\begin{center}
\begin{picture}(4,4)(0,0)
\thinlines

\put(0,0){\line(1,0){3.5}}
\put(0,0){\line(0,1){4}}
\put(0,.5){\line(1,0){3.5}}
\put(.5,0){\line(0,1){4}}

\put(0,1){\line(1,0){.5}}
\put(0,1.5){\line(1,0){.5}}
\put(0,2){\line(1,0){.5}}
\put(0,2.5){\line(1,0){.5}}
\put(0,3){\line(1,0){.5}}
\put(0,3.5){\line(1,0){.5}}
\put(0,4){\line(1,0){.5}}

\put(1,0){\line(0,1){.5}}
\put(1.5,0){\line(0,1){.5}}
\put(2,0){\line(0,1){.5}}
\put(2.5,0){\line(0,1){.5}}
\put(3,0){\line(0,1){.5}}
\put(3.5,0){\line(0,1){.5}}

\put(.2,2.7){$\bullet$}
\put(.2,3.2){$\bullet$}
\put(.2,3.7){$\bullet$}
\put(2.2,.2){$\bullet$}
\put(2.7,.2){$\bullet$}
\put(3.2,.2){$\bullet$}
\put(1.7,.2){$\bullet$}

\end{picture}
\end{center}

  Notice any $k$-ribbon shape corresponding to a $2$ can be covered by a $k$-ribbon of height $j$ stacked on top of a $k$-ribbon of height $k+1-j.$

 A {\it $k$-ribbon tiling} is a placement of $k$-ribbons into a domino Fibonacci shape such that all squares are covered.  All $k$-ribbon Fibonacci shapes have at least one valid $k$-ribbon tiling (by tiling squares created from a $1_j$ with a $k$-ribbon of height $j$ and squares created from a $2$ with a $k$-ribbon of height $j$, for some $j$, on top of a $k$-ribbon of height $k+1-j$).  For example, each of the following is a valid domino tiling of the shape corresponding to $1_4 1_1 2 2 1_5 2 1_2$:

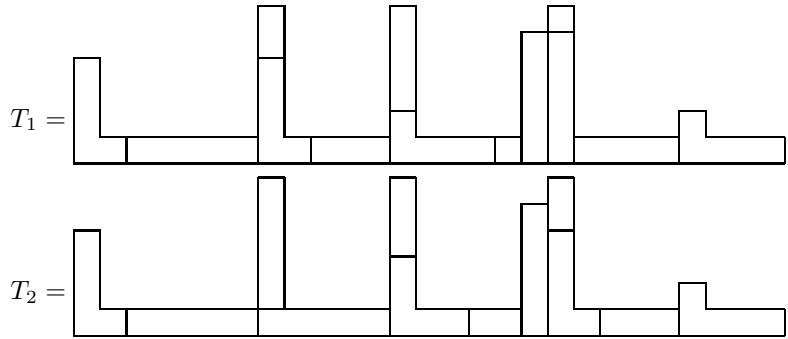
\setlength{\unitlength}{.7cm}
\begin{center}
\begin{picture}(14,3)(0,0)
\thinlines

\put(-.7,.7){$T_1=$}

\put(.5,0){\line(1,0){13.5}}
\put(.5,0){\line(0,1){2}}
\put(.5,2){\line(1,0){.5}}
\put(1,.5){\line(1,0){3}}

\put(1,.5){\line(0,1){1.5}}
\put(1.5,0){\line(0,1){.5}}
\put(4,0){\line(0,1){3}}

\put(4,2){\line(1,0){.5}}
\put(4,3){\line(1,0){.5}}
\put(4.5,.5){\line(0,1){2.5}}
\put(4.5,.5){\line(1,0){2}}

\put(5,0){\line(0,1){.5}}
\put(6.5,0){\line(0,1){3}}

\put(6.5,1){\line(1,0){.5}}
\put(6.5,3){\line(1,0){.5}}
\put(7,.5){\line(0,1){2.5}}
\put(7,.5){\line(1,0){2}}

\put(8.5,0){\line(0,1){.5}}
\put(9,0){\line(0,1){2.5}}

\put(9,2.5){\line(1,0){1}}
\put(9.5,3){\line(1,0){.5}}
\put(9.5,0){\line(0,1){3}}
\put(10,0){\line(0,1){3}}
\put(10,.5){\line(1,0){2}}

\put(12,0){\line(0,1){1}}
\put(12,1){\line(1,0){.5}}
\put(12.5,.5){\line(0,1){.5}}
\put(12.5,.5){\line(1,0){1.5}}
\put(14,0){\line(0,1){.5}}

\end{picture}
\end{center}

\setlength{\unitlength}{.7cm}
\begin{center}
\begin{picture}(14,3)(0,0)
\thinlines

\put(-.7,.7){$T_2=$}

\put(.5,0){\line(1,0){13.5}}
\put(.5,0){\line(0,1){2}}
\put(.5,2){\line(1,0){.5}}
\put(1,.5){\line(1,0){5.5}}

\put(1,.5){\line(0,1){1.5}}
\put(1.5,0){\line(0,1){.5}}
\put(4,0){\line(0,1){3}}

\put(4,3){\line(1,0){.5}}
\put(4.5,.5){\line(0,1){2.5}}
\put(4.5,.5){\line(1,0){2}}
\put(6.5,0){\line(0,1){3}}

\put(6.5,1.5){\line(1,0){.5}}
\put(6.5,3){\line(1,0){.5}}
\put(7,.5){\line(0,1){2.5}}
\put(7,.5){\line(1,0){2}}

\put(8,0){\line(0,1){.5}}
\put(9,0){\line(0,1){2.5}}

\put(9,2.5){\line(1,0){.5}}
\put(9.5,2){\line(1,0){.5}}
\put(9.5,3){\line(1,0){.5}}
\put(9.5,0){\line(0,1){3}}
\put(10,.5){\line(0,1){2.5}}
\put(10.5,0){\line(0,1){.5}}
\put(10,.5){\line(1,0){2}}

\put(12,0){\line(0,1){1}}
\put(12,1){\line(1,0){.5}}
\put(12.5,.5){\line(0,1){.5}}
\put(12.5,.5){\line(1,0){1.5}}
\put(14,0){\line(0,1){.5}}

\end{picture}
\end{center}

>From this point on, we shall use the term ``column" only when referring to the set of squares occupied by the Fibonacci shape corresponding to a $1_j$ or a $2.$  A ``column" corresponding to a $1_j$ will be called a ``column of height 1" and a column corresponding to a $2$ will be called a ``column of height 2".  
For instance, the tilings pictured above each have seven columns.  Notice that given a $k$-ribbon shape, the number of columns in any $k$-ribbon tiling of that shape will be the same.

Now we define ${\mathcal{F}_k}$ to be the poset of $k$-ribbon Fibonacci shapes together with cover relations inherited from $Z(k)$.  ${\mathcal{F}_k}$ is naturally isomorphic to $Z(k)$.  A chain $\nu = (\emptyset, \nu_1, \nu_2, \cdots, \nu_k = \nu)$ in $Z(k)$ can be translated into a {\it $k$-ribbon Fibonacci path tableau} by placing $i$'s in $\nu_i/\nu_{i-1}$, i.e. in each of the $k$ new squares created at the $i$th step.  For example, the chain
\[
\nu = (\emptyset, \ 1_2, \ 1_5 1_2, \ 1_2 1_5 1_2, \ 2 1_5 1_2, \ 2 2 1_2, \ 2 1_5 2 1_2, \ 1_4 2 1_5 2 1_2, \ 2 2 1_5 2 1_2, \ 1_1 2 2 1_5 2 1_2)
\]

in $Z(5)$ corresponds to the following $5$-ribbon Fibonacci path tableau:

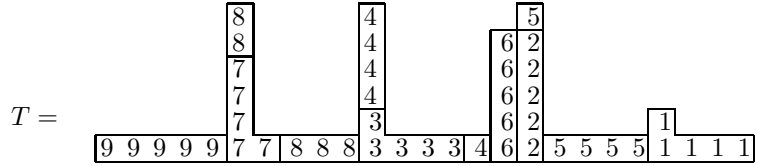
\setlength{\unitlength}{.7cm}
\begin{center}
\begin{picture}(14,3)(0,0)
\thinlines

\put(-.1,.7){$T=$}

\put(1.5,0){\line(1,0){12.5}}
\put(1.5,.5){\line(1,0){2.5}}
\put(1.5,0){\line(0,1){.5}}
\put(4,0){\line(0,1){3}}

\put(4,2){\line(1,0){.5}}
\put(4,3){\line(1,0){.5}}
\put(4.5,.5){\line(0,1){2.5}}
\put(4.5,.5){\line(1,0){2}}

\put(5,0){\line(0,1){.5}}
\put(6.5,0){\line(0,1){3}}

\put(6.5,1){\line(1,0){.5}}
\put(6.5,3){\line(1,0){.5}}
\put(7,.5){\line(0,1){2.5}}
\put(7,.5){\line(1,0){2}}

\put(8.5,0){\line(0,1){.5}}
\put(9,0){\line(0,1){2.5}}

\put(9,2.5){\line(1,0){1}}
\put(9.5,3){\line(1,0){.5}}
\put(9.5,0){\line(0,1){3}}
\put(10,0){\line(0,1){3}}
\put(10,.5){\line(1,0){2}}

\put(12,0){\line(0,1){1}}
\put(12,1){\line(1,0){.5}}
\put(12.5,.5){\line(0,1){.5}}
\put(12.5,.5){\line(1,0){1.5}}
\put(14,0){\line(0,1){.5}}

\put(12.2, .1){$1$}
\put(12.2, .6){$1$}
\put(12.7,.1){$1$}
\put(13.2,.1){$1$}
\put(13.7,.1){$1$}

\put(9.7,.1){$2$}
\put(9.7,.6){$2$}
\put(9.7,1.1){$2$}
\put(9.7,1.6){$2$}
\put(9.7,2.1){$2$}

\put(6.7,.1){$3$}
\put(6.7,.6){$3$}
\put(7.2,.1){$3$}
\put(7.7,.1){$3$}
\put(8.2,.1){$3$}

\put(6.6,1.1){$4$}
\put(6.6,1.6){$4$}
\put(6.6,2.1){$4$}
\put(6.6,2.6){$4$}
\put(8.7,.1){$4$}

\put(9.7,2.6){$5$}
\put(10.2,.1){$5$}
\put(10.7,.1){$5$}
\put(11.2,.1){$5$}
\put(11.7,.1){$5$}

\put(9.2,.1){$6$}
\put(9.2,.6){$6$}
\put(9.2,1.1){$6$}
\put(9.2,1.6){$6$}
\put(9.2,2.1){$6$}

\put(4.1,.1){$7$}
\put(4.1,.6){$7$}
\put(4.1,1.1){$7$}
\put(4.1,1.6){$7$}
\put(4.6,.1){$7$}

\put(4.1,2.1){$8$}
\put(4.1,2.6){$8$}
\put(5.2,.1){$8$}
\put(5.7,.1){$8$}
\put(6.2,.1){$8$}

\put(1.6,.1){$9$}
\put(2.1,.1){$9$}
\put(2.6,.1){$9$}
\put(3.1,.1){$9$}
\put(3.6,.1){$9$}

\end{picture}
\end{center}

As seen in Section \ref{background}, Fomin's method gives a bijection between $k$-colored permutations and pairs of chains in $Z(k)$, each of which can now be represented by a $k$-ribbon Fibonacci path tableau.  We will call the $k$-ribbon Fibonacci path tableau obtained from the right edge of the diagram $\hat{P}$ and the one obtained from the top edge of the diagram $\hat{Q}$.  From our previous growth diagram:

\setlength{\unitlength}{1.5cm}
\begin{center}
\begin{picture}(7,6.5)(1,1)
\thinlines
\put(1,0){\line(0,1){7}}
\put(2,0){\line(0,1){7}}
\put(3,0){\line(0,1){7}}
\put(4,0){\line(0,1){7}}
\put(5,0){\line(0,1){7}}
\put(6,0){\line(0,1){7}}
\put(7,0){\line(0,1){7}}
\put(8,0){\line(0,1){7}}
\put(1,0){\line(1,0){7}}
\put(1,1){\line(1,0){7}}
\put(1,2){\line(1,0){7}}
\put(1,3){\line(1,0){7}}
\put(1,4){\line(1,0){7}}
\put(1,5){\line(1,0){7}}
\put(1,6){\line(1,0){7}}
\put(1,7){\line(1,0){7}}
\put(.8,-.3){$\emptyset$}
\put(.8,.7){$\emptyset$}
\put(.8,1.7){$\emptyset$}
\put(.8,2.7){$\emptyset$}
\put(.8,3.7){$\emptyset$}
\put(.8,4.7){$\emptyset$}
\put(.8,5.7){$\emptyset$}
\put(.8,6.7){$\emptyset$}
\put(1.8,-.3){$\emptyset$}
\put(2.8,-.3){$\emptyset$}
\put(3.8,-.3){$\emptyset$}
\put(4.8,-.3){$\emptyset$}
\put(5.8,-.3){$\emptyset$}
\put(6.8,-.3){$\emptyset$}
\put(7.8,-.3){$\emptyset$}
\put(1.4, 1.4){$X^3$}
\put(2.4, 6.3){$X^1$}
\put(3.4, .4){$X^1$}
\put(4.4, 4.4){$X^4$}
\put(5.3, 5.3){$X^3$}
\put(6.3, 3.3){$X^2$}
\put(7.4, 2.4){$X^4$}
\put(1.8, .7){$\emptyset$}
\put(2.8, .7){$\emptyset$}
\put(3.7, .7){$1_1$}
\put(4.7, .7){$1_1$}
\put(5.7, .7){$1_1$}
\put(6.7, .7){$1_1$}
\put(7.7, .7){$1_1$}
\put(1.7, 1.7){$1_3$}
\put(2.7, 1.7){$1_3$}
\put(3.7, 1.7){2}
\put(4.7, 1.7){2}
\put(5.7, 1.7){2}
\put(6.7, 1.7){2}
\put(7.7, 1.7){2}
\put(1.7, 2.7){$1_3$}
\put(2.7, 2.7){$1_3$}
\put(3.7, 2.7){2}
\put(4.7, 2.7){2}
\put(5.7, 2.7){2}
\put(6.7, 2.7){2}
\put(7.6, 2.7){$1_42$}
\put(1.7, 3.7){$1_3$}
\put(2.7, 3.7){$1_3$}
\put(3.7, 3.7){2}
\put(4.7, 3.7){2}
\put(5.7, 3.7){2}
\put(6.5, 3.7){$1_2 2$}
\put(7.6, 3.7){22}
\put(1.7, 4.7){$1_3$}
\put(2.7, 4.7){$1_3$}
\put(3.7, 4.7){2}
\put(4.6, 4.7){$1_4 2$}
\put(5.5, 4.7){$1_4 2$}
\put(6.6, 4.7){22}
\put(7.3, 4.7){$2 1_2 2$}
\put(1.7, 5.7){$1_3$}
\put(2.7, 5.7){$1_3$}
\put(3.7, 5.7){2}
\put(4.5, 5.7){$1_4 2$}
\put(5.2, 5.7){$1_3 1_4 2$}
\put(6.3, 5.7){$2 1_4 2$}
\put(7.4, 5.7){222}
\put(1.7, 6.7){$1_3$}
\put(2.4, 6.7){$1_1 1_3$}
\put(3.5, 6.7){$2 1_3$}
\put(4.6, 6.7){22}
\put(5.3, 6.7){$2 1_4 2$}
\put(6.2, 6.7){$2 1_3 1_4 2$}
\put(7.2, 6.7){$2 2 1_4 2$}

\put(4.5,7.3){$\hat{Q}$}
\put(8.3,3.5){$\hat{P}$}
\end{picture}
\end{center}
\vspace{1in}

we have

\setlength{\unitlength}{.7cm}
\begin{center}
\begin{picture}(9.5,3)(0,0)
\thinlines

\put(-.7,1.2){$\hat{P} = $}

\put(.5,0){\line(1,0){8.5}}
\put(.5,0){\line(0,1){3}}
\put(.5,2){\line(1,0){.5}}
\put(.5,3){\line(1,0){.5}}

\put(1,.5){\line(0,1){2.5}}
\put(1,.5){\line(1,0){2}}
\put(1.5,0){\line(0,1){.5}}
\put(3,0){\line(0,1){3}}

\put(3,1){\line(1,0){.5}}
\put(3,3){\line(1,0){.5}}
\put(3.5,.5){\line(0,1){2.5}}
\put(3.5,.5){\line(1,0){2}}
\put(5,0){\line(0,1){.5}}

\put(5.5,0){\line(0,1){2}}
\put(5.5,2){\line(1,0){.5}}
\put(6,.5){\line(0,1){1.5}}
\put(6,.5){\line(1,0){3}}

\put(6.5,0){\line(0,1){3}}
\put(6.5,3){\line(1,0){.5}}
\put(7,.5){\line(0,1){2.5}}
\put(9,0){\line(0,1){.5}}

\put(.6,.1){$3$}
\put(.6,.6){$3$}
\put(.6,1.1){$3$}
\put(.6,1.6){$3$}
\put(.6,2.1){$4$}
\put(.6,2.6){$4$}

\put(1.1,.1){$3$}
\put(1.6,.1){$4$}
\put(2.1,.1){$4$}
\put(2.6,.1){$4$}

\put(3.1,.1){$5$}
\put(3.1,.6){$5$}
\put(3.1,1.1){$6$}
\put(3.1,1.6){$6$}
\put(3.1,2.1){$6$}
\put(3.1,2.6){$6$}

\put(3.6,.1){$5$}
\put(4.1,.1){$5$}
\put(4.6,.1){$5$}
\put(5.1,.1){$6$}

\put(5.7,.1){$7$}
\put(5.7,.6){$7$}
\put(5.7,1.1){$7$}
\put(5.7,1.6){$7$}
\put(6.2,.1){$7$}

\put(6.6,.1){$1$}
\put(6.6,.6){$2$}
\put(6.6,1.1){$2$}
\put(6.6,1.6){$2$}
\put(6.6,2.1){$2$}
\put(6.6,2.6){$2$}

\put(7.1,.1){$1$}
\put(7.6,.1){$1$}
\put(8.1,.1){$1$}
\put(8.6,.1){$1$}

\end{picture}
\end{center}

\setlength{\unitlength}{.7cm}
\begin{center}
\begin{picture}(9.5,3)(0,0)
\thinlines

\put(-.7,1.2){$\hat{Q} = $}

\put(.5,0){\line(1,0){8.5}}
\put(.5,0){\line(0,1){3}}
\put(.5,3){\line(1,0){.5}}

\put(1,.5){\line(0,1){2.5}}
\put(.5,.5){\line(1,0){2.5}}
\put(3,0){\line(0,1){3}}

\put(3,1.5){\line(1,0){.5}}
\put(3,3){\line(1,0){.5}}
\put(3.5,.5){\line(0,1){2.5}}
\put(3.5,.5){\line(1,0){2}}
\put(4.5,0){\line(0,1){.5}}

\put(5.5,0){\line(0,1){2}}
\put(5.5,2){\line(1,0){.5}}
\put(6,.5){\line(0,1){1.5}}
\put(6,.5){\line(1,0){.5}}

\put(6.5,0){\line(0,1){3}}
\put(6.5,1.5){\line(1,0){.5}}
\put(6.5,3){\line(1,0){.5}}
\put(7,.5){\line(0,1){2.5}}
\put(7,.5){\line(1,0){2}}
\put(8,0){\line(0,1){.5}}
\put(9,0){\line(0,1){.5}}

\put(.6,.1){$2$}
\put(.6,.6){$3$}
\put(.6,1.1){$3$}
\put(.6,1.6){$3$}
\put(.6,2.1){$3$}
\put(.6,2.6){$3$}

\put(1.1,.1){$2$}
\put(1.6,.1){$2$}
\put(2.1,.1){$2$}
\put(2.6,.1){$2$}

\put(3.1,.1){$6$}
\put(3.1,.6){$6$}
\put(3.1,1.1){$6$}
\put(3.1,1.6){$7$}
\put(3.1,2.1){$7$}
\put(3.1,2.6){$7$}

\put(3.6,.1){$6$}
\put(4.1,.1){$6$}
\put(4.6,.1){$7$}
\put(5.1,.1){$7$}

\put(5.7,.1){$5$}
\put(5.7,.6){$5$}
\put(5.7,1.1){$5$}
\put(5.7,1.6){$5$}
\put(6.2,.1){$5$}

\put(6.6,.1){$1$}
\put(6.6,.6){$1$}
\put(6.6,1.1){$1$}
\put(6.6,1.6){$4$}
\put(6.6,2.1){$4$}
\put(6.6,2.6){$4$}

\put(7.1,.1){$1$}
\put(7.6,.1){$1$}
\put(8.1,.1){$4$}
\put(8.6,.1){$4$}

\end{picture}
\end{center}

We now define a {\it $k$-ribbon Fibonacci tableau} to be a filling of the $k$-ribbons in a $k$-ribbon Fibonacci shape with $k$ 1's, $k$ 2's, $k$ 3's, $\dots$, $k$ $n$'s such that each number appears in exactly one $k$-ribbon and each $k$-ribbon contains $k$ of the same number.  A {\it standard} $k$-ribbon Fibonacci tableau is a $k$-ribbon Fibonacci tableau with the additional properties that (i) the $k$-ribbon containing the leftmost square in the bottom row is the $k$-ribbon containing $n,$ and (ii) for every $j$, the $k$-ribbon containing $j$ is either appended as a $k$-ribbon to the shape containing $i$'s for $j<i\leq n$ or is stacked as a $k$-ribbon on top of a single $k$-ribbon containing $i$ for $j<i\leq n.$   For example, the following is a standard $k$-ribbon Fibonacci tableau:

\setlength{\unitlength}{.7cm}
\begin{center}
\begin{picture}(14,3)(0,0)
\thinlines

\put(-.1,.7){$T=$}

\put(1.5,0){\line(1,0){12.5}}
\put(1.5,.5){\line(1,0){2.5}}
\put(1.5,0){\line(0,1){.5}}
\put(4,0){\line(0,1){3}}

\put(4,2){\line(1,0){.5}}
\put(4,3){\line(1,0){.5}}
\put(4.5,.5){\line(0,1){2.5}}
\put(4.5,.5){\line(1,0){2}}

\put(5,0){\line(0,1){.5}}
\put(6.5,0){\line(0,1){3}}

\put(6.5,1){\line(1,0){.5}}
\put(6.5,3){\line(1,0){.5}}
\put(7,.5){\line(0,1){2.5}}
\put(7,.5){\line(1,0){2}}

\put(8.5,0){\line(0,1){.5}}
\put(9,0){\line(0,1){2.5}}

\put(9,2.5){\line(1,0){1}}
\put(9.5,3){\line(1,0){.5}}
\put(9.5,0){\line(0,1){3}}
\put(10,0){\line(0,1){3}}
\put(10,.5){\line(1,0){2}}

\put(12,0){\line(0,1){1}}
\put(12,1){\line(1,0){.5}}
\put(12.5,.5){\line(0,1){.5}}
\put(12.5,.5){\line(1,0){1.5}}
\put(14,0){\line(0,1){.5}}

\put(12.2, .1){$3$}
\put(12.2, .6){$3$}
\put(12.7,.1){$3$}
\put(13.2,.1){$3$}
\put(13.7,.1){$3$}

\put(9.7,.1){$4$}
\put(9.7,.6){$4$}
\put(9.7,1.1){$4$}
\put(9.7,1.6){$4$}
\put(9.7,2.1){$4$}

\put(6.7,.1){$7$}
\put(6.7,.6){$7$}
\put(7.2,.1){$7$}
\put(7.7,.1){$7$}
\put(8.2,.1){$7$}

\put(6.6,1.1){$2$}
\put(6.6,1.6){$2$}
\put(6.6,2.1){$2$}
\put(6.6,2.6){$2$}
\put(8.7,.1){$2$}

\put(9.7,2.6){$1$}
\put(10.2,.1){$1$}
\put(10.7,.1){$1$}
\put(11.2,.1){$1$}
\put(11.7,.1){$1$}

\put(9.2,.1){$5$}
\put(9.2,.6){$5$}
\put(9.2,1.1){$5$}
\put(9.2,1.6){$5$}
\put(9.2,2.1){$5$}

\put(4.1,.1){$8$}
\put(4.1,.6){$8$}
\put(4.1,1.1){$8$}
\put(4.1,1.6){$8$}
\put(4.6,.1){$8$}

\put(4.1,2.1){$6$}
\put(4.1,2.6){$6$}
\put(5.2,.1){$6$}
\put(5.7,.1){$6$}
\put(6.2,.1){$6$}

\put(1.6,.1){$9$}
\put(2.1,.1){$9$}
\put(2.6,.1){$9$}
\put(3.1,.1){$9$}
\put(3.6,.1){$9$}

\end{picture}
\end{center}

As is the case for $k$-ribbon Fibonacci path tableaux, it is possible to think of a standard $k$-ribbon Fibonacci tableau as a chain in a poset.  Define $S(k)$ to be a new partial order on the set of words in the alphabet $\{ 1_1, 1_2,\ldots,1_k, 2\}$ in which an element $z$ is covered by an element $w$ if $w$ is obtained from $z$ by appending a $1_i$ for some $i=1,\ldots,k$ or if $w$ is obtained from $z$ by replacing a $1_i$ by a $2$.  Then any standard $k$-ribbon Fibonacci tableau of shape $w$ corresponds to the path tableau represented by a maximal chain from $\emptyset$ to $w$ in $S(k)$, with $j$'s placed in the $k$-ribbon created at the $(n-j+1)$st step.  

The evacuation algorithm presented in section 5 can be used to show that the number of standard $k$-ribbon Fibonacci tableaux is equal to the number of $k$-ribbon Fibonacci path tableaux.

\section{A $k$-Ribbon Fibonacci Insertion}
\label{insertion}

For a permutation in $S_n$, Fomin's growth diagrams can be used to obtain a pair of chains in Young's lattice which can be represented as a pair of standard Young tableaux, $(\tilde{P}, \tilde{Q})$, of the same shape $\lambda$.  The Schensted correspondence is an insertion algorithm which takes any permutation to a pair of standard Young tableaux, $(P, Q)$, of the same shape $\lambda$.  For Young's lattice, the pair of standard Young tableaux given through these two methods are the same. 

In addition to Young's lattice, Fomin's growth diagrams can be used to give a bijection between a ($k$-colored) permutation and a pair of chains in the Fibonacci poset $Z(k)$ which can be represented as a pair of $k$-ribbon Fibonacci path tableaux $(\hat{P}, \hat{Q})$.  In the case of $Z(1)$ and $Z(2)$, there is an insertion algorithm which gives a bijection between a 1- or 2-colored permutation and a pair of tableaux $(P, Q)$ of the same shape where $P$ is a standard (domino) Fibonacci tableau and $Q$ is a (domino) Fibonacci path tableau.  Unlike Young's lattice, the two pair of tableaux obtained from these two methods are not the same.  While $\hat{Q} = Q$, $\hat{P}$ is not equal to $P$.  In each case, however, there is an evacuation map, $ev$, for Fibonacci tableaux such that $ev(P) = \hat{P}$.  See \cite{CamKill}, \cite{Kil} for details.

We now describe a $k$-ribbon insertion algorithm which takes a $k$-colored permutation to a pair $(P,Q)$ of $k$-ribbon Fibonacci tableaux.  The $P$ tableau that is created will be a standard $k$-ribbon Fibonacci tableau and the $Q$ tableau that is created will be a $k$-ribbon Fibonacci path tableau.  To apply our algorithm to a colored permutation, we will construct a sequence $\{ (P_i, Q_i) \}_{i=0}^n$ where $(P_0, Q_0) = (\emptyset, \emptyset)$ and $(P_i, Q_i)$ are the tableaux obtained from the insertion of $x_i^j$ into $P_{i-1}$ in the following manner.  The general principle of the insertion algorithm is that each element $x_i^j$ of the permutation (with value $x_i$ and color $j$) is inserted into $P_{i-1}$ as a $k$-ribbon of height $j$.  Recall that $i$ ranges between $1$ and $n$ and $j$ ranges between $1$ and $k$.

\begin{enumerate}

\item  Compare the value of $x_i$ to the value $t_1$ in the $k$-ribbon containing the leftmost square in the bottom row of $P_{i-1}$.  

\item  If $x_i > t_1$, add a $k$-ribbon of height $j$ containing $x_i$'s to the left of the square containing $t_1$ in the bottom row.  Call this new tableau $P_i$.  To form $Q_i$, a tableau of the same shape as $P_i$, place $i$'s in this newly created $k$-ribbon of height $j$.

\item  If $x_i < t_1$ then place a $k$-ribbon of height $j$ on top of the $k$-ribbon containing $t_1$, while forcing the $k$-ribbon containing $t_1$ to become a $k$-ribbon of height $k+1-j.$

If there were no $k$-ribbon on top of the $k$-ribbon containing $t_1$ in $P_{i-1}$, then this new tableau is $P_i$.  Obtain $Q_i$ by placing $i$'s into the $k$ new squares that were created in the first stack of 2 $k$-ribbons.

If there were a $k$-ribbon of height $l$ containing $b$'s on top of the $k$-ribbon containing $t_1$ in $P_{i-1}$, then the $k$-ribbon of height $l$ containing $b$'s is bumped out of the first stack.  Continue by inductively inserting the $k$-ribbon of height $l$ containing $b$'s into the tableau to the right of the first stack by comparing $b$ to the element $t_2$ in the $k$-ribbon in the bottom row of the next stack and repeating the insertion algorithm.

\end{enumerate}

\begin{Example}  When applying the insertion algorithm to the $5$-colored permutation $\pi = 2^3 7^1 1^1 5^4 6^3 4^2 3^4$ that was used to form the square diagram in Section 2, we obtain the following:
\end{Example}

\setlength{\unitlength}{.7cm}
\begin{center}
\begin{picture}(10,3)(0,0)
\thinlines

\put(0,.5){$P_i$ :}

\put(1.5,0){\line(1,0){1.5}}
\put(1.5,0){\line(0,1){1.5}}
\put(1.5,1.5){\line(1,0){.5}}
\put(2,.5){\line(1,0){1}}
\put(2,.5){\line(0,1){1}}
\put(3,0){\line(0,1){.5}}
\put(3.2,0){,}

\put(1.6,.1){$2$}
\put(1.6,.6){$2$}
\put(1.6,1.1){$2$}
\put(2.1,.1){$2$}
\put(2.6,.1){$2$}

\put(4,0){\line(1,0){4}}
\put(4,0){\line(0,1){.5}}
\put(4,.5){\line(1,0){2.5}}
\put(6.5,0){\line(0,1){1.5}}
\put(6.5,1.5){\line(1,0){.5}}
\put(7,.5){\line(0,1){1}}
\put(7,.5){\line(1,0){1}}
\put(8,0){\line(0,1){.5}}
\put(8.2,0){,}

\put(4.1,.1){$7$}
\put(4.6,.1){$7$}
\put(5.1,.1){$7$}
\put(5.6,.1){$7$}
\put(6.1,.1){$7$}
\put(6.6,.1){$2$}
\put(6.6,.6){$2$}
\put(6.6,1.1){$2$}
\put(7.1,.1){$2$}
\put(7.6,.1){$2$}

\end{picture}
\end{center}

\setlength{\unitlength}{.7cm}
\begin{center}
\begin{picture}(11,3.5)(0,0)
\thinlines

\put(0,0){\line(1,0){4}}
\put(0,0){\line(0,1){3}}
\put(0,2.5){\line(1,0){.5}}
\put(0,3){\line(1,0){.5}}
\put(.5,0){\line(0,1){3}}
\put(.5,.5){\line(1,0){2}}
\put(2.5,0){\line(0,1){1.5}}
\put(2.5,1.5){\line(1,0){.5}}
\put(3,.5){\line(1,0){1}}
\put(3,.5){\line(0,1){1}}
\put(4,0){\line(0,1){.5}}
\put(4.2,0){,}

\put(.1,.1){$7$}
\put(.1,.6){$7$}
\put(.1,1.1){$7$}
\put(.1,1.6){$7$}
\put(.1,2.1){$7$}
\put(.1,2.6){$1$}
\put(.6,.1){$1$}
\put(1.1,.1){$1$}
\put(1.6,.1){$1$}
\put(2.1,.1){$1$}
\put(2.6,.1){$2$}
\put(2.6,.6){$2$}
\put(2.6,1.1){$2$}
\put(3.1,.1){$2$}
\put(3.6,.1){$2$}

\put(5.5,0){\line(1,0){5}}
\put(5.5,0){\line(0,1){3}}
\put(5.5,1){\line(1,0){.5}}
\put(5.5,3){\line(1,0){.5}}
\put(6,.5){\line(1,0){2}}
\put(6,.5){\line(0,1){2.5}}
\put(7.5,0){\line(0,1){.5}}
\put(8,0){\line(0,1){3}}
\put(8,2.5){\line(1,0){.5}}
\put(8,3){\line(1,0){.5}}
\put(8.5,0){\line(0,1){3}}
\put(8.5,.5){\line(1,0){2}}
\put(10.5,0){\line(0,1){.5}}
\put(10.7,0){,}

\put(5.6,.1){$7$}
\put(5.6,.6){$7$}
\put(5.6,1.1){$5$}
\put(5.6,1.6){$5$}
\put(5.6,2.1){$5$}
\put(5.6,2.6){$5$}
\put(6.1,.1){$7$}
\put(6.6,.1){$7$}
\put(7.1,.1){$7$}
\put(7.6,.1){$5$}
\put(8.1,.1){$2$}
\put(8.1,.6){$2$}
\put(8.1,1.1){$2$}
\put(8.1,1.6){$2$}
\put(8.1,2.1){$2$}
\put(8.1,2.6){$1$}
\put(8.6,.1){$1$}
\put(9.1,.1){$1$}
\put(9.6,.1){$1$}
\put(10.1,.1){$1$}

\end{picture}
\end{center}

\setlength{\unitlength}{.7cm}
\begin{center}
\begin{picture}(15,3.5)(0,0)
\thinlines

\put(0,0){\line(1,0){6}}
\put(0,0){\line(0,1){3}}
\put(0,1.5){\line(1,0){.5}}
\put(0,3){\line(1,0){.5}}
\put(.5,.5){\line(0,1){2.5}}
\put(.5,.5){\line(1,0){2}}
\put(1.5,0){\line(0,1){.5}}
\put(2.5,0){\line(0,1){2}}
\put(2.5,2){\line(1,0){.5}}
\put(3,.5){\line(0,1){1.5}}
\put(3,.5){\line(1,0){.5}}
\put(3.5,0){\line(0,1){3}}
\put(3.5,2.5){\line(1,0){.5}}
\put(3.5,3){\line(1,0){.5}}
\put(4,0){\line(0,1){3}}
\put(4,.5){\line(1,0){2}}
\put(6,0){\line(0,1){.5}}
\put(6.2,0){,}

\put(.1,.1){$7$}
\put(.1,.6){$7$}
\put(.1,1.1){$7$}
\put(.1,1.6){$6$}
\put(.1,2.1){$6$}
\put(.1,2.6){$6$}
\put(.6,.1){$7$}
\put(1.1,.1){$7$}
\put(1.6,.1){$6$}
\put(2.1,.1){$6$}
\put(2.6,.1){$5$}
\put(2.6,.6){$5$}
\put(2.6,1.1){$5$}
\put(2.6,1.6){$5$}
\put(3.1,.1){$5$}
\put(3.6,.1){$2$}
\put(3.6,.6){$2$}
\put(3.6,1.1){$2$}
\put(3.6,1.6){$2$}
\put(3.6,2.1){$2$}
\put(3.6,2.6){$1$}
\put(4.1,.1){$1$}
\put(4.6,.1){$1$}
\put(5.1,.1){$1$}
\put(5.6,.1){$1$}

\put(7,0){\line(1,0){7.5}}
\put(7,0){\line(0,1){3}}
\put(7.5,.5){\line(1,0){2}}
\put(7,2){\line(1,0){.5}}
\put(7,3){\line(1,0){.5}}
\put(7.5,.5){\line(0,1){2.5}}
\put(8,0){\line(0,1){.5}}
\put(9.5,0){\line(0,1){1.5}}
\put(9.5,1.5){\line(1,0){.5}}
\put(10,.5){\line(0,1){1}}
\put(10,.5){\line(1,0){1}}
\put(11,0){\line(0,1){2}}
\put(11,2){\line(1,0){.5}}
\put(11.5,.5){\line(0,1){1.5}}
\put(11.5,.5){\line(1,0){.5}}
\put(12,0){\line(0,1){3}}
\put(12,2.5){\line(1,0){.5}}
\put(12,3){\line(1,0){.5}}
\put(12.5,0){\line(0,1){3}}
\put(12.5,.5){\line(1,0){2}}
\put(14.5,0){\line(0,1){.5}}
\put(14.7,0){,}

\put(7.1,.1){$7$}
\put(7.1,.6){$7$}
\put(7.1,1.1){$7$}
\put(7.1,1.6){$7$}
\put(7.1,2.1){$4$}
\put(7.1,2.6){$4$}
\put(7.6,.1){$7$}
\put(8.1,.1){$4$}
\put(8.6,.1){$4$}
\put(9.1,.1){$4$}
\put(9.6,.1){$6$}
\put(9.6,.6){$6$}
\put(9.6,1.1){$6$}
\put(10.1,.1){$6$}
\put(10.6,.1){$6$}
\put(11.1,.1){$5$}
\put(11.1,.6){$5$}
\put(11.1,1.1){$5$}
\put(11.1,1.6){$5$}
\put(11.6,.1){$5$}
\put(12.1,.1){$2$}
\put(12.1,.6){$2$}
\put(12.1,1.1){$2$}
\put(12.1,1.6){$2$}
\put(12.1,2.1){$2$}
\put(12.1,2.6){$1$}
\put(12.6,.1){$1$}
\put(13.1,.1){$1$}
\put(13.6,.1){$1$}
\put(14.1,.1){$1$}

\end{picture}
\end{center}

\setlength{\unitlength}{.7cm}
\begin{center}
\begin{picture}(9,3.5)(0,0)
\thinlines

\put(-1.3,1){$P=$}
\put(0,0){\line(1,0){8.5}}
\put(0,0){\line(0,1){3}}
\put(.5,.5){\line(1,0){2}}
\put(0,1){\line(1,0){.5}}
\put(0,3){\line(1,0){.5}}
\put(.5,.5){\line(0,1){2.5}}
\put(2,0){\line(0,1){.5}}
\put(2.5,0){\line(0,1){3}}
\put(2.5,2){\line(1,0){.5}}
\put(2.5,3){\line(1,0){.5}}
\put(3,.5){\line(0,1){2.5}}
\put(3,.5){\line(1,0){2}}
\put(3.5,0){\line(0,1){.5}}
\put(5,0){\line(0,1){2}}
\put(5,2){\line(1,0){.5}}
\put(5.5,.5){\line(0,1){1.5}}
\put(5.5,.5){\line(1,0){.5}}
\put(6,0){\line(0,1){3}}
\put(6,2.5){\line(1,0){.5}}
\put(6,3){\line(1,0){.5}}
\put(6.5,0){\line(0,1){3}}
\put(6.5,.5){\line(1,0){2}}
\put(8.5,0){\line(0,1){.5}}

\put(.1,.1){$7$}
\put(.1,.6){$7$}
\put(.1,1.1){$3$}
\put(.1,1.6){$3$}
\put(.1,2.1){$3$}
\put(.1,2.6){$3$}
\put(.6,.1){$7$}
\put(1.1,.1){$7$}
\put(1.6,.1){$7$}
\put(2.1,.1){$3$}
\put(2.6,.1){$6$}
\put(2.6,.6){$6$}
\put(2.6,1.1){$6$}
\put(2.6,1.6){$6$}
\put(2.6,2.1){$4$}
\put(2.6,2.6){$4$}
\put(3.1,.1){$6$}
\put(3.6,.1){$4$}
\put(4.1,.1){$4$}
\put(4.6,.1){$4$}
\put(5.1,.1){$5$}
\put(5.1,.6){$5$}
\put(5.1,1.1){$5$}
\put(5.1,1.6){$5$}
\put(5.6,.1){$5$}
\put(6.1,.1){$2$}
\put(6.1,.6){$2$}
\put(6.1,1.1){$2$}
\put(6.1,1.6){$2$}
\put(6.1,2.1){$2$}
\put(6.1,2.6){$1$}
\put(6.6,.1){$1$}
\put(7.1,.1){$1$}
\put(7.6,.1){$1$}
\put(8.1,.1){$1$}

\end{picture}
\end{center}

\setlength{\unitlength}{.7cm}
\begin{center}
\begin{picture}(10,3)(0,0)
\thinlines

\put(0,.5){$Q_i$ :}

\put(1.5,0){\line(1,0){1.5}}
\put(1.5,0){\line(0,1){1.5}}
\put(1.5,1.5){\line(1,0){.5}}
\put(2,.5){\line(1,0){1}}
\put(2,.5){\line(0,1){1}}
\put(3,0){\line(0,1){.5}}
\put(3.2,0){,}

\put(1.6,.1){$1$}
\put(1.6,.6){$1$}
\put(1.6,1.1){$1$}
\put(2.1,.1){$1$}
\put(2.6,.1){$1$}

\put(4,0){\line(1,0){4}}
\put(4,0){\line(0,1){.5}}
\put(4,.5){\line(1,0){2.5}}
\put(6.5,0){\line(0,1){1.5}}
\put(6.5,1.5){\line(1,0){.5}}
\put(7,.5){\line(0,1){1}}
\put(7,.5){\line(1,0){1}}
\put(8,0){\line(0,1){.5}}
\put(8.2,0){,}

\put(4.1,.1){$2$}
\put(4.6,.1){$2$}
\put(5.1,.1){$2$}
\put(5.6,.1){$2$}
\put(6.1,.1){$2$}
\put(6.6,.1){$1$}
\put(6.6,.6){$1$}
\put(6.6,1.1){$1$}
\put(7.1,.1){$1$}
\put(7.6,.1){$1$}

\end{picture}
\end{center}

\setlength{\unitlength}{.7cm}
\begin{center}
\begin{picture}(11,3.5)(0,0)
\thinlines

\put(0,0){\line(1,0){4}}
\put(0,0){\line(0,1){3}}

\put(0,3){\line(1,0){.5}}
\put(.5,.5){\line(0,1){2.5}}
\put(0,.5){\line(1,0){2.5}}
\put(2.5,0){\line(0,1){1.5}}
\put(2.5,1.5){\line(1,0){.5}}
\put(3,.5){\line(1,0){1}}
\put(3,.5){\line(0,1){1}}
\put(4,0){\line(0,1){.5}}
\put(4.2,0){,}

\put(.1,.1){$2$}
\put(.1,.6){$3$}
\put(.1,1.1){$3$}
\put(.1,1.6){$3$}
\put(.1,2.1){$3$}
\put(.1,2.6){$3$}
\put(.6,.1){$2$}
\put(1.1,.1){$2$}
\put(1.6,.1){$2$}
\put(2.1,.1){$2$}
\put(2.6,.1){$1$}
\put(2.6,.6){$1$}
\put(2.6,1.1){$1$}
\put(3.1,.1){$1$}
\put(3.6,.1){$1$}

\put(5.5,0){\line(1,0){5}}
\put(5.5,0){\line(0,1){3}}

\put(5.5,3){\line(1,0){.5}}
\put(5.5,.5){\line(1,0){2.5}}
\put(6,.5){\line(0,1){2.5}}

\put(8,0){\line(0,1){3}}
\put(8,1.5){\line(1,0){.5}}
\put(8,3){\line(1,0){.5}}
\put(8.5,.5){\line(0,1){2.5}}
\put(8.5,.5){\line(1,0){2}}
\put(9.5,0){\line(0,1){.5}}
\put(10.5,0){\line(0,1){.5}}
\put(10.7,0){,}

\put(5.6,.1){$2$}
\put(5.6,.6){$3$}
\put(5.6,1.1){$3$}
\put(5.6,1.6){$3$}
\put(5.6,2.1){$3$}
\put(5.6,2.6){$3$}
\put(6.1,.1){$2$}
\put(6.6,.1){$2$}
\put(7.1,.1){$2$}
\put(7.6,.1){$2$}
\put(8.1,.1){$1$}
\put(8.1,.6){$1$}
\put(8.1,1.1){$1$}
\put(8.1,1.6){$4$}
\put(8.1,2.1){$4$}
\put(8.1,2.6){$4$}
\put(8.6,.1){$1$}
\put(9.1,.1){$1$}
\put(9.6,.1){$4$}
\put(10.1,.1){$4$}

\end{picture}
\end{center}

\setlength{\unitlength}{.7cm}
\begin{center}
\begin{picture}(15,3.5)(0,0)
\thinlines

\put(0,0){\line(1,0){6}}
\put(0,0){\line(0,1){3}}
\put(0,.5){\line(1,0){.5}}
\put(0,3){\line(1,0){.5}}
\put(.5,.5){\line(0,1){2.5}}
\put(.5,.5){\line(1,0){2}}

\put(2.5,0){\line(0,1){2}}
\put(2.5,2){\line(1,0){.5}}
\put(3,.5){\line(0,1){1.5}}
\put(3,.5){\line(1,0){.5}}
\put(3.5,0){\line(0,1){3}}
\put(3.5,1.5){\line(1,0){.5}}
\put(3.5,3){\line(1,0){.5}}
\put(4,.5){\line(0,1){2.5}}
\put(4,.5){\line(1,0){2}}
\put(5,0){\line(0,1){.5}}
\put(6,0){\line(0,1){.5}}
\put(6.2,0){,}

\put(.1,.1){$2$}
\put(.1,.6){$3$}
\put(.1,1.1){$3$}
\put(.1,1.6){$3$}
\put(.1,2.1){$3$}
\put(.1,2.6){$3$}
\put(.6,.1){$2$}
\put(1.1,.1){$2$}
\put(1.6,.1){$2$}
\put(2.1,.1){$2$}
\put(2.6,.1){$5$}
\put(2.6,.6){$5$}
\put(2.6,1.1){$5$}
\put(2.6,1.6){$5$}
\put(3.1,.1){$5$}
\put(3.6,.1){$1$}
\put(3.6,.6){$1$}
\put(3.6,1.1){$1$}
\put(3.6,1.6){$4$}
\put(3.6,2.1){$4$}
\put(3.6,2.6){$4$}
\put(4.1,.1){$1$}
\put(4.6,.1){$1$}
\put(5.1,.1){$4$}
\put(5.6,.1){$4$}

\put(7,0){\line(1,0){7.5}}
\put(7,0){\line(0,1){3}}
\put(7.5,.5){\line(1,0){2}}
\put(7,.5){\line(1,0){.5}}
\put(7,3){\line(1,0){.5}}
\put(7.5,.5){\line(0,1){2.5}}

\put(9.5,0){\line(0,1){1.5}}
\put(9.5,1.5){\line(1,0){.5}}
\put(10,.5){\line(0,1){1}}
\put(10,.5){\line(1,0){1}}
\put(11,0){\line(0,1){2}}
\put(11,2){\line(1,0){.5}}
\put(11.5,.5){\line(0,1){1.5}}
\put(11.5,.5){\line(1,0){.5}}
\put(12,0){\line(0,1){3}}
\put(12,1.5){\line(1,0){.5}}
\put(12,3){\line(1,0){.5}}
\put(12.5,.5){\line(0,1){2.5}}
\put(12.5,.5){\line(1,0){2}}
\put(13.5,0){\line(0,1){.5}}
\put(14.5,0){\line(0,1){.5}}
\put(14.7,0){,}

\put(7.1,.1){$2$}
\put(7.1,.6){$3$}
\put(7.1,1.1){$3$}
\put(7.1,1.6){$3$}
\put(7.1,2.1){$3$}
\put(7.1,2.6){$3$}
\put(7.6,.1){$2$}
\put(8.1,.1){$2$}
\put(8.6,.1){$2$}
\put(9.1,.1){$2$}
\put(9.6,.1){$6$}
\put(9.6,.6){$6$}
\put(9.6,1.1){$6$}
\put(10.1,.1){$6$}
\put(10.6,.1){$6$}
\put(11.1,.1){$5$}
\put(11.1,.6){$5$}
\put(11.1,1.1){$5$}
\put(11.1,1.6){$5$}
\put(11.6,.1){$5$}
\put(12.1,.1){$1$}
\put(12.1,.6){$1$}
\put(12.1,1.1){$1$}
\put(12.1,1.6){$4$}
\put(12.1,2.1){$4$}
\put(12.1,2.6){$4$}
\put(12.6,.1){$1$}
\put(13.1,.1){$1$}
\put(13.6,.1){$4$}
\put(14.1,.1){$4$}

\end{picture}
\end{center}

\setlength{\unitlength}{.7cm}
\begin{center}
\begin{picture}(9,3.5)(0,0)
\thinlines

\put(-1.3,1){$Q=$}
\put(0,0){\line(1,0){8.5}}
\put(0,0){\line(0,1){3}}
\put(.5,.5){\line(1,0){2}}
\put(0,.5){\line(1,0){.5}}
\put(0,3){\line(1,0){.5}}
\put(.5,.5){\line(0,1){2.5}}

\put(2.5,0){\line(0,1){3}}
\put(2.5,1.5){\line(1,0){.5}}
\put(2.5,3){\line(1,0){.5}}
\put(3,.5){\line(0,1){2.5}}
\put(3,.5){\line(1,0){2}}
\put(4,0){\line(0,1){.5}}

\put(5,0){\line(0,1){2}}
\put(5,2){\line(1,0){.5}}
\put(5.5,.5){\line(0,1){1.5}}
\put(5.5,.5){\line(1,0){.5}}
\put(6,0){\line(0,1){3}}
\put(6,1.5){\line(1,0){.5}}
\put(6,3){\line(1,0){.5}}
\put(6.5,.5){\line(0,1){2.5}}
\put(6.5,.5){\line(1,0){2}}
\put(7.5,0){\line(0,1){.5}}
\put(8.5,0){\line(0,1){.5}}

\put(.1,.1){$2$}
\put(.1,.6){$3$}
\put(.1,1.1){$3$}
\put(.1,1.6){$3$}
\put(.1,2.1){$3$}
\put(.1,2.6){$3$}
\put(.6,.1){$2$}
\put(1.1,.1){$2$}
\put(1.6,.1){$2$}
\put(2.1,.1){$2$}
\put(2.6,.1){$6$}
\put(2.6,.6){$6$}
\put(2.6,1.1){$6$}
\put(2.6,1.6){$7$}
\put(2.6,2.1){$7$}
\put(2.6,2.6){$7$}
\put(3.1,.1){$6$}
\put(3.6,.1){$6$}
\put(4.1,.1){$7$}
\put(4.6,.1){$7$}
\put(5.1,.1){$5$}
\put(5.1,.6){$5$}
\put(5.1,1.1){$5$}
\put(5.1,1.6){$5$}
\put(5.6,.1){$5$}
\put(6.1,.1){$1$}
\put(6.1,.6){$1$}
\put(6.1,1.1){$1$}
\put(6.1,1.6){$4$}
\put(6.1,2.1){$4$}
\put(6.1,2.6){$4$}
\put(6.6,.1){$1$}
\put(7.1,.1){$1$}
\put(7.6,.1){$4$}
\put(8.1,.1){$4$}

\end{picture}
\end{center}

Notice that the $Q$ tableau obtained from the insertion method is the same as the $\hat{Q}$ tableaux obtained from the square diagram, but that $P \neq \hat{P}$.

\begin{Theorem}
The $k$-ribbon insertion algorithm is a bijection between the set of $k$-colored permutations and the set of pairs $(P, Q)$ with $P$ a standard $k$-ribbon Fibonacci tableau and $Q$ a $k$-ribbon Fibonacci path tableau.  
\end{Theorem}

\begin{proof}
We claim that the insertion procedure defined above is invertible.  At the $i$th stage of the insertion, the $Q$ tableau tells us which $k$-ribbon was the most recently created in the tableau $P_i$.  If this $k$-ribbon was added on top of another $k$-ribbon, then the shape of $P_i$ must have had a shape bijectively equivalent to $2^i \omega$ for some word $\omega$ of $1_1$'s, $1_2$'s, $\dots$, $1_k$'s, and $2$'s.  Then to work backwards, each $k$-ribbon in the top row will bump to the left, preserving their height, until the leftmost $k$-ribbon in the top row is bumped out of the tableau as a $k$-ribbon of height $j$.  If this $k$-ribbon contained $x_i$'s then $x_i^j$ was the element inserted at this step of the algorithm.  

If the $k$-ribbon that was most recently created was not added on top of another $k$-ribbon, then the shape of $P_i$ is bijectively equivalent to $2^{l-1} 1_m \omega$ for some $m$.  In both cases, the element inside the newly created $k$-ribbon, say $t_i$, is smaller than the element inside the bottom $k$-ribbon to the left of it.  When we reverse the bumping algorithm, the $k$-ribbon containing $t_i$ will bump to the top $k$-ribbon in the stack to the left of it and each $k$-ribbon in the top row will bump to the left, preserving their height until the leftmost $k$-ribbon in the top row is bumped out of the tableau, again preserving its height.   If $l=1$, then the newly created $k$-ribbon in the first row is itself bumped out of the tableau.  If this $k$-ribbon that is bumped has height $j$ and contains $x_i$'s, then $x_i^j$ is the element that was inserted at this step.  

In either case, we obtain the originally inserted element along with its shape and $P_{i-1}$.
\end{proof}

\section{A $k$-Ribbon Fibonacci Evacuation}
\label{evacuation}

An evacuation algorithm for standard Fibonacci tableaux in $Z(1)$ has been given in by Killpatrick \cite{Kil} and a more general evacuation algorithm for standard domino Fibonacci tableaux has been given by Cameron and Killpatrick \cite{CamKill}.  Here we give an evacuation algorithm for standard $k$-ribbon Fibonacci tableaux which, for $k=2$, restricts to the evacuation algorithm given for domino tableaux.  This evacuation algorithm will give a bijection between standard $k$-ribbon Fibonacci tableaux and $k$-ribbon Fibonacci path tableaux.  

Let $P$ be a standard $k$-ribbon Fibonacci tableaux.  To compute the evacuation of $P$, do the following:

\begin{enumerate}

\item
Erase the number in the bottom $k$-ribbon of the leftmost column of $P$.  By definition of a standard $k$-ribbon Fibonacci tableaux, this will be the largest number in $P$.  

\item
As long as there is a $k$-ribbon containing $a$'s above the empty $k$-ribbon, compare $a$ with the number in the bottom $k$-ribbon in the column to the right of the empty $k$-ribbon, call it $b$.  

\begin{enumerate}

\item
If $a>b$, replace the empty $k$-ribbon with a $k$-ribbon of the same height as the $k$-ribbon containing $a$.  This leaves an empty $k$-ribbon of some height on the top of the first column.  

\item
If $a<b$, then fill the empty $k$-ribbon with $b$'s, leaving an empty $k$-ribbon in the second column.

\end{enumerate}

\item
Continue the algorithm until reaching an empty $k$-ribbon that has no $k$-ribbon immediately above it.  Remove the empty $k$-ribbon from the tableau and, if necessary, slide all remaining columns one column to the left so that the result has the shape of a $k$-ribbon Fibonacci tableau.  Call this remaining tableau $P^{(1)}$.  

\item
In a new tableau with the same shape as $P$, called $\tilde{P}$, put $n$'s into the position of the last empty $k$-ribbon.

\item
Create $P^{(2)}$ by iterating the above procedure, starting with $P^{(1)}$.  At the fourth step, place $n-1$'s into $\tilde{P}$ in the position of the last empty $k$-ribbon.  Continue the procedure until $P^{(n)} = \emptyset$ and $\tilde{P}$ is a $k$-ribbon Fibonacci path tableau.  The tableau $\tilde{P}$ will be called the evacuation tableau of $P$ and will be denoted by $ev(P)$.

\end{enumerate}

For example, if we begin with the standard $k$-ribbon Fibonacci tableau

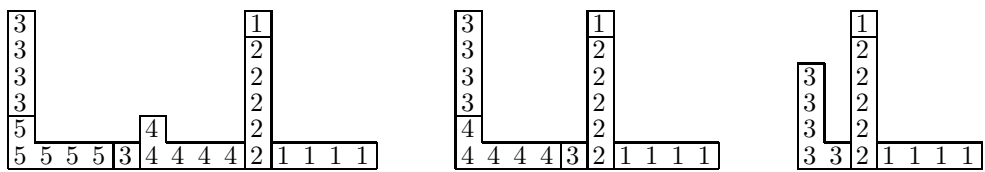
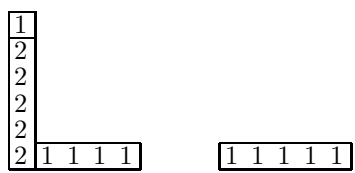
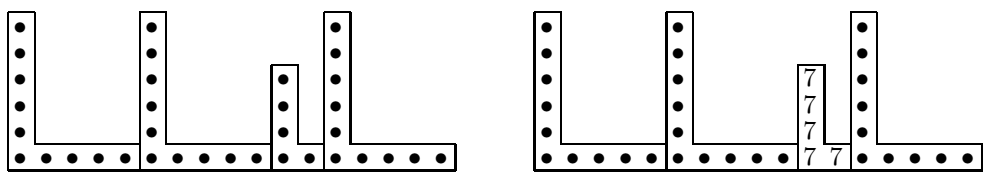
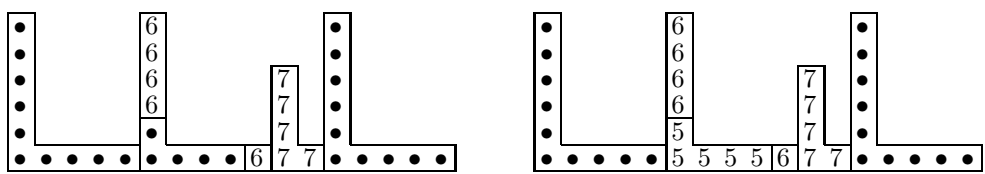
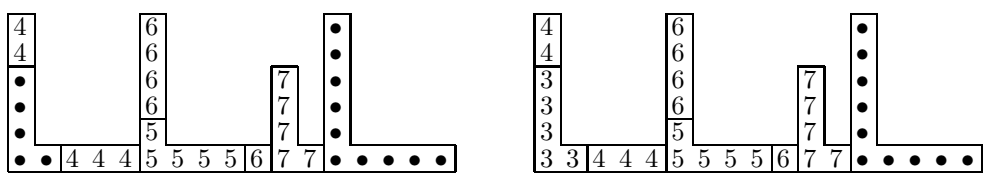
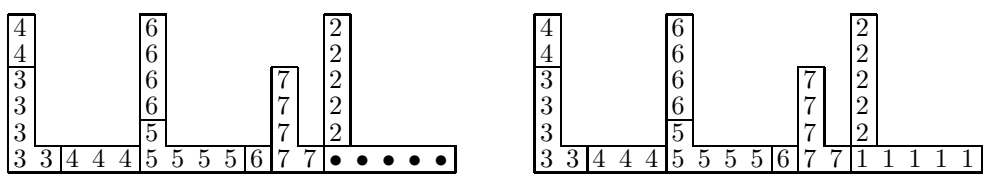
\setlength{\unitlength}{.7cm}
\begin{center}
\begin{picture}(9,3.5)(0,0)
\thinlines

\put(-2,1){$P(\pi)=$}

\put(0,0){\line(1,0){8.5}}
\put(0,0){\line(0,1){3}}
\put(.5,.5){\line(1,0){2}}
\put(0,1){\line(1,0){.5}}
\put(0,3){\line(1,0){.5}}
\put(.5,.5){\line(0,1){2.5}}
\put(2,0){\line(0,1){.5}}

\put(2.5,0){\line(0,1){3}}
\put(2.5,2){\line(1,0){.5}}
\put(2.5,3){\line(1,0){.5}}
\put(3,.5){\line(0,1){2.5}}
\put(3,.5){\line(1,0){2}}
\put(3.5,0){\line(0,1){.5}}

\put(5,0){\line(0,1){2}}
\put(5,2){\line(1,0){.5}}
\put(5.5,.5){\line(0,1){1.5}}
\put(5.5, .5){\line(1,0){.5}}

\put(6,0){\line(0,1){3}}
\put(6,2.5){\line(1,0){.5}}
\put(6,3){\line(1,0){.5}}
\put(6.5,0){\line(0,1){3}}
\put(6.5,.5){\line(1,0){2}}
\put(8.5,0){\line(0,1){.5}}

\put(.1,.1){$7$}
\put(.1,.6){$7$}
\put(.1,1.1){$3$}
\put(.1,1.6){$3$}
\put(.1,2.1){$3$}
\put(.1,2.6){$3$}
\put(.6,.1){$7$}
\put(1.1,.1){$7$}
\put(1.6,.1){$7$}
\put(2.1,.1){$3$}
\put(2.6,.1){$6$}
\put(2.6,.6){$6$}
\put(2.6,1.1){$6$}
\put(2.6,1.6){$6$}
\put(2.6,2.1){$4$}
\put(2.6,2.6){$4$}
\put(3.1,.1){$6$}
\put(3.6,.1){$4$}
\put(4.1,.1){$4$}
\put(4.6,.1){$4$}
\put(5.1,.1){$5$}
\put(5.1,.6){$5$}
\put(5.1,1.1){$5$}
\put(5.1,1.6){$5$}
\put(5.6,.1){$5$}
\put(6.1,.1){$2$}
\put(6.1,.6){$2$}
\put(6.1,1.1){$2$}
\put(6.1,1.6){$2$}
\put(6.1,2.1){$2$}
\put(6.1,2.6){$1$}
\put(6.6,.1){$1$}
\put(7.1,.1){$1$}
\put(7.6,.1){$1$}
\put(8.1,.1){$1$}

\end{picture}
\end{center}

then the first sequence of steps in the evacuation algorithm is

\setlength{\unitlength}{.7cm}
\begin{center}
\begin{picture}(20,3.5)(0,0)
\thinlines

\put(0,0){\line(1,0){8.5}}
\put(0,0){\line(0,1){3}}
\put(.5,.5){\line(1,0){2}}
\put(0,1){\line(1,0){.5}}
\put(0,3){\line(1,0){.5}}
\put(.5,.5){\line(0,1){2.5}}
\put(2,0){\line(0,1){.5}}

\put(2.5,0){\line(0,1){3}}
\put(2.5,2){\line(1,0){.5}}
\put(2.5,3){\line(1,0){.5}}
\put(3,.5){\line(0,1){2.5}}
\put(3,.5){\line(1,0){2}}
\put(3.5,0){\line(0,1){.5}}

\put(5,0){\line(0,1){2}}
\put(5,2){\line(1,0){.5}}
\put(5.5,.5){\line(0,1){1.5}}
\put(5.5, .5){\line(1,0){.5}}

\put(6,0){\line(0,1){3}}
\put(6,2.5){\line(1,0){.5}}
\put(6,3){\line(1,0){.5}}
\put(6.5,0){\line(0,1){3}}
\put(6.5,.5){\line(1,0){2}}
\put(8.5,0){\line(0,1){.5}}

\put(.1,.1){$\bullet$}
\put(.1,.6){$\bullet$}
\put(.1,1.1){$3$}
\put(.1,1.6){$3$}
\put(.1,2.1){$3$}
\put(.1,2.6){$3$}
\put(.6,.1){$\bullet$}
\put(1.1,.1){$\bullet$}
\put(1.6,.1){$\bullet$}
\put(2.1,.1){$3$}
\put(2.6,.1){$6$}
\put(2.6,.6){$6$}
\put(2.6,1.1){$6$}
\put(2.6,1.6){$6$}
\put(2.6,2.1){$4$}
\put(2.6,2.6){$4$}
\put(3.1,.1){$6$}
\put(3.6,.1){$4$}
\put(4.1,.1){$4$}
\put(4.6,.1){$4$}
\put(5.1,.1){$5$}
\put(5.1,.6){$5$}
\put(5.1,1.1){$5$}
\put(5.1,1.6){$5$}
\put(5.6,.1){$5$}
\put(6.1,.1){$2$}
\put(6.1,.6){$2$}
\put(6.1,1.1){$2$}
\put(6.1,1.6){$2$}
\put(6.1,2.1){$2$}
\put(6.1,2.6){$1$}
\put(6.6,.1){$1$}
\put(7.1,.1){$1$}
\put(7.6,.1){$1$}
\put(8.1,.1){$1$}

\put(10,0){\line(1,0){8.5}}
\put(10,0){\line(0,1){3}}
\put(10.5,.5){\line(1,0){2}}
\put(10,1){\line(1,0){.5}}
\put(10,3){\line(1,0){.5}}
\put(10.5,.5){\line(0,1){2.5}}
\put(12,0){\line(0,1){.5}}

\put(12.5,0){\line(0,1){3}}
\put(12.5,2){\line(1,0){.5}}
\put(12.5,3){\line(1,0){.5}}
\put(13,.5){\line(0,1){2.5}}
\put(13,.5){\line(1,0){2}}
\put(13.5,0){\line(0,1){.5}}

\put(15,0){\line(0,1){2}}
\put(15,2){\line(1,0){.5}}
\put(15.5,.5){\line(0,1){1.5}}
\put(15.5, .5){\line(1,0){.5}}

\put(16,0){\line(0,1){3}}
\put(16,2.5){\line(1,0){.5}}
\put(16,3){\line(1,0){.5}}
\put(16.5,0){\line(0,1){3}}
\put(16.5,.5){\line(1,0){2}}
\put(18.5,0){\line(0,1){.5}}

\put(10.1,.1){$6$}
\put(10.1,.6){$6$}
\put(10.1,1.1){$3$}
\put(10.1,1.6){$3$}
\put(10.1,2.1){$3$}
\put(10.1,2.6){$3$}
\put(10.6,.1){$6$}
\put(11.1,.1){$6$}
\put(11.6,.1){$6$}
\put(12.1,.1){$3$}
\put(12.6,.1){$\bullet$}
\put(12.6,.6){$\bullet$}
\put(12.6,1.1){$\bullet$}
\put(12.6,1.6){$\bullet$}
\put(12.6,2.1){$4$}
\put(12.6,2.6){$4$}
\put(13.1,.1){$\bullet$}
\put(13.6,.1){$4$}
\put(14.1,.1){$4$}
\put(14.6,.1){$4$}
\put(15.1,.1){$5$}
\put(15.1,.6){$5$}
\put(15.1,1.1){$5$}
\put(15.1,1.6){$5$}
\put(15.6,.1){$5$}
\put(16.1,.1){$2$}
\put(16.1,.6){$2$}
\put(16.1,1.1){$2$}
\put(16.1,1.6){$2$}
\put(16.1,2.1){$2$}
\put(16.1,2.6){$1$}
\put(16.6,.1){$1$}
\put(17.1,.1){$1$}
\put(17.6,.1){$1$}
\put(18.1,.1){$1$}

\end{picture}
\end{center}

\setlength{\unitlength}{.7cm}
\begin{center}
\begin{picture}(20,3.5)(0,0)
\thinlines

\put(0,0){\line(1,0){8.5}}
\put(0,0){\line(0,1){3}}
\put(.5,.5){\line(1,0){2}}
\put(0,1){\line(1,0){.5}}
\put(0,3){\line(1,0){.5}}
\put(.5,.5){\line(0,1){2.5}}
\put(2,0){\line(0,1){.5}}

\put(2.5,0){\line(0,1){3}}
\put(2.5,2){\line(1,0){.5}}
\put(2.5,3){\line(1,0){.5}}
\put(3,.5){\line(0,1){2.5}}
\put(3,.5){\line(1,0){2}}
\put(3.5,0){\line(0,1){.5}}

\put(5,0){\line(0,1){2}}
\put(5,2){\line(1,0){.5}}
\put(5.5,.5){\line(0,1){1.5}}
\put(5.5, .5){\line(1,0){.5}}

\put(6,0){\line(0,1){3}}
\put(6,2.5){\line(1,0){.5}}
\put(6,3){\line(1,0){.5}}
\put(6.5,0){\line(0,1){3}}
\put(6.5,.5){\line(1,0){2}}
\put(8.5,0){\line(0,1){.5}}

\put(.1,.1){$6$}
\put(.1,.6){$6$}
\put(.1,1.1){$3$}
\put(.1,1.6){$3$}
\put(.1,2.1){$3$}
\put(.1,2.6){$3$}
\put(.6,.1){$6$}
\put(1.1,.1){$6$}
\put(1.6,.1){$6$}
\put(2.1,.1){$3$}
\put(2.6,.1){$5$}
\put(2.6,.6){$5$}
\put(2.6,1.1){$5$}
\put(2.6,1.6){$5$}
\put(2.6,2.1){$4$}
\put(2.6,2.6){$4$}
\put(3.1,.1){$5$}
\put(3.6,.1){$4$}
\put(4.1,.1){$4$}
\put(4.6,.1){$4$}
\put(5.1,.1){$\bullet$}
\put(5.1,.6){$\bullet$}
\put(5.1,1.1){$\bullet$}
\put(5.1,1.6){$\bullet$}
\put(5.6,.1){$\bullet$}
\put(6.1,.1){$2$}
\put(6.1,.6){$2$}
\put(6.1,1.1){$2$}
\put(6.1,1.6){$2$}
\put(6.1,2.1){$2$}
\put(6.1,2.6){$1$}
\put(6.6,.1){$1$}
\put(7.1,.1){$1$}
\put(7.6,.1){$1$}
\put(8.1,.1){$1$}

\put(9.2,1){$P^{(1)}=$}

\put(11,0){\line(1,0){7.5}}
\put(11,0){\line(0,1){3}}
\put(11.5,.5){\line(1,0){2}}
\put(11,1){\line(1,0){.5}}
\put(11,3){\line(1,0){.5}}
\put(11.5,.5){\line(0,1){2.5}}
\put(13,0){\line(0,1){.5}}

\put(13.5,0){\line(0,1){3}}
\put(13.5,2){\line(1,0){.5}}
\put(13.5,3){\line(1,0){.5}}
\put(14,.5){\line(0,1){2.5}}
\put(14,.5){\line(1,0){2}}
\put(14.5,0){\line(0,1){.5}}

\put(16,0){\line(0,1){3}}
\put(16,2.5){\line(1,0){.5}}
\put(16,3){\line(1,0){.5}}
\put(16.5,0){\line(0,1){3}}
\put(16.5,.5){\line(1,0){2}}
\put(18.5,0){\line(0,1){.5}}

\put(11.1,.1){$6$}
\put(11.1,.6){$6$}
\put(11.1,1.1){$3$}
\put(11.1,1.6){$3$}
\put(11.1,2.1){$3$}
\put(11.1,2.6){$3$}
\put(11.6,.1){$6$}
\put(12.1,.1){$6$}
\put(12.6,.1){$6$}
\put(13.1,.1){$3$}
\put(13.6,.1){$5$}
\put(13.6,.6){$5$}
\put(13.6,1.1){$5$}
\put(13.6,1.6){$5$}
\put(13.6,2.1){$4$}
\put(13.6,2.6){$4$}
\put(14.1,.1){$5$}
\put(14.6,.1){$4$}
\put(15.1,.1){$4$}
\put(15.6,.1){$4$}

\put(16.1,.1){$2$}
\put(16.1,.6){$2$}
\put(16.1,1.1){$2$}
\put(16.1,1.6){$2$}
\put(16.1,2.1){$2$}
\put(16.1,2.6){$1$}
\put(16.6,.1){$1$}
\put(17.1,.1){$1$}
\put(17.6,.1){$1$}
\put(18.1,.1){$1$}

\end{picture}
\end{center}

After this first sequence of steps, the evacuation tabelaux $ev(P)$ looks like

\setlength{\unitlength}{.7cm}
\begin{center}
\begin{picture}(9,3.5)(0,0)
\thinlines

\put(0,0){\line(1,0){8.5}}
\put(0,0){\line(0,1){3}}
\put(.5,.5){\line(1,0){2}}

\put(0,3){\line(1,0){.5}}
\put(.5,.5){\line(0,1){2.5}}

\put(2.5,0){\line(0,1){3}}

\put(2.5,3){\line(1,0){.5}}
\put(3,.5){\line(0,1){2.5}}
\put(3,.5){\line(1,0){2}}

\put(5,0){\line(0,1){2}}
\put(5,2){\line(1,0){.5}}
\put(5.5,.5){\line(0,1){1.5}}
\put(5.5, .5){\line(1,0){.5}}

\put(6,0){\line(0,1){3}}

\put(6,3){\line(1,0){.5}}
\put(6.5,.5){\line(0,1){2.5}}
\put(6.5,.5){\line(1,0){2}}
\put(8.5,0){\line(0,1){.5}}

\put(.1,.1){$\bullet$}
\put(.1,.6){$\bullet$}
\put(.1,1.1){$\bullet$}
\put(.1,1.6){$\bullet$}
\put(.1,2.1){$\bullet$}
\put(.1,2.6){$\bullet$}
\put(.6,.1){$\bullet$}
\put(1.1,.1){$\bullet$}
\put(1.6,.1){$\bullet$}
\put(2.1,.1){$\bullet$}
\put(2.6,.1){$\bullet$}
\put(2.6,.6){$\bullet$}
\put(2.6,1.1){$\bullet$}
\put(2.6,1.6){$\bullet$}
\put(2.6,2.1){$\bullet$}
\put(2.6,2.6){$\bullet$}
\put(3.1,.1){$\bullet$}
\put(3.6,.1){$\bullet$}
\put(4.1,.1){$\bullet$}
\put(4.6,.1){$\bullet$}
\put(5.1,.1){$7$}
\put(5.1,.6){$7$}
\put(5.1,1.1){$7$}
\put(5.1,1.6){$7$}
\put(5.6,.1){$7$}
\put(6.1,.1){$\bullet$}
\put(6.1,.6){$\bullet$}
\put(6.1,1.1){$\bullet$}
\put(6.1,1.6){$\bullet$}
\put(6.1,2.1){$\bullet$}
\put(6.1,2.6){$\bullet$}
\put(6.6,.1){$\bullet$}
\put(7.1,.1){$\bullet$}
\put(7.6,.1){$\bullet$}
\put(8.1,.1){$\bullet$}

\end{picture}
\end{center}

We now give all of the steps in the development of $ev(P)$:

\setlength{\unitlength}{.7cm}
\begin{center}
\begin{picture}(20,3.5)(0,0)
\thinlines

\put(-1.5, 1){$P^{(k)}:$}

\put(1,0){\line(1,0){8.5}}
\put(1,0){\line(0,1){3}}
\put(1.5,.5){\line(1,0){2}}
\put(1,1){\line(1,0){.5}}
\put(1,3){\line(1,0){.5}}
\put(1.5,.5){\line(0,1){2.5}}
\put(3,0){\line(0,1){.5}}

\put(3.5,0){\line(0,1){3}}
\put(3.5,2){\line(1,0){.5}}
\put(3.5,3){\line(1,0){.5}}
\put(4,.5){\line(0,1){2.5}}
\put(4,.5){\line(1,0){2}}
\put(4.5,0){\line(0,1){.5}}

\put(6,0){\line(0,1){2}}
\put(6,2){\line(1,0){.5}}
\put(6.5,.5){\line(0,1){1.5}}
\put(6.5, .5){\line(1,0){.5}}

\put(7,0){\line(0,1){3}}
\put(7,2.5){\line(1,0){.5}}
\put(7,3){\line(1,0){.5}}
\put(7.5,0){\line(0,1){3}}
\put(7.5,.5){\line(1,0){2}}
\put(9.5,0){\line(0,1){.5}}

\put(1.1,.1){$7$}
\put(1.1,.6){$7$}
\put(1.1,1.1){$3$}
\put(1.1,1.6){$3$}
\put(1.1,2.1){$3$}
\put(1.1,2.6){$3$}
\put(1.6,.1){$7$}
\put(2.1,.1){$7$}
\put(2.6,.1){$7$}
\put(3.1,.1){$3$}
\put(3.6,.1){$6$}
\put(3.6,.6){$6$}
\put(3.6,1.1){$6$}
\put(3.6,1.6){$6$}
\put(3.6,2.1){$4$}
\put(3.6,2.6){$4$}
\put(4.1,.1){$6$}
\put(4.6,.1){$4$}
\put(5.1,.1){$4$}
\put(5.6,.1){$4$}
\put(6.1,.1){$5$}
\put(6.1,.6){$5$}
\put(6.1,1.1){$5$}
\put(6.1,1.6){$5$}
\put(6.6,.1){$5$}
\put(7.1,.1){$2$}
\put(7.1,.6){$2$}
\put(7.1,1.1){$2$}
\put(7.1,1.6){$2$}
\put(7.1,2.1){$2$}
\put(7.1,2.6){$1$}
\put(7.6,.1){$1$}
\put(8.1,.1){$1$}
\put(8.6,.1){$1$}
\put(9.1,.1){$1$}

\put(11,0){\line(1,0){7.5}}
\put(11,0){\line(0,1){3}}
\put(11.5,.5){\line(1,0){2}}
\put(11,1){\line(1,0){.5}}
\put(11,3){\line(1,0){.5}}
\put(11.5,.5){\line(0,1){2.5}}
\put(13,0){\line(0,1){.5}}

\put(13.5,0){\line(0,1){3}}
\put(13.5,2){\line(1,0){.5}}
\put(13.5,3){\line(1,0){.5}}
\put(14,.5){\line(0,1){2.5}}
\put(14,.5){\line(1,0){2}}
\put(14.5,0){\line(0,1){.5}}

\put(16,0){\line(0,1){3}}
\put(16,2.5){\line(1,0){.5}}
\put(16,3){\line(1,0){.5}}
\put(16.5,0){\line(0,1){3}}
\put(16.5,.5){\line(1,0){2}}
\put(18.5,0){\line(0,1){.5}}

\put(11.1,.1){$6$}
\put(11.1,.6){$6$}
\put(11.1,1.1){$3$}
\put(11.1,1.6){$3$}
\put(11.1,2.1){$3$}
\put(11.1,2.6){$3$}
\put(11.6,.1){$6$}
\put(12.1,.1){$6$}
\put(12.6,.1){$6$}
\put(13.1,.1){$3$}
\put(13.6,.1){$5$}
\put(13.6,.6){$5$}
\put(13.6,1.1){$5$}
\put(13.6,1.6){$5$}
\put(13.6,2.1){$4$}
\put(13.6,2.6){$4$}
\put(14.1,.1){$5$}
\put(14.6,.1){$4$}
\put(15.1,.1){$4$}
\put(15.6,.1){$4$}

\put(16.1,.1){$2$}
\put(16.1,.6){$2$}
\put(16.1,1.1){$2$}
\put(16.1,1.6){$2$}
\put(16.1,2.1){$2$}
\put(16.1,2.6){$1$}
\put(16.6,.1){$1$}
\put(17.1,.1){$1$}
\put(17.6,.1){$1$}
\put(18.1,.1){$1$}

\end{picture}
\end{center}

\setlength{\unitlength}{.7cm}
\begin{center}
\begin{picture}(20,3.5)(0,0)
\thinlines

\put(0,0){\line(1,0){7}}
\put(0,0){\line(0,1){3}}
\put(.5,.5){\line(1,0){2}}
\put(0,1){\line(1,0){.5}}
\put(0,3){\line(1,0){.5}}
\put(0.5,.5){\line(0,1){2.5}}
\put(2,0){\line(0,1){.5}}

\put(2.5,0){\line(0,1){1}}
\put(2.5,1){\line(1,0){.5}}

\put(3,.5){\line(0,1){.5}}
\put(3,.5){\line(1,0){1.5}}

\put(4.5,0){\line(0,1){3}}
\put(4.5,2.5){\line(1,0){.5}}
\put(4.5,3){\line(1,0){.5}}
\put(5,0){\line(0,1){3}}
\put(5,.5){\line(1,0){2}}
\put(7,0){\line(0,1){.5}}

\put(.1,.1){$5$}
\put(.1,.6){$5$}
\put(.1,1.1){$3$}
\put(.1,1.6){$3$}
\put(.1,2.1){$3$}
\put(.1,2.6){$3$}
\put(.6,.1){$5$}
\put(1.1,.1){$5$}
\put(1.6,.1){$5$}
\put(2.1,.1){$3$}
\put(2.6,.1){$4$}
\put(2.6,.6){$4$}
\put(3.1,.1){$4$}
\put(3.6,.1){$4$}
\put(4.1,.1){$4$}

\put(4.6,.1){$2$}
\put(4.6,.6){$2$}
\put(4.6,1.1){$2$}
\put(4.6,1.6){$2$}
\put(4.6,2.1){$2$}
\put(4.6,2.6){$1$}
\put(5.1,.1){$1$}
\put(5.6,.1){$1$}
\put(6.1,.1){$1$}
\put(6.6,.1){$1$}

\put(8.5,0){\line(1,0){5}}
\put(8.5,0){\line(0,1){3}}
\put(9,.5){\line(1,0){2}}
\put(8.5,1){\line(1,0){.5}}
\put(8.5,3){\line(1,0){.5}}
\put(9,.5){\line(0,1){2.5}}
\put(10.5,0){\line(0,1){.5}}

\put(11,0){\line(0,1){3}}
\put(11,2.5){\line(1,0){.5}}
\put(11,3){\line(1,0){.5}}
\put(11.5,0){\line(0,1){3}}
\put(11.5,.5){\line(1,0){2}}
\put(13.5,0){\line(0,1){.5}}

\put(8.6,.1){$4$}
\put(8.6,.6){$4$}
\put(8.6,1.1){$3$}
\put(8.6,1.6){$3$}
\put(8.6,2.1){$3$}
\put(8.6,2.6){$3$}
\put(9.1,.1){$4$}
\put(9.6,.1){$4$}
\put(10.1,.1){$4$}
\put(10.6,.1){$3$}

\put(11.1,.1){$2$}
\put(11.1,.6){$2$}
\put(11.1,1.1){$2$}
\put(11.1,1.6){$2$}
\put(11.1,2.1){$2$}
\put(11.1,2.6){$1$}
\put(11.6,.1){$1$}
\put(12.1,.1){$1$}
\put(12.6,.1){$1$}
\put(13.1,.1){$1$}

\put(15,0){\line(1,0){3.5}}
\put(15,0){\line(0,1){2}}
\put(15.5,.5){\line(1,0){.5}}
\put(15,2){\line(1,0){.5}}
\put(15.5,.5){\line(0,1){1.5}}

\put(16,0){\line(0,1){3}}
\put(16,2.5){\line(1,0){.5}}
\put(16,3){\line(1,0){.5}}
\put(16.5,0){\line(0,1){3}}
\put(16.5,.5){\line(1,0){2}}
\put(18.5,0){\line(0,1){.5}}

\put(15.1,.1){$3$}
\put(15.1,.6){$3$}
\put(15.1,1.1){$3$}
\put(15.1,1.6){$3$}
\put(15.6,.1){$3$}

\put(16.1,.1){$2$}
\put(16.1,.6){$2$}
\put(16.1,1.1){$2$}
\put(16.1,1.6){$2$}
\put(16.1,2.1){$2$}
\put(16.1,2.6){$1$}
\put(16.6,.1){$1$}
\put(17.1,.1){$1$}
\put(17.6,.1){$1$}
\put(18.1,.1){$1$}

\end{picture}
\end{center}

\setlength{\unitlength}{.7cm}
\begin{center}
\begin{picture}(20,3.5)(0,0)
\thinlines

\put(6,0){\line(1,0){2.5}}
\put(6,0){\line(0,1){3}}
\put(6,2.5){\line(1,0){.5}}
\put(6,3){\line(1,0){.5}}
\put(6.5,0){\line(0,1){3}}
\put(6.5,.5){\line(1,0){2}}
\put(8.5,0){\line(0,1){.5}}

\put(6.1,.1){$2$}
\put(6.1,.6){$2$}
\put(6.1,1.1){$2$}
\put(6.1,1.6){$2$}
\put(6.1,2.1){$2$}
\put(6.1,2.6){$1$}
\put(6.6,.1){$1$}
\put(7.1,.1){$1$}
\put(7.6,.1){$1$}
\put(8.1,.1){$1$}

\put(10,0){\line(1,0){2.5}}
\put(10,0){\line(0,1){.5}}
\put(10,.5){\line(1,0){2.5}}
\put(12.5,0){\line(0,1){.5}}

\put(10.1,.1){$1$}
\put(10.6,.1){$1$}
\put(11.1,.1){$1$}
\put(11.6,.1){$1$}
\put(12.1,.1){$1$}

\end{picture}
\end{center}

\setlength{\unitlength}{.7cm}
\begin{center}
\begin{picture}(20,3.5)(0,0)
\thinlines

\put(-1.5, 1){$\tilde{P}:$}

\put(0,0){\line(1,0){8.5}}
\put(0,0){\line(0,1){3}}
\put(.5,.5){\line(1,0){2}}
\put(0,3){\line(1,0){.5}}
\put(.5,.5){\line(0,1){2.5}}

\put(2.5,0){\line(0,1){3}}
\put(2.5,3){\line(1,0){.5}}
\put(3,.5){\line(0,1){2.5}}
\put(3,.5){\line(1,0){2}}

\put(5,0){\line(0,1){2}}
\put(5,2){\line(1,0){.5}}
\put(5.5,.5){\line(0,1){1.5}}
\put(5.5, .5){\line(1,0){.5}}

\put(6,0){\line(0,1){3}}
\put(6,3){\line(1,0){.5}}
\put(6.5,.5){\line(0,1){2.5}}
\put(6.5,.5){\line(1,0){2}}
\put(8.5,0){\line(0,1){.5}}

\put(.1,.1){$\bullet$}
\put(.1,.6){$\bullet$}
\put(.1,1.1){$\bullet$}
\put(.1,1.6){$\bullet$}
\put(.1,2.1){$\bullet$}
\put(.1,2.6){$\bullet$}
\put(.6,.1){$\bullet$}
\put(1.1,.1){$\bullet$}
\put(1.6,.1){$\bullet$}
\put(2.1,.1){$\bullet$}
\put(2.6,.1){$\bullet$}
\put(2.6,.6){$\bullet$}
\put(2.6,1.1){$\bullet$}
\put(2.6,1.6){$\bullet$}
\put(2.6,2.1){$\bullet$}
\put(2.6,2.6){$\bullet$}
\put(3.1,.1){$\bullet$}
\put(3.6,.1){$\bullet$}
\put(4.1,.1){$\bullet$}
\put(4.6,.1){$\bullet$}
\put(5.1,.1){$\bullet$}
\put(5.1,.6){$\bullet$}
\put(5.1,1.1){$\bullet$}
\put(5.1,1.6){$\bullet$}
\put(5.6,.1){$\bullet$}
\put(6.1,.1){$\bullet$}
\put(6.1,.6){$\bullet$}
\put(6.1,1.1){$\bullet$}
\put(6.1,1.6){$\bullet$}
\put(6.1,2.1){$\bullet$}
\put(6.1,2.6){$\bullet$}
\put(6.6,.1){$\bullet$}
\put(7.1,.1){$\bullet$}
\put(7.6,.1){$\bullet$}
\put(8.1,.1){$\bullet$}

\put(10,0){\line(1,0){8.5}}
\put(10,0){\line(0,1){3}}
\put(10.5,.5){\line(1,0){2}}
\put(10,3){\line(1,0){.5}}
\put(10.5,.5){\line(0,1){2.5}}

\put(12.5,0){\line(0,1){3}}
\put(12.5,3){\line(1,0){.5}}
\put(13,.5){\line(0,1){2.5}}
\put(13,.5){\line(1,0){2}}

\put(15,0){\line(0,1){2}}
\put(15,2){\line(1,0){.5}}
\put(15.5,.5){\line(0,1){1.5}}
\put(15.5, .5){\line(1,0){.5}}

\put(16,0){\line(0,1){3}}
\put(16,3){\line(1,0){.5}}
\put(16.5,.5){\line(0,1){2.5}}
\put(16.5,.5){\line(1,0){2}}
\put(18.5,0){\line(0,1){.5}}

\put(10.1,.1){$\bullet$}
\put(10.1,.6){$\bullet$}
\put(10.1,1.1){$\bullet$}
\put(10.1,1.6){$\bullet$}
\put(10.1,2.1){$\bullet$}
\put(10.1,2.6){$\bullet$}
\put(10.6,.1){$\bullet$}
\put(11.1,.1){$\bullet$}
\put(11.6,.1){$\bullet$}
\put(12.1,.1){$\bullet$}
\put(12.6,.1){$\bullet$}
\put(12.6,.6){$\bullet$}
\put(12.6,1.1){$\bullet$}
\put(12.6,1.6){$\bullet$}
\put(12.6,2.1){$\bullet$}
\put(12.6,2.6){$\bullet$}
\put(13.1,.1){$\bullet$}
\put(13.6,.1){$\bullet$}
\put(14.1,.1){$\bullet$}
\put(14.6,.1){$\bullet$}
\put(15.1,.1){$7$}
\put(15.1,.6){$7$}
\put(15.1,1.1){$7$}
\put(15.1,1.6){$7$}
\put(15.6,.1){$7$}
\put(16.1,.1){$\bullet$}
\put(16.1,.6){$\bullet$}
\put(16.1,1.1){$\bullet$}
\put(16.1,1.6){$\bullet$}
\put(16.1,2.1){$\bullet$}
\put(16.1,2.6){$\bullet$}
\put(16.6,.1){$\bullet$}
\put(17.1,.1){$\bullet$}
\put(17.6,.1){$\bullet$}
\put(18.1,.1){$\bullet$}

\end{picture}
\end{center}

\setlength{\unitlength}{.7cm}
\begin{center}
\begin{picture}(20,3.5)(0,0)
\thinlines

\put(0,0){\line(1,0){8.5}}
\put(0,0){\line(0,1){3}}
\put(.5,.5){\line(1,0){2}}
\put(0,3){\line(1,0){.5}}
\put(.5,.5){\line(0,1){2.5}}

\put(2.5,0){\line(0,1){3}}
\put(2.5,1){\line(1,0){.5}}
\put(2.5,3){\line(1,0){.5}}
\put(3,.5){\line(0,1){2.5}}
\put(3,.5){\line(1,0){2}}
\put(4.5,0){\line(0,1){.5}}

\put(5,0){\line(0,1){2}}
\put(5,2){\line(1,0){.5}}
\put(5.5,.5){\line(0,1){1.5}}
\put(5.5, .5){\line(1,0){.5}}

\put(6,0){\line(0,1){3}}
\put(6,3){\line(1,0){.5}}
\put(6.5,.5){\line(0,1){2.5}}
\put(6.5,.5){\line(1,0){2}}
\put(8.5,0){\line(0,1){.5}}

\put(.1,.1){$\bullet$}
\put(.1,.6){$\bullet$}
\put(.1,1.1){$\bullet$}
\put(.1,1.6){$\bullet$}
\put(.1,2.1){$\bullet$}
\put(.1,2.6){$\bullet$}
\put(.6,.1){$\bullet$}
\put(1.1,.1){$\bullet$}
\put(1.6,.1){$\bullet$}
\put(2.1,.1){$\bullet$}
\put(2.6,.1){$\bullet$}
\put(2.6,.6){$\bullet$}
\put(2.6,1.1){$6$}
\put(2.6,1.6){$6$}
\put(2.6,2.1){$6$}
\put(2.6,2.6){$6$}
\put(3.1,.1){$\bullet$}
\put(3.6,.1){$\bullet$}
\put(4.1,.1){$\bullet$}
\put(4.6,.1){$6$}
\put(5.1,.1){$7$}
\put(5.1,.6){$7$}
\put(5.1,1.1){$7$}
\put(5.1,1.6){$7$}
\put(5.6,.1){$7$}
\put(6.1,.1){$\bullet$}
\put(6.1,.6){$\bullet$}
\put(6.1,1.1){$\bullet$}
\put(6.1,1.6){$\bullet$}
\put(6.1,2.1){$\bullet$}
\put(6.1,2.6){$\bullet$}
\put(6.6,.1){$\bullet$}
\put(7.1,.1){$\bullet$}
\put(7.6,.1){$\bullet$}
\put(8.1,.1){$\bullet$}

\put(10,0){\line(1,0){8.5}}
\put(10,0){\line(0,1){3}}
\put(10.5,.5){\line(1,0){2}}
\put(10,3){\line(1,0){.5}}
\put(10.5,.5){\line(0,1){2.5}}

\put(12.5,0){\line(0,1){3}}
\put(12.5,1){\line(1,0){.5}}
\put(12.5,3){\line(1,0){.5}}
\put(13,.5){\line(0,1){2.5}}
\put(13,.5){\line(1,0){2}}
\put(14.5,0){\line(0,1){.5}}

\put(15,0){\line(0,1){2}}
\put(15,2){\line(1,0){.5}}
\put(15.5,.5){\line(0,1){1.5}}
\put(15.5, .5){\line(1,0){.5}}

\put(16,0){\line(0,1){3}}
\put(16,3){\line(1,0){.5}}
\put(16.5,.5){\line(0,1){2.5}}
\put(16.5,.5){\line(1,0){2}}
\put(18.5,0){\line(0,1){.5}}

\put(10.1,.1){$\bullet$}
\put(10.1,.6){$\bullet$}
\put(10.1,1.1){$\bullet$}
\put(10.1,1.6){$\bullet$}
\put(10.1,2.1){$\bullet$}
\put(10.1,2.6){$\bullet$}
\put(10.6,.1){$\bullet$}
\put(11.1,.1){$\bullet$}
\put(11.6,.1){$\bullet$}
\put(12.1,.1){$\bullet$}
\put(12.6,.1){$5$}
\put(12.6,.6){$5$}
\put(12.6,1.1){$6$}
\put(12.6,1.6){$6$}
\put(12.6,2.1){$6$}
\put(12.6,2.6){$6$}
\put(13.1,.1){$5$}
\put(13.6,.1){$5$}
\put(14.1,.1){$5$}
\put(14.6,.1){$6$}
\put(15.1,.1){$7$}
\put(15.1,.6){$7$}
\put(15.1,1.1){$7$}
\put(15.1,1.6){$7$}
\put(15.6,.1){$7$}
\put(16.1,.1){$\bullet$}
\put(16.1,.6){$\bullet$}
\put(16.1,1.1){$\bullet$}
\put(16.1,1.6){$\bullet$}
\put(16.1,2.1){$\bullet$}
\put(16.1,2.6){$\bullet$}
\put(16.6,.1){$\bullet$}
\put(17.1,.1){$\bullet$}
\put(17.6,.1){$\bullet$}
\put(18.1,.1){$\bullet$}

\end{picture}
\end{center}

\setlength{\unitlength}{.7cm}
\begin{center}
\begin{picture}(20,3.5)(0,0)
\thinlines

\put(0,0){\line(1,0){8.5}}
\put(0,0){\line(0,1){3}}
\put(.5,.5){\line(1,0){2}}
\put(0,2){\line(1,0){.5}}
\put(0,3){\line(1,0){.5}}
\put(.5,.5){\line(0,1){2.5}}
\put(1,0){\line(0,1){.5}}

\put(2.5,0){\line(0,1){3}}
\put(2.5,1){\line(1,0){.5}}
\put(2.5,3){\line(1,0){.5}}
\put(3,.5){\line(0,1){2.5}}
\put(3,.5){\line(1,0){2}}
\put(4.5,0){\line(0,1){.5}}

\put(5,0){\line(0,1){2}}
\put(5,2){\line(1,0){.5}}
\put(5.5,.5){\line(0,1){1.5}}
\put(5.5, .5){\line(1,0){.5}}

\put(6,0){\line(0,1){3}}
\put(6,3){\line(1,0){.5}}
\put(6.5,.5){\line(0,1){2.5}}
\put(6.5,.5){\line(1,0){2}}
\put(8.5,0){\line(0,1){.5}}

\put(.1,.1){$\bullet$}
\put(.1,.6){$\bullet$}
\put(.1,1.1){$\bullet$}
\put(.1,1.6){$\bullet$}
\put(.1,2.1){$4$}
\put(.1,2.6){$4$}
\put(.6,.1){$\bullet$}
\put(1.1,.1){$4$}
\put(1.6,.1){$4$}
\put(2.1,.1){$4$}
\put(2.6,.1){$5$}
\put(2.6,.6){$5$}
\put(2.6,1.1){$6$}
\put(2.6,1.6){$6$}
\put(2.6,2.1){$6$}
\put(2.6,2.6){$6$}
\put(3.1,.1){$5$}
\put(3.6,.1){$5$}
\put(4.1,.1){$5$}
\put(4.6,.1){$6$}
\put(5.1,.1){$7$}
\put(5.1,.6){$7$}
\put(5.1,1.1){$7$}
\put(5.1,1.6){$7$}
\put(5.6,.1){$7$}
\put(6.1,.1){$\bullet$}
\put(6.1,.6){$\bullet$}
\put(6.1,1.1){$\bullet$}
\put(6.1,1.6){$\bullet$}
\put(6.1,2.1){$\bullet$}
\put(6.1,2.6){$\bullet$}
\put(6.6,.1){$\bullet$}
\put(7.1,.1){$\bullet$}
\put(7.6,.1){$\bullet$}
\put(8.1,.1){$\bullet$}

\put(10,0){\line(1,0){8.5}}
\put(10,0){\line(0,1){3}}
\put(10.5,.5){\line(1,0){2}}
\put(10,2){\line(1,0){.5}}
\put(10,3){\line(1,0){.5}}
\put(10.5,.5){\line(0,1){2.5}}
\put(11,0){\line(0,1){.5}}

\put(12.5,0){\line(0,1){3}}
\put(12.5,1){\line(1,0){.5}}
\put(12.5,3){\line(1,0){.5}}
\put(13,.5){\line(0,1){2.5}}
\put(13,.5){\line(1,0){2}}
\put(14.5,0){\line(0,1){.5}}

\put(15,0){\line(0,1){2}}
\put(15,2){\line(1,0){.5}}
\put(15.5,.5){\line(0,1){1.5}}
\put(15.5, .5){\line(1,0){.5}}

\put(16,0){\line(0,1){3}}
\put(16,3){\line(1,0){.5}}
\put(16.5,.5){\line(0,1){2.5}}
\put(16.5,.5){\line(1,0){2}}
\put(18.5,0){\line(0,1){.5}}

\put(10.1,.1){$3$}
\put(10.1,.6){$3$}
\put(10.1,1.1){$3$}
\put(10.1,1.6){$3$}
\put(10.1,2.1){$4$}
\put(10.1,2.6){$4$}
\put(10.6,.1){$3$}
\put(11.1,.1){$4$}
\put(11.6,.1){$4$}
\put(12.1,.1){$4$}
\put(12.6,.1){$5$}
\put(12.6,.6){$5$}
\put(12.6,1.1){$6$}
\put(12.6,1.6){$6$}
\put(12.6,2.1){$6$}
\put(12.6,2.6){$6$}
\put(13.1,.1){$5$}
\put(13.6,.1){$5$}
\put(14.1,.1){$5$}
\put(14.6,.1){$6$}
\put(15.1,.1){$7$}
\put(15.1,.6){$7$}
\put(15.1,1.1){$7$}
\put(15.1,1.6){$7$}
\put(15.6,.1){$7$}
\put(16.1,.1){$\bullet$}
\put(16.1,.6){$\bullet$}
\put(16.1,1.1){$\bullet$}
\put(16.1,1.6){$\bullet$}
\put(16.1,2.1){$\bullet$}
\put(16.1,2.6){$\bullet$}
\put(16.6,.1){$\bullet$}
\put(17.1,.1){$\bullet$}
\put(17.6,.1){$\bullet$}
\put(18.1,.1){$\bullet$}

\end{picture}
\end{center}

\setlength{\unitlength}{.7cm}
\begin{center}
\begin{picture}(20,3.5)(0,0)
\thinlines

\put(0,0){\line(1,0){8.5}}
\put(0,0){\line(0,1){3}}
\put(.5,.5){\line(1,0){2}}
\put(0,2){\line(1,0){.5}}
\put(0,3){\line(1,0){.5}}
\put(.5,.5){\line(0,1){2.5}}
\put(1,0){\line(0,1){.5}}

\put(2.5,0){\line(0,1){3}}
\put(2.5,1){\line(1,0){.5}}
\put(2.5,3){\line(1,0){.5}}
\put(3,.5){\line(0,1){2.5}}
\put(3,.5){\line(1,0){2}}
\put(4.5,0){\line(0,1){.5}}

\put(5,0){\line(0,1){2}}
\put(5,2){\line(1,0){.5}}
\put(5.5,.5){\line(0,1){1.5}}
\put(5.5, .5){\line(1,0){.5}}

\put(6,0){\line(0,1){3}}
\put(6,3){\line(1,0){.5}}
\put(6.5,.5){\line(0,1){2.5}}
\put(6,.5){\line(1,0){2.5}}
\put(8.5,0){\line(0,1){.5}}

\put(.1,.1){$3$}
\put(.1,.6){$3$}
\put(.1,1.1){$3$}
\put(.1,1.6){$3$}
\put(.1,2.1){$4$}
\put(.1,2.6){$4$}
\put(.6,.1){$3$}
\put(1.1,.1){$4$}
\put(1.6,.1){$4$}
\put(2.1,.1){$4$}
\put(2.6,.1){$5$}
\put(2.6,.6){$5$}
\put(2.6,1.1){$6$}
\put(2.6,1.6){$6$}
\put(2.6,2.1){$6$}
\put(2.6,2.6){$6$}
\put(3.1,.1){$5$}
\put(3.6,.1){$5$}
\put(4.1,.1){$5$}
\put(4.6,.1){$6$}
\put(5.1,.1){$7$}
\put(5.1,.6){$7$}
\put(5.1,1.1){$7$}
\put(5.1,1.6){$7$}
\put(5.6,.1){$7$}
\put(6.1,.1){$\bullet$}
\put(6.1,.6){$2$}
\put(6.1,1.1){$2$}
\put(6.1,1.6){$2$}
\put(6.1,2.1){$2$}
\put(6.1,2.6){$2$}
\put(6.6,.1){$\bullet$}
\put(7.1,.1){$\bullet$}
\put(7.6,.1){$\bullet$}
\put(8.1,.1){$\bullet$}

\put(10,0){\line(1,0){8.5}}
\put(10,0){\line(0,1){3}}
\put(10.5,.5){\line(1,0){2}}
\put(10,2){\line(1,0){.5}}
\put(10,3){\line(1,0){.5}}
\put(10.5,.5){\line(0,1){2.5}}
\put(11,0){\line(0,1){.5}}

\put(12.5,0){\line(0,1){3}}
\put(12.5,1){\line(1,0){.5}}
\put(12.5,3){\line(1,0){.5}}
\put(13,.5){\line(0,1){2.5}}
\put(13,.5){\line(1,0){2}}
\put(14.5,0){\line(0,1){.5}}

\put(15,0){\line(0,1){2}}
\put(15,2){\line(1,0){.5}}
\put(15.5,.5){\line(0,1){1.5}}
\put(15.5, .5){\line(1,0){.5}}

\put(16,0){\line(0,1){3}}
\put(16,3){\line(1,0){.5}}
\put(16.5,.5){\line(0,1){2.5}}
\put(16,.5){\line(1,0){2.5}}
\put(18.5,0){\line(0,1){.5}}

\put(10.1,.1){$3$}
\put(10.1,.6){$3$}
\put(10.1,1.1){$3$}
\put(10.1,1.6){$3$}
\put(10.1,2.1){$4$}
\put(10.1,2.6){$4$}
\put(10.6,.1){$3$}
\put(11.1,.1){$4$}
\put(11.6,.1){$4$}
\put(12.1,.1){$4$}
\put(12.6,.1){$5$}
\put(12.6,.6){$5$}
\put(12.6,1.1){$6$}
\put(12.6,1.6){$6$}
\put(12.6,2.1){$6$}
\put(12.6,2.6){$6$}
\put(13.1,.1){$5$}
\put(13.6,.1){$5$}
\put(14.1,.1){$5$}
\put(14.6,.1){$6$}
\put(15.1,.1){$7$}
\put(15.1,.6){$7$}
\put(15.1,1.1){$7$}
\put(15.1,1.6){$7$}
\put(15.6,.1){$7$}
\put(16.1,.1){$1$}
\put(16.1,.6){$2$}
\put(16.1,1.1){$2$}
\put(16.1,1.6){$2$}
\put(16.1,2.1){$2$}
\put(16.1,2.6){$2$}
\put(16.6,.1){$1$}
\put(17.1,.1){$1$}
\put(17.6,.1){$1$}
\put(18.1,.1){$1$}

\end{picture}
\end{center}

The last tableau above is the evacuation tableau for $P$:

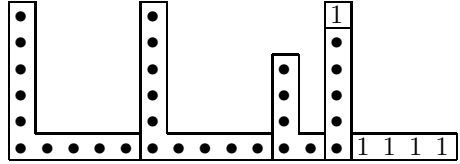
\setlength{\unitlength}{.7cm}
\begin{center}
\begin{picture}(20,3.5)(0,0)
\thinlines

\put(3.8,1){$ev(P)=$}

\put(6,0){\line(1,0){8.5}}
\put(6,0){\line(0,1){3}}
\put(6.5,.5){\line(1,0){2}}
\put(6,2){\line(1,0){.5}}
\put(6,3){\line(1,0){.5}}
\put(6.5,.5){\line(0,1){2.5}}
\put(7,0){\line(0,1){.5}}

\put(8.5,0){\line(0,1){3}}
\put(8.5,1){\line(1,0){.5}}
\put(8.5,3){\line(1,0){.5}}
\put(9,.5){\line(0,1){2.5}}
\put(9,.5){\line(1,0){2}}
\put(10.5,0){\line(0,1){.5}}

\put(11,0){\line(0,1){2}}
\put(11,2){\line(1,0){.5}}
\put(11.5,.5){\line(0,1){1.5}}
\put(11.5, .5){\line(1,0){.5}}

\put(12,0){\line(0,1){3}}
\put(12,3){\line(1,0){.5}}
\put(12.5,.5){\line(0,1){2.5}}
\put(12,.5){\line(1,0){2.5}}
\put(14.5,0){\line(0,1){.5}}

\put(6.1,.1){$3$}
\put(6.1,.6){$3$}
\put(6.1,1.1){$3$}
\put(6.1,1.6){$3$}
\put(6.1,2.1){$4$}
\put(6.1,2.6){$4$}
\put(6.6,.1){$3$}
\put(7.1,.1){$4$}
\put(7.6,.1){$4$}
\put(8.1,.1){$4$}
\put(8.6,.1){$5$}
\put(8.6,.6){$5$}
\put(8.6,1.1){$6$}
\put(8.6,1.6){$6$}
\put(8.6,2.1){$6$}
\put(8.6,2.6){$6$}
\put(9.1,.1){$5$}
\put(9.6,.1){$5$}
\put(10.1,.1){$5$}
\put(10.6,.1){$6$}
\put(11.1,.1){$7$}
\put(11.1,.6){$7$}
\put(11.1,1.1){$7$}
\put(11.1,1.6){$7$}
\put(11.6,.1){$7$}
\put(12.1,.1){$1$}
\put(12.1,.6){$2$}
\put(12.1,1.1){$2$}
\put(12.1,1.6){$2$}
\put(12.1,2.1){$2$}
\put(12.1,2.6){$2$}
\put(12.6,.1){$1$}
\put(13.1,.1){$1$}
\put(13.6,.1){$1$}
\put(14.1,.1){$1$}

\end{picture}
\end{center}

In Section \ref{theorems}, we will prove that $ev(P)=\hat{P}$, where $P$ is the tableau obtained from the the $k$-ribbon Fibonacci insertion algorithm described in Section \ref{insertion} and $\hat{P}$ is the $k$-ribbon Fibonacci path tableau obtained from Fomin's growth diagram.  As a consequence in Section \ref{theorems}, we will also see that the evacuation algorithm is a bijection between standard $k$-ribbon Fibonacci tableaux and $k$-ribbon Fibonacci path tableaux.  

To complete this section, we will describe the inverse of the evacuation algorithm.  Given a $k$-ribbon Fibonacci path tableau of shape $\lambda$, denote the column containing the $k$-ribbon labeled with $1$'s as column $c$.  Remove this $k$-ribbon from the tableau and decrease all remaining values in the tableau by 1.  If there is no $k$-ribbon in column $c$ once the $k$-ribbon containing $1$'s has been removed, then stop and place a $k$-ribbon containing $1$'s in column $c$ of an empty $k$-ribbon tableau of the same shape $\lambda$.  

If there is a $k$-ribbon of height $j$ in column $c$ once the $k$-ribbon containing $1$'s has been removed, first slide this $k$-ribbon down so that it has the shape of a single $k$-ribbon of height $j$ in column $c$, leaving an empty $k$-ribbon of height $k+1-j$ at the top of column $c$.  In an empty tableau of shape $\lambda$, place $1$'s in the position of this empty $k$-ribbon.  Next cycle the values in column $c$ and all columns to the right of $c$, leaving all orientations of $k$-ribbons fixed, in the following manner.  If $a_1 < a_2 < \cdots < a_k$ are the values in the $k$-ribbons in column $c$ and the columns to the right of $c$, then replace $a_1$ with $a_k$, $a_2$ with $a_1$, $a_3$ with $a_2$, and so on.  This creates a new $k$-ribbon Fibonacci path tableau that is one $k$-ribbon smaller than the shape $\lambda$. 

Repeat this process on the new smaller $k$-ribbon Fibonacci path tableau.  At the $i$th step of the iteration, place $i$'s into the partially empty tableau of shape $\lambda$.  Continue until $\lambda$ is a standard $k$-ribbon Fibonacci tableau containing no empty $k$-ribbons.  Interestingly, the tiling of the standard $k$-ribbon Fibonacci tableau and the tiling of the evacuation of that tableau are related by simply swapping the heights of the $k$-ribbons in the columns corresponding to a 2, though the numbers in these $k$-ribbons may be different.

We will describe the first three steps of the inverse algorithm using the evacuation tableau $ev(P)$ obtained from the earlier example.  In the first step we erase the $k$-ribbon containing $1$'s and reduce all values in the remaining tableau by 1.  Since there are no values in any columns to the right, then the inverse algorithm stops here and the first step of the inverse algorithm gives:

\setlength{\unitlength}{.7cm}
\begin{center}
\begin{picture}(20,3.5)(0,0)
\thinlines

\put(6,0){\line(1,0){8.5}}
\put(6,0){\line(0,1){3}}
\put(6.5,.5){\line(1,0){2}}
\put(6,3){\line(1,0){.5}}
\put(6.5,.5){\line(0,1){2.5}}

\put(8.5,0){\line(0,1){3}}
\put(8.5,3){\line(1,0){.5}}
\put(9,.5){\line(0,1){2.5}}
\put(9,.5){\line(1,0){2}}

\put(11,0){\line(0,1){2}}
\put(11,2){\line(1,0){.5}}
\put(11.5,.5){\line(0,1){1.5}}
\put(11.5, .5){\line(1,0){.5}}

\put(12,0){\line(0,1){3}}
\put(12,2.5){\line(1,0){.5}}
\put(12,3){\line(1,0){.5}}
\put(12.5,0){\line(0,1){3}}
\put(12.5,.5){\line(1,0){2}}
\put(14.5,0){\line(0,1){.5}}

\put(6.1,.1){$\bullet$}
\put(6.1,.6){$\bullet$}
\put(6.1,1.1){$\bullet$}
\put(6.1,1.6){$\bullet$}
\put(6.1,2.1){$\bullet$}
\put(6.1,2.6){$\bullet$}
\put(6.6,.1){$\bullet$}
\put(7.1,.1){$\bullet$}
\put(7.6,.1){$\bullet$}
\put(8.1,.1){$\bullet$}
\put(8.6,.1){$\bullet$}
\put(8.6,.6){$\bullet$}
\put(8.6,1.1){$\bullet$}
\put(8.6,1.6){$\bullet$}
\put(8.6,2.1){$\bullet$}
\put(8.6,2.6){$\bullet$}
\put(9.1,.1){$\bullet$}
\put(9.6,.1){$\bullet$}
\put(10.1,.1){$\bullet$}
\put(10.6,.1){$\bullet$}
\put(11.1,.1){$\bullet$}
\put(11.1,.6){$\bullet$}
\put(11.1,1.1){$\bullet$}
\put(11.1,1.6){$\bullet$}
\put(11.6,.1){$\bullet$}
\put(12.1,.1){$\bullet$}
\put(12.1,.6){$\bullet$}
\put(12.1,1.1){$\bullet$}
\put(12.1,1.6){$\bullet$}
\put(12.1,2.1){$\bullet$}
\put(12.1,2.6){$1$}
\put(12.6,.1){$1$}
\put(13.1,.1){$1$}
\put(13.6,.1){$1$}
\put(14.1,.1){$1$}

\end{picture}
\end{center}

Similarly, in the second step of the algorithm we erase the $1$'s, reduce all values in the tableau by 1 and since there is no $k$-ribbon now present in column $c$ we stop here and the tableau obtained from the inverse algorithm looks like:

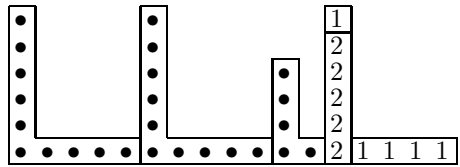
\setlength{\unitlength}{.7cm}
\begin{center}
\begin{picture}(20,3.5)(0,0)
\thinlines

\put(6,0){\line(1,0){8.5}}
\put(6,0){\line(0,1){3}}
\put(6.5,.5){\line(1,0){2}}
\put(6,3){\line(1,0){.5}}
\put(6.5,.5){\line(0,1){2.5}}

\put(8.5,0){\line(0,1){3}}
\put(8.5,3){\line(1,0){.5}}
\put(9,.5){\line(0,1){2.5}}
\put(9,.5){\line(1,0){2}}

\put(11,0){\line(0,1){2}}
\put(11,2){\line(1,0){.5}}
\put(11.5,.5){\line(0,1){1.5}}
\put(11.5, .5){\line(1,0){.5}}

\put(12,0){\line(0,1){3}}
\put(12,2.5){\line(1,0){.5}}
\put(12,3){\line(1,0){.5}}
\put(12.5,0){\line(0,1){3}}
\put(12.5,.5){\line(1,0){2}}
\put(14.5,0){\line(0,1){.5}}

\put(6.1,.1){$\bullet$}
\put(6.1,.6){$\bullet$}
\put(6.1,1.1){$\bullet$}
\put(6.1,1.6){$\bullet$}
\put(6.1,2.1){$\bullet$}
\put(6.1,2.6){$\bullet$}
\put(6.6,.1){$\bullet$}
\put(7.1,.1){$\bullet$}
\put(7.6,.1){$\bullet$}
\put(8.1,.1){$\bullet$}
\put(8.6,.1){$\bullet$}
\put(8.6,.6){$\bullet$}
\put(8.6,1.1){$\bullet$}
\put(8.6,1.6){$\bullet$}
\put(8.6,2.1){$\bullet$}
\put(8.6,2.6){$\bullet$}
\put(9.1,.1){$\bullet$}
\put(9.6,.1){$\bullet$}
\put(10.1,.1){$\bullet$}
\put(10.6,.1){$\bullet$}
\put(11.1,.1){$\bullet$}
\put(11.1,.6){$\bullet$}
\put(11.1,1.1){$\bullet$}
\put(11.1,1.6){$\bullet$}
\put(11.6,.1){$\bullet$}
\put(12.1,.1){$2$}
\put(12.1,.6){$2$}
\put(12.1,1.1){$2$}
\put(12.1,1.6){$2$}
\put(12.1,2.1){$2$}
\put(12.1,2.6){$1$}
\put(12.6,.1){$1$}
\put(13.1,.1){$1$}
\put(13.6,.1){$1$}
\put(14.1,.1){$1$}

\end{picture}
\end{center}

At the third step, the path tableau we are applying the inverse algorithm to looks like:

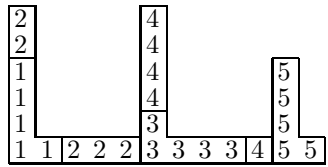
\setlength{\unitlength}{.7cm}
\begin{center}
\begin{picture}(20,3.5)(0,0)
\thinlines

\put(6,0){\line(1,0){6}}
\put(6,0){\line(0,1){3}}
\put(6.5,.5){\line(1,0){2}}
\put(6,2){\line(1,0){.5}}
\put(6,3){\line(1,0){.5}}
\put(6.5,.5){\line(0,1){2.5}}
\put(7,0){\line(0,1){.5}}

\put(8.5,0){\line(0,1){3}}
\put(8.5,1){\line(1,0){.5}}
\put(8.5,3){\line(1,0){.5}}
\put(9,.5){\line(0,1){2.5}}
\put(9,.5){\line(1,0){2}}
\put(10.5,0){\line(0,1){.5}}

\put(11,0){\line(0,1){2}}
\put(11,2){\line(1,0){.5}}
\put(11.5,.5){\line(0,1){1.5}}
\put(11.5,.5){\line(1,0){.5}}
\put(12,0){\line(0,1){.5}}

\put(6.1,.1){$1$}
\put(6.1,.6){$1$}
\put(6.1,1.1){$1$}
\put(6.1,1.6){$1$}
\put(6.1,2.1){$2$}
\put(6.1,2.6){$2$}
\put(6.6,.1){$1$}
\put(7.1,.1){$2$}
\put(7.6,.1){$2$}
\put(8.1,.1){$2$}
\put(8.6,.1){$3$}
\put(8.6,.6){$3$}
\put(8.6,1.1){$4$}
\put(8.6,1.6){$4$}
\put(8.6,2.1){$4$}
\put(8.6,2.6){$4$}
\put(9.1,.1){$3$}
\put(9.6,.1){$3$}
\put(10.1,.1){$3$}
\put(10.6,.1){$4$}
\put(11.1,.1){$5$}
\put(11.1,.6){$5$}
\put(11.1,1.1){$5$}
\put(11.1,1.6){$5$}
\put(11.6,.1){$5$}

\end{picture}
\end{center}

Following the inverse algorithm, we erase the $1$'s and reduce the each of the remaining values by 1 while preserving the shape of each $k$-ribbon to obtain:

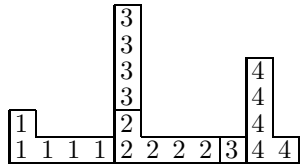
\setlength{\unitlength}{.7cm}
\begin{center}
\begin{picture}(20,3.5)(0,0)
\thinlines

\put(6.5,0){\line(1,0){5.5}}
\put(6.5,0){\line(0,1){1}}
\put(7,.5){\line(1,0){1.5}}
\put(6.5,1){\line(1,0){.5}}
\put(7,.5){\line(0,1){.5}}

\put(8.5,0){\line(0,1){3}}
\put(8.5,1){\line(1,0){.5}}
\put(8.5,3){\line(1,0){.5}}
\put(9,.5){\line(0,1){2.5}}
\put(9,.5){\line(1,0){2}}
\put(10.5,0){\line(0,1){.5}}

\put(11,0){\line(0,1){2}}
\put(11,2){\line(1,0){.5}}
\put(11.5,.5){\line(0,1){1.5}}
\put(11.5,.5){\line(1,0){.5}}
\put(12,0){\line(0,1){.5}}

\put(6.6,.1){$1$}
\put(6.6,.6){$1$}
\put(7.1,.1){$1$}
\put(7.6,.1){$1$}
\put(8.1,.1){$1$}
\put(8.6,.1){$2$}
\put(8.6,.6){$2$}
\put(8.6,1.1){$3$}
\put(8.6,1.6){$3$}
\put(8.6,2.1){$3$}
\put(8.6,2.6){$3$}
\put(9.1,.1){$2$}
\put(9.6,.1){$2$}
\put(10.1,.1){$2$}
\put(10.6,.1){$3$}
\put(11.1,.1){$4$}
\put(11.1,.6){$4$}
\put(11.1,1.1){$4$}
\put(11.1,1.6){$4$}
\put(11.6,.1){$4$}

\end{picture}
\end{center}

Now we cycle the values as described to obtain a new path tableau (shown on the left) and the tableau corresponding the third step of the inverse algorithm (shown on the right): 

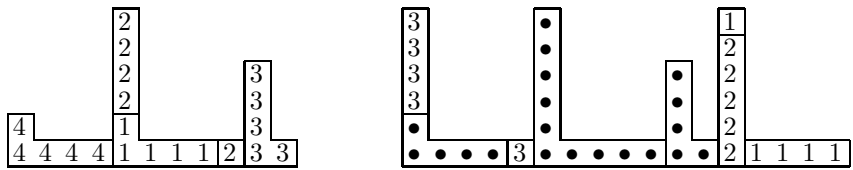
\setlength{\unitlength}{.7cm}
\begin{center}
\begin{picture}(20,3.5)(0,0)
\thinlines

\put(.5,0){\line(1,0){5.5}}
\put(.5,0){\line(0,1){1}}
\put(1,.5){\line(1,0){1.5}}
\put(.5,1){\line(1,0){.5}}
\put(1,.5){\line(0,1){.5}}

\put(2.5,0){\line(0,1){3}}
\put(2.5,1){\line(1,0){.5}}
\put(2.5,3){\line(1,0){.5}}
\put(3,.5){\line(0,1){2.5}}
\put(3,.5){\line(1,0){2}}
\put(4.5,0){\line(0,1){.5}}

\put(5,0){\line(0,1){2}}
\put(5,2){\line(1,0){.5}}
\put(5.5,.5){\line(0,1){1.5}}
\put(5.5,.5){\line(1,0){.5}}
\put(6,0){\line(0,1){.5}}

\put(.6,.1){$4$}
\put(.6,.6){$4$}
\put(1.1,.1){$4$}
\put(1.6,.1){$4$}
\put(2.1,.1){$4$}
\put(2.6,.1){$1$}
\put(2.6,.6){$1$}
\put(2.6,1.1){$2$}
\put(2.6,1.6){$2$}
\put(2.6,2.1){$2$}
\put(2.6,2.6){$2$}
\put(3.1,.1){$1$}
\put(3.6,.1){$1$}
\put(4.1,.1){$1$}
\put(4.6,.1){$2$}
\put(5.1,.1){$3$}
\put(5.1,.6){$3$}
\put(5.1,1.1){$3$}
\put(5.1,1.6){$3$}
\put(5.6,.1){$3$}

\put(8,0){\line(1,0){8.5}}
\put(8,0){\line(0,1){3}}
\put(8,1){\line(1,0){.5}}
\put(8.5,.5){\line(1,0){2}}
\put(8,3){\line(1,0){.5}}
\put(8.5,.5){\line(0,1){2.5}}
\put(10,0){\line(0,1){.5}}

\put(10.5,0){\line(0,1){3}}
\put(10.5,3){\line(1,0){.5}}
\put(11,.5){\line(0,1){2.5}}
\put(11,.5){\line(1,0){2}}

\put(13,0){\line(0,1){2}}
\put(13,2){\line(1,0){.5}}
\put(13.5,.5){\line(0,1){1.5}}
\put(13.5, .5){\line(1,0){.5}}

\put(14,0){\line(0,1){3}}
\put(14,2.5){\line(1,0){.5}}
\put(14,3){\line(1,0){.5}}
\put(14.5,0){\line(0,1){3}}
\put(14.5,.5){\line(1,0){2}}
\put(16.5,0){\line(0,1){.5}}

\put(8.1,.1){$\bullet$}
\put(8.1,.6){$\bullet$}
\put(8.1,1.1){$3$}
\put(8.1,1.6){$3$}
\put(8.1,2.1){$3$}
\put(8.1,2.6){$3$}
\put(8.6,.1){$\bullet$}
\put(9.1,.1){$\bullet$}
\put(9.6,.1){$\bullet$}
\put(10.1,.1){$3$}
\put(10.6,.1){$\bullet$}
\put(10.6,.6){$\bullet$}
\put(10.6,1.1){$\bullet$}
\put(10.6,1.6){$\bullet$}
\put(10.6,2.1){$\bullet$}
\put(10.6,2.6){$\bullet$}
\put(11.1,.1){$\bullet$}
\put(11.6,.1){$\bullet$}
\put(12.1,.1){$\bullet$}
\put(12.6,.1){$\bullet$}
\put(13.1,.1){$\bullet$}
\put(13.1,.6){$\bullet$}
\put(13.1,1.1){$\bullet$}
\put(13.1,1.6){$\bullet$}
\put(13.6,.1){$\bullet$}
\put(14.1,.1){$2$}
\put(14.1,.6){$2$}
\put(14.1,1.1){$2$}
\put(14.1,1.6){$2$}
\put(14.1,2.1){$2$}
\put(14.1,2.6){$1$}
\put(14.6,.1){$1$}
\put(15.1,.1){$1$}
\put(15.6,.1){$1$}
\put(16.1,.1){$1$}

\end{picture}
\end{center}

\section{A Geometric Interpretation for $k$-ribbon Fibonacci Tableaux}
\label{geometric}

In  \cite{CamKill}, \cite{Kil}, there is a description of shadow lines for the square diagram of a colored permutation using 1 or 2 colors which can be used to determine the standard Fibonacci tableau obtained through the insertion algorithm.  This description easily generalizes to $k$-colored permutations and the $k$-ribbon insertion algorithm.  In this section, we include that description and prove that it produces the $P$ tableau obtained by the $k$-ribbon insertion algorithm presented in Section 4.

Starting with the square diagram of a $k$-colored permutation $\pi$ (see Section 2), we draw shadow lines to determine its associated $k$-ribbon tableau.  To draw the shadow lines, $L_1,L_2\ldots$, for the square diagram of $\pi$, start at the top row and draw a broken line $L_1$ through the $X^{j_1}$ in the top row and the $X^{l_1}$ in the rightmost column.  The second broken line $L_2$ will be drawn through the row containing the highest $X^{j_2}$ not already on a line and the rightmost column containing an $X^{l_2}$ not already on a line.  Continue in this manner until there are no more X's available.  For example, for the $5$-colored permutation $\pi = 2^37^11^15^46^34^23^4$, the lines look like:

\setlength{\unitlength}{1cm}
\begin{center}
\begin{picture}(7,6.5)(1,1)
\thinlines

\put(1,0){\line(0,1){7}}
\put(2,0){\line(0,1){7}}
\put(3,0){\line(0,1){7}}
\put(4,0){\line(0,1){7}}
\put(5,0){\line(0,1){7}}
\put(6,0){\line(0,1){7}}
\put(7,0){\line(0,1){7}}
\put(8,0){\line(0,1){7}}
\put(1,0){\line(1,0){7}}
\put(1,1){\line(1,0){7}}
\put(1,2){\line(1,0){7}}
\put(1,3){\line(1,0){7}}
\put(1,4){\line(1,0){7}}
\put(1,5){\line(1,0){7}}
\put(1,6){\line(1,0){7}}
\put(1,7){\line(1,0){7}}
\put(1.3, 1.3){$X^3$}
\put(2.3, 6.3){$X^1$}
\put(3.3, .3){$X^1$}
\put(4.3, 4.3){$X^4$}
\put(5.3, 5.3){$X^3$}
\put(6.3, 3.3){$X^2$}
\put(7.3, 2.3){$X^4$}

\dashline{.1}(1,6.4)(7.5,6.4) 
\dashline{.1}(7.5,6.4)(7.5,0.5) 
\put(0,6.4){$L_1$}

\dashline{.1}(1,5.4)(6.5,5.4) 
\dashline{.1}(6.5,5.4)(6.5,0.5) 
\put(0,5.4){$L_2$}

\dashline{.1}(1,4.4)(4.5,4.4) 
\dashline{.1}(4.5,4.4)(4.5,0.5) 
\put(0,4.4){$L_3$}

\dashline{.1}(1,1.4)(3.5,1.4) 
\dashline{.1}(3.5,1.4)(3.5,0.5) 
\put(0,1.4){$L_4$}

\end{picture}
\end{center}
\vspace{.5in}

\begin{Theorem}\label{Ptableau}
The shadow lines $L_1,L_2,\ldots$ in the square diagram of a $k$-colored permutation $\pi$ have the following relationship with the $k$-ribbon Fibonacci path tableau $P$ obtained by the insertion algorithm:
\begin{enumerate}
\item The row numbers of the $X$'s on each $L_i$ give the numbers in the $k$-ribbons in the $i$th column of $P$.  
\item If there is a single $X^j$ on the line $L_i$, then the $k$-ribbon in the $i$th column of $P$ is a single $k$-ribbon of height $j$.  
\item If there are two $X$'s on the line $L_i$, say $X^j$ and $X^l$ (with $X^j$ in a higher row than $X^l$), then the smaller row number is contained in the $k$-ribbon of height $l$ stacked on top of a $k$-ribbon of height $k+1-l$.
\end{enumerate}
\end{Theorem}

For the example above, the shadow lines give the $P$ tableau:

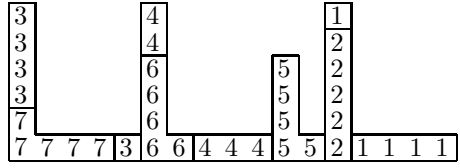
\setlength{\unitlength}{.7cm}
\begin{center}
\begin{picture}(9,3.5)(0,0)
\thinlines

\put(0,0){\line(1,0){8.5}}
\put(0,0){\line(0,1){3}}
\put(.5,.5){\line(1,0){2}}
\put(0,1){\line(1,0){.5}}
\put(0,3){\line(1,0){.5}}
\put(.5,.5){\line(0,1){2.5}}
\put(2,0){\line(0,1){.5}}
\put(2.5,0){\line(0,1){3}}
\put(2.5,2){\line(1,0){.5}}
\put(2.5,3){\line(1,0){.5}}
\put(3,.5){\line(0,1){2.5}}
\put(3,.5){\line(1,0){2}}
\put(3.5,0){\line(0,1){.5}}
\put(5,0){\line(0,1){2}}
\put(5,2){\line(1,0){.5}}
\put(5.5,.5){\line(0,1){1.5}}
\put(5.5,.5){\line(1,0){.5}}
\put(6,0){\line(0,1){3}}
\put(6,2.5){\line(1,0){.5}}
\put(6,3){\line(1,0){.5}}
\put(6.5,0){\line(0,1){3}}
\put(6.5,.5){\line(1,0){2}}
\put(8.5,0){\line(0,1){.5}}

\put(.1,.1){$7$}
\put(.1,.6){$7$}
\put(.1,1.1){$3$}
\put(.1,1.6){$3$}
\put(.1,2.1){$3$}
\put(.1,2.6){$3$}
\put(.6,.1){$7$}
\put(1.1,.1){$7$}
\put(1.6,.1){$7$}
\put(2.1,.1){$3$}
\put(2.6,.1){$6$}
\put(2.6,.6){$6$}
\put(2.6,1.1){$6$}
\put(2.6,1.6){$6$}
\put(2.6,2.1){$4$}
\put(2.6,2.6){$4$}
\put(3.1,.1){$6$}
\put(3.6,.1){$4$}
\put(4.1,.1){$4$}
\put(4.6,.1){$4$}
\put(5.1,.1){$5$}
\put(5.1,.6){$5$}
\put(5.1,1.1){$5$}
\put(5.1,1.6){$5$}
\put(5.6,.1){$5$}
\put(6.1,.1){$2$}
\put(6.1,.6){$2$}
\put(6.1,1.1){$2$}
\put(6.1,1.6){$2$}
\put(6.1,2.1){$2$}
\put(6.1,2.6){$1$}
\put(6.6,.1){$1$}
\put(7.1,.1){$1$}
\put(7.6,.1){$1$}
\put(8.1,.1){$1$}

\end{picture}
\end{center}

\begin{proof}
We will prove this result by induction on the size $n$ of the $k$-colored permutation $\pi$.  If $n=1$ then 
$\pi=1^j$ for $1 \leq j \leq k$  In this case, the shadow lines and the insertion algorithm both produce a $P$ tableau consisting of a single $k$-ribbon of height $j$ containing $1$'s.

If $n=2$ then there are $2k^2$ $k$-colored permutations.  It is not hard to check that if $\pi=1^i2^j$, then the shadow lines and insertion produce a tableau with two single $k$-ribbons of heights $j$ and $i$, containing $2$'s and $1$'s respectively, located in adjacent columns. On the other hand, if $\pi=2^i1^j$, then both the shadow lines and insertion produce a tableau consisting of a $k$-ribbon of height $j$ containing $1$'s stacked on top of a $k$-ribbon of height $k+1-j$ containing $2$'s.

Now assume that the tableau determined by the shadow lines for the $k$-colored permutation $\sigma$ of size $m$ with $m<n$ is equal to the insertion tableau $P(\sigma)$.  Let $\pi$ be a $k$-colored permutation of $n$.  Represent the permutation $\pi$ with a square diagram and draw $L_1$.  
\vspace{.2in}
 
\text{\underline{\it{Case 1:}}  }  If there is an $X^j$ in the upper right corner of the square diagram, then $L_1$ only passes through one X.  Since an $X^j$ in the upper right corner implies that $n^j$ is the last number in the permutation $\pi$, we can write $\pi = \pi_{n-1} n^j$, where $\pi_{n-1}$ represents the first $n-1$ digits in the $k$-colored permutation $\pi$.  Since $n^j$ is the last number in the permutation, when we apply the insertion algorithm to $\pi$, $n^j$ is the last number inserted into the tableau.  Thus the insertion tableau $P$ consists of a $k$-ribbon of height $j$ containing $n$'s followed by $P_{n-1}$, where $P_{n-1}$ is the insertion tableaux for $\pi_{n-1}$.  Thus the fact that the line $L_1$ drawn in the $n$th row and $n$th column only passes through one $X^j$ corresponds to the fact that there is only one $k$-ribbon at the beginning of the $P$ tableau and the shape of that $k$-ribbon is determined by the color $j$.  Then the insertion tableau $P$ and the tableau obtained from the shadow lines agree in the first column and by induction, they agree in the remaining positions.
\vspace{.1in}

\text{\underline{\it{Case 2:}}  }  If there is no X in the upper right square, then $L_1$ passes through two X's, an $X^i$ in row $n$ and an $X^j$ in column $n$ and row $a$ (counting from the bottom) with $a < n$.  This means that $a^j$, where $j$ is some number between 1 and $k$ is the last element in the permutation $\pi$, and $a^j$ is the last element inserted into the $P$ tableau.  Due to the method of insertion, the element $n$, which corresponds to the $X^i$ in the top row, is always in a $k$-ribbon located in the lower left position of $P$.  Thus when $a^j$ is inserted into the tableau, it is inserted as a $k$-ribbon of height $j$ above the $k$-ribbon containing $n$'s, possibly bumping an element $b$ to the second column.  The resulting $P$ tableau has a $k$-ribbon of height $j$ containing $a$'s on top of a $k$-ribbon of height $k+1-j$ containing $n$'s in the first column, corresponding to the fact that $L_1$ passes through two X's, one in row $n$ and one in row $a$.  

It remains to show that the rest of the $P$ tableau can be determined by removing the $n$th row and the $n$th column from the square diagram, since these elements are in the first column of $P$, and applying the inductive hypothesis to the remaining diagram.  Let the permutation $\pi$ be written as
\[
\begin{array}{ccc}
\pi&=&\begin{array}{ccccccccc}
1&2&\cdots&i-1&i&i+1&\cdots&n-1&n\\
x_1^{j_1}&x_2^{j_2}&\cdots&x_{i-1}^{j_{i-1}}&n^{i}&x_{i+1}^{j_{i+1}}&\cdots&x_{n-1}^{j_{n-1}}&a^j
\end{array}
\end{array}
\]

Recall that $P_i$ is the insertion tableau of the first $i$ elements $x_1^{j_1}x_2^{j_2}\cdots x_{i-1}^{j_{i-1}}n^{i}$.  By definition of the insertion algorithm, $P_i$ contains a $k$-ribbon of height $i$ containing $n$'s followed by $P_{i-1}$.  Since $x_{i+1} < n$, we know that $P_{i+1}$ has a $k$-ribbon of height $j_{i+1}$ containing $x_{i+1}$ on top of the $k$-ribbon containing $n$'s, now of shape $k+1-j_{i+1}$, followed by $P_{i-1}$.

When $x_{i+2}^{j_{i+2}}$ is inserted into $P_{i+1}$, $x_{i+1}$ is bumped out of the first stack of $k$-ribbons and inserted into the tableau to the right, which is $P_{i-1}$, preserving its shape.  When $x_{i+3}^{j_{i+3}}$ is inserted, $x_{i+2}$ is bumped out of the first stack of $k$-ribbons and inserted into the tableau to the right.  At the last step, $a^j$ bumps $x_{n-1}$ from the first stack of $k$-ribbons and $x_{n-1}^{j_{n-1}}$ is then inserted into the tableau to the right.  In the insertion algorithm, the shape of each stack of $k$-ribbons is determined by the color of the element in the top $k$-ribbon.  So the resulting tableau is the same as the tableau obtained by stacking the $k$-ribbon of height $j$ containing $a$ on top of the $k$-ribbon containing $n$'s followed by the tableau obtained from the insertion of 
\[
\begin{array}{ccc}
\sigma&=&\begin{array}{ccccccc}
1&2&\cdots&i-1&i&\cdots&n-2\\
x_1^{j_1}&x_2^{j_2}&\cdots&x_{i-1}^{j_{i-1}}&x_{i+1}^{j_{x+1}}&\cdots&x_{n-1}^{j_{n-1}}\end{array},
\end{array}
\]
(i.e., the permutation of size $n-2$ obtained by removing $n^{j_i}$ and $a^j$ from $\pi$).  The square diagram for $\sigma$ is the same as the square diagram for $\pi$ with the top row and rightmost column removed and any empty rows and columns removed (since empty rows and empty columns do not affect the growth diagram).  Inductively, we can now apply the above conditions to this new square diagram and continue to determine the complete insertion tableau $P(\pi)$.  
\end{proof}

\section{A Connection Between $(P,Q)$ and $(\hat{P}, \hat{Q})$}
\label{theorems}

In the previous sections, we have seen methods for obtaining a pair of $k$-ribbon Fibonacci tableaux from a $k$-colored permutation.  Given a $k$-colored permutation $\pi$, the insertion method of Section \ref{insertion} produces a pair $(P,Q)$, where $P$ is a standard $k$-ribbon Fibonacci tableau and $Q$ is a $k$-ribbon Fibonacci path tableau.  Fomin's growth diagram described in Section \ref{background} produces a pair $(\hat{P},\hat{Q})$, where $\hat{P}$ and $\hat{Q}$ are both $k$-ribbon Fibonacci path tableaux.  In this section, we prove that the evacuation method described in Section \ref{evacuation} is a bijection between standard $k$-ribbon Fibonacci tableaux and $k$-ribbon Fibonacci path tableaux by showing that $ev(P(\pi))=\hat{P}(\pi)$ while $Q(\pi)=\hat{Q}(\pi).$ 
\begin{Theorem}
For $\pi$ a $k$-colored permutation of length $n$,
$ev(P(\pi)) = \hat{P}(\pi)$.
\end{Theorem}

\begin{proof}
We will prove that $ev(P(\pi))=\hat{P}(\pi)$ by induction.  If the length of $\pi$ is 1, then the path tableau $\hat{P}$ is a single $k$-ribbon of some height $j_1$ and the insertion tableau $P$ is the same, so $\hat{P}(1^{j_1}) = ev(P(1^{j_1}))$.

Assume that for $\sigma$ a $k$-colored permutation of length $l$ with $l<n$, $ev(P(\sigma)) = \hat{P}(\sigma)$ and let $\pi$ be a $k$-colored permutation of length $n$.  
\vspace{.2in}

\text{\underline{\it{Case 1:}}  }  Suppose the square in the uppermost, rightmost corner of the square diagram for $\pi$ contains an $X^{j_1}$.

An $X^{j_1}$ in this square implies that $n^{j_1}$ is the last element in the permutation $\pi$, so $\pi = \pi_{n-1} n^{j_1}$ where $\pi_{n-1}$ represents the first $n-1$ digits in the $k$-colored permutation $\pi$.  From the square diagram, we have that $\hat{P} = n^{j_1} \hat{P}_{n-1}$ where $\hat{P}_{n-1}$ is the path tableau of shape $\nu$ obtained from $\pi_{n-1}$.  Since $n^{j_1}$ is the last number in the permutation $\pi$, when we apply the insertion algorithm, $n^{j_1}$ is the last number inserted into the tableau.  Thus the insertion tableau $P$ is a $k$-ribbon of height $j_1$ containing $n$'s followed by $P_{n-1}$ where $P_{n-1}$ is the insertion tableaux for $\pi_{n-1}$.  Following the evacuation procedure, the $k$-ribbon containing $n$'s is simply removed from $P$ and $ev(P)$ is a $k$-ribbon of height $j_1$ containing $n$'s followed by $ev(P_{n-1})$.  Since $\pi = \pi_{n-1} n^{j_1}$ where $\pi_{n-1}$ is a $k$-colored permutation of length $n-1$ and $\hat{P}_{n-1}$ is the path tableau obtained from $\pi_{n-1}$, then the inductive hypothesis implies that $\hat{P}_{n-1} = ev(P_{n-1})$.  Thus
\[
ev(P)= \hat{P}.
\]   
\vspace{.2in}

\text{\underline{\it{Case 2:}}  }  Suppose the $X^{j_n}$ in the $n$th column of the square diagram is in row $n-1$.  In this case, the permutation $\pi$ looks like:
\[
\begin{array}{ccc}
\pi&=&\begin{array}{cccccccc}
1&2&\cdots&i&i+1&\cdots&n-1&n\\
x_1^{j_1}&x_2^{j_2}&\cdots&n^{j_i}&x_{i+1}^{j_{i+1}}&\cdots&x_{n-1}^{j_{n-1}}&n-1^{j_n}.
\end{array}
\end{array}
\]
The top two squares in the last column of the growth diagram look like the following:

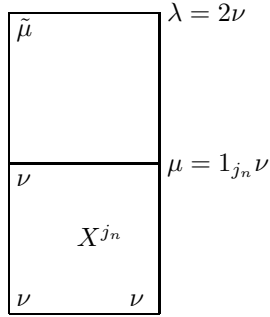
\setlength{\unitlength}{1cm}
\begin{center}
\begin{picture}(3,4)(0,0)
\thinlines

\put(0,0){\line(0,1){4}}
\put(2,0){\line(0,1){4}}
\put(0,0){\line(1,0){2}}
\put(0,2){\line(1,0){2}}
\put(0,4){\line(1,0){2}}
\put(.1,.1){$\nu$}
\put(.1, 1.7){$\nu$}
\put(.1, 3.7){$\tilde{\mu}$}
\put(1.6, .1){$\nu$}
\put(.9, .9){$X^{j_n}$}
\put(2.1, 1.9){$\mu = 1_{j_n} \nu$}
\put(2.1, 3.9){$\lambda = 2 \nu$}

\end{picture}
\end{center}

Here $\mu$ and $\lambda$ differ by a $k$-ribbon of height $k+1-j_n$ in the initial column of height 2 and $\mu$ and $\nu$ differ by a single $k$-ribbon of height $j_n$ in the initial column of height 1.  Thus the first column of $\hat{P}$ is a $k$-ribbon of height $k+1-j_n$ containing $n$'s on top of a $k$-ribbon of height $j_n$ containing $n-1$'s.  

 The first $n-2$ rows of the square diagram for $\pi$ have columns $i$ and $n$ empty, where $i$ is the column containing $X^{j_i}$ in the $n$th row of the square diagram, and these first $n-2$ rows are the growth diagram for 
\[
\begin{array}{ccc}
\sigma&=&\begin{array}{ccccccc}
1&2&\cdots&i-1&i&\cdots&n-2\\
x_1^{j_1}&x_2^{j_2}&\cdots&x_{i-1}^{j_{i-1}}&x_{i+1}^{j_{i+1}}&\cdots&x_{n-1}^{j_{n-1}}
\end{array}
\end{array}
\]
once empty columns have been removed.  Note that $\sigma$ is a $k$-colored permutation of length ${n-2}$.

By Theorem \ref{Ptableau}, the insertion tableau for $\pi$ can be determined by the shadow lines of the square diagram.  Since there is no $X^{j_l}$ in the upper right corner, the $X^{j_i}$ in the uppermost row is paired with the $X^{j_n}$ in row $n-1$ of the $n$th column.  Thus, the insertion tableau $P$ begins with a column of height two with a $k$-ribbon of height $j_n$ containing $n-1$'s on top of a $k$-ribbon of height $k+1-j_n$ containing $n$'s.  When $P$ is evacuated, the shape of the top $k$-ribbon is preserved, leaving an empty $k$-ribbon of height $k+1-j_n$ at the top of the column.  Thus the initial column of height 2 of $ev(P)$ has a $k$-ribbon of height $k+1-j_n$ containing $n$'s on top of a $k$-ribbon of height $j_n$, which is the same as the placement of the $k$-ribbon containing $n$'s in $\hat{P}$.   

At the second step of the evacuation process, the $k$-ribbon containing $n-1$'s is removed from the first column, leaving an empty $k$-ribbon of height $j_n$.  Thus $ev(P)$ has as its first column a $k$-ribbon of height $k+1-j_n$ containing $n$'s on top of a $k$-ribbon of height $j_n$ containing $n-1$'s.
Comparing $\hat{P}$ and $ev(P)$ we can see that they agree in the first column of height two.
 As shown in the proof of Theorem \ref{Ptableau}, $P(\pi)$ has a column of height 2 followed by $P(\sigma)$ where $\sigma$ is as given above.  
Since $\sigma \in S_{n-2}$, we can use our inductive hypothesis to obtain $ev(P(\sigma))=\hat{P}(\sigma)$, thus
\[
ev(P(\pi))= \hat{P}(\pi).
\]
\vspace{.2in}

\text{\underline{\it{Case 3:}}  }  Suppose the $X^{j_n}$ in column $n$ is in row $a_1 < n-1$.  In this case, $\pi$ is given by:
\[
\begin{array}{ccc}
\pi&=&\begin{array}{ccccccc}
1&2&\cdots&i&\cdots&n-1&n\\
x_1^{j_1}&x_2^{j_2}&\cdots&n^{j_i}&\cdots&x_{n-1}^{j_{n-1}}&a_1^{j_n}.
\end{array}
\end{array}
\]
The top two squares in the rightmost column of the growth diagram look like:

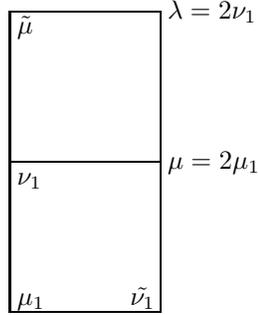
\setlength{\unitlength}{1cm}
\begin{center}
\begin{picture}(2,4)(0,0)
\thinlines
\put(0,0){\line(0,1){4}}
\put(2,0){\line(0,1){4}}
\put(0,0){\line(1,0){2}}
\put(0,2){\line(1,0){2}}
\put(0,4){\line(1,0){2}}
\put(.1,.1){$\mu_1$}
\put(.1, 1.7){$\nu_1$}
\put(.1, 3.7){$\tilde{\mu}$}
\put(1.6, .1){$\tilde{\nu_1}$}
\put(2.1, 1.9){$\mu = 2 \mu_1$}
\put(2.1, 3.9){$\lambda = 2 \nu_1$}
\end{picture}
\end{center}

Since $\lambda = 2 \nu_1$ and $\mu=2 \mu_1$, then $\lambda$ and $\mu$ differ by the same $k$-ribbon as $\nu_1$ and $\mu_1$.  If we remove the upper row and rightmost column, as well as any empty rows and columns, then the partial growth diagram of the new upper right square looks like 

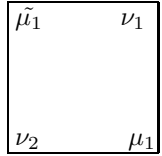
\setlength{\unitlength}{1cm}
\begin{center}
\begin{picture}(2,2)(0,0)
\thinlines
\put(0,0){\line(0,1){2}}
\put(2,0){\line(0,1){2}}
\put(0,0){\line(1,0){2}}
\put(0,2){\line(1,0){2}}
\put(.1, .1){$\nu_2$}
\put(.1, 1.7){$\tilde{\mu_1}$}
\put(1.6, .1){$\mu_1$}
\put(1.5, 1.7){$\nu_1$}
\end{picture}
\end{center}

As before, if there is an $X^{j_{n-1}}$ in the new upper right square, then $\nu_1 = 1_{j_{n-1}} \mu_1$.
If there is an $X^{j_{n-1}}$ in the square below this one, then $\nu_1$ and $\mu_1$ differ by a $k$-ribbon of height $k+1-j_{n-1}$ in the initial column of height 2. 

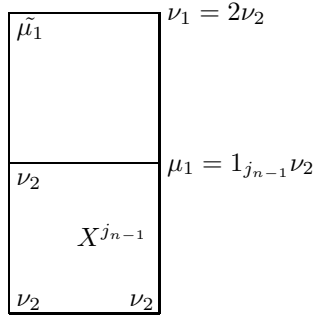
\setlength{\unitlength}{1cm}
\begin{center}
\begin{picture}(3,4)(0,0)
\thinlines

\put(0,0){\line(0,1){4}}
\put(2,0){\line(0,1){4}}
\put(0,0){\line(1,0){2}}
\put(0,2){\line(1,0){2}}
\put(0,4){\line(1,0){2}}
\put(.1,.1){$\nu_2$}
\put(.1, 1.7){$\nu_2$}
\put(1.6, .1){$\nu_2$}
\put(.1, 3.7){$\tilde{\mu_1}$}
\put(.9, .9){$X^{j_{n-1}}$}
\put(2.1, 1.9){$\mu_1 = 1_{j_{n-1}} \nu_2$}
\put(2.1, 3.9){$\nu_1 = 2 \nu_2$}

\end{picture}
\end{center}

If there is no X in either square, then the growth diagram looks like

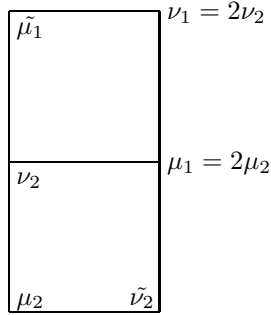
\setlength{\unitlength}{1cm}
\begin{center}
\begin{picture}(3,4)(0,0)
\thinlines
\put(0,0){\line(0,1){4}}
\put(2,0){\line(0,1){4}}
\put(0,0){\line(1,0){2}}
\put(0,2){\line(1,0){2}}
\put(0,4){\line(1,0){2}}
\put(.1,.1){$\mu_2$}
\put(.1, 1.7){$\nu_2$}
\put(1.6, .1){$\tilde{\nu_2}$}
\put(.1, 3.7){$\tilde{\mu_1}$}
\put(2.1, 1.9){$\mu_1 = 2 \mu_2$}
\put(2.1, 3.9){$\nu_1 = 2 \nu_2$}
\end{picture}
\end{center}

We can continue this procedure until $\mu_i$ and $\nu_i$ differ by a $k$-ribbon in the first column which implies that $\lambda$ and $\mu$ differ by a $k$-ribbon in the $(i+1)$st column.  If there is an $X^{j_{n-i}}$ in the upper right square then the $(i+1)$st column of $\hat{P}$ has a single $k$-ribbon of height $j_{n-i}$ containing $n$'s.  If there is an $X^{j_{n-i}}$ in the square below the upper right square at this step, then the $(i+1)$st column of $\hat{P}$ has height 2 and has a top $k$-ribbon of height $k+1-j_{n-i}$ containing $n$'s.  

We now show that the evacuation tableau $ev(P)$ has a $k$-ribbon containing $n$'s in the same place in the tableau as $\hat{P}$.  If there is not an $X^{j_n}$ in the $n$th row or $(n-1)$st row of the $n$th column of the growth diagram, then by Theorem \ref{Ptableau} the first column of $P$ has a $k$-ribbon of height $j_n$ containing $a_1$'s, with $a_1 < n-1$, on top of a $k$-ribbon of height $k+1-j_n$ containing $n$'s.  

After removing the $n$th row and the $n$th column and any empty rows and columns from the growth diagram, if there is not an $X^{j_{n-1}}$ in one of the top two rows of the rightmost column of the new growth diagram, then the second column of $P$ is a $k$-ribbon of height $j_{n-1}$ containing $a_2$'s, with $a_2 < n-2$, on top of a $k$-ribbon of height $k+1-j_{n-1}$ containing $n-1$'s.  We can continue in this manner until one of two things happens.
\vspace{.1in}

{\it{Subcase a:}}  Suppose after $i$ iterations of this process, there is an $X^{j_{n-i}}$ in the uppermost corner of the growth diagram.  In this case, the insertion tableau $P$ has $i$ columns of height $2$ followed by a $k$-ribbon of height $j_{n-i}$.  These first $i+1$ columns look like 

\begin{align*}\label{columnheight1}
&\begin{array}{cccccc}
a_1&a_2&a_3&\cdots&a_i \\n&n-1&n-2&\cdots&n-(i-1)&n-i\end{array}
\end{align*}
with $a_1 < n$, $a_2<n-1$, $a_3<n-2$, $\dots$, $a_i<n-(i-1)$ and where the column $\shortstack{$a_k$\\$n-(k-1)$}$ represents a $k$-ribbon of height $j_{n-k+1}$ containing $a_k$'s on top of a $k$-ribbon of height $k+1-j_{n-k+1}$ containing $n-(k-1)$'s.  
At the first step of evacuation for $P$, the $k$-ribbon containing $n-1$ slides one column to the left into the empty $k$-ribbon evacuated by $n$, $n-2$ slides one column to the left, and so on until $n-i$ slides one column to the left and the evacuation process terminates with an empty single $k$-ribbon of height $j_{n-i}$ in column $i+1$.  Thus $ev(P)$ has $n$'s in the single $k$-ribbon of height $j_{n-i}$ in column $i+1$, the same as $\hat{P}$, and after one step of the evacuation procedure the first $i$ columns of the $P$ tableau look like:
\[
\begin{array}{ccccc}
a_1&a_2&a_3&\cdots&a_i\\
n-1&n-2&n-3&\cdots&n-i
\end{array}.
\]
where again $\shortstack{$a_i$\\$n-i$}$ represents a stack of two $k$-ribbons and the shape of the $k$-ribbons are again determined by the top $k$-ribbon containing $a_i$'s.
The rest of the $P$ tableau remains unchanged by the evacuation procedure.
\vspace{.1in}

{\it{Subcase b:}}  Suppose after $i$ iterations of this process there is an $X^{j_{n-i}}$ in the second row from the top.  In this case, the first $i+1$ columns of the insertion tableau $P$ have height 2.  These first $i+1$ columns look like 
\begin{align*}\label{columnheight2}
&\begin{array}{cccccc}
a_1&a_2&a_3&\cdots&a_i&n-(i+1)\\n&n-1&n-2&\cdots&n-(i+1)&n-i\end{array}
\end{align*}
with $a_1<n-1$, $a_2<n-2$, $\dots$, $a_i<n-i$, where columns represent stacks of $k$-ribbons as in Subcase a.  In the evacuation process, the $k$-ribbons containing $n-1$ through $n-i$ all move one column to the left with the shape of the column determined by the top $k$-ribbon.  The $k$-ribbon containing $n-(i+1)$ becomes a single $k$-ribbon of height $j_{n-i}$ in column $i+1$.  This leaves an empty split $k$-ribbon of height $k+1-j_{n-i}$ in column $i+1$.  Thus $ev(P)$ has a $k$-ribbon of height $k+1-j_{n-i}$ containing $n$'s as the top $k$-ribbon in column $i+1$, as does $\hat{P}$.  
The part of the $P$ tableau to the right of the $(i+1)$st column remains the same.
\vspace{.1in}

In both subcases, we can now remove the $k$-ribbon containing $n$ from the $(i+1)$st column of $\hat{P}$ to obtain $\hat{P}_{n-1}$ of shape $\mu$.  The path tableau $\hat{P}_{n-1}$ is the path tableau obtained from the first $n-1$ rows of the square diagram, which come from the colored permutation
\[
\begin{array}{ccc}
\tau&=&\begin{array}{ccccccc}
1&2&\cdots&i-1&i&\cdots&n-1\\
x_1^{j_1}&x_2^{j_2}&\cdots&x_{i-1}^{j_{i-1}}&x_{i+1}^{j_{i+1}}&\cdots&a^{j_n}
\end{array}.
\end{array}
\]
Note that $\tau$ has length ${n-1}$.  
In order to use our inductive hypothesis, it remains to show that after one step of the evacuation of $P$, we obtain $P(\tau)$.  In the proof of Theorem \ref{Ptableau}, we proved that $P(\pi)$ is equal to a column of height 2 that has a split $k$-ribbon of height $j_n$ containing $a$'s on top of a $k$-ribbon of height $k+1-j_n$ containing $n$'s followed by $P(\sigma)$ where 
\[
\begin{array}{ccc}
\sigma&=&\begin{array}{ccccccc}
1&2&\cdots&i-1&i&\cdots&n-2\\
x_1^{j_1}&x_2^{j_2}&\cdots&x_{i-1}^{j_{i-1}}&x_{i+1}^{j_{i+1}}&\cdots&x_{n-1}^{j_{n-1}}
\end{array}.
\end{array}
\]
To obtain $P(\tau)$ we must insert $a^{j_n}$, into $P(\sigma)$.  

In {\it{Subcase a}}, $P(\sigma)$ looks like
\[
\begin{array}{cccccc}
a_2&a_3&\cdots&a_{i-1}&a_i& \\
n-1&n-2&\cdots&n-(i-2)&n-(i-1)&n-i
\end{array}
\]
and $a^{j_n}$ inserted into this tableau gives
\[
\begin{array}{ccccc}
a_1&a_2&\cdots&a_{i-1}&a_i\\
n-1&n-2&\cdots&n-(i-1)&n-i
\end{array}
\]
for the first $i$ columns and does not change the remaining tableau.  Again the shape of each column of height 2 is determined by the shape of the $k$-ribbon in the top row.  This is exactly what $P$ looks like after one step of the evacuation procedure.  

In {\it{Subcase b}}, $P(\sigma)$ looks like
\[
\begin{array}{ccccc}
a_2&a_3&\cdots&a_i&a_{n-(i+1)}\\
n-1&n-2&\cdots&n-(i-1)&n-i
\end{array}
\]
and $a^{j_n}$ inserted into this tableau gives
\[
\begin{array}{cccccc}
a_1&a_2&a_3&\cdots&a_i& \\
n-1&n-2&n-3&\cdots&n-i&a_{n-(i+1)}
\end{array}
\]
for the first $i+1$ columns and does not change the remaining tableau.  This is again exactly what $P$ looks like after one step of the evacuation procedure.  By induction, $ev(P(\tau)) = \hat{P}(\tau)$ and since $ev(P(\pi))$ and $\hat{P}(\pi)$ agree in the position of the $k$-ribbon containing $n$, then $ev(P(\pi)) = \hat{P}(\pi)$.  

\end{proof}

\begin{Theorem}
$Q(\pi)=\hat{Q}(\pi)$.  
\label{Qtableau}
\end{Theorem}

\begin{proof}
We will prove this result by induction on the size of $Q(\pi)$.  If $\pi$ is a permutation of a single element, then $\pi = 1^{j_1}$ so there is an $X^{j_1}$ in the single square in the growth diagram for $\pi$ and $\hat{Q(\pi)}$ is a $k$-ribbon of height $j_1$ with $1$'s in it.  One can easily check that the tableau $Q(\pi)$ for the insertion of this single element agrees with $\hat{Q(\pi)}$.  

Now assume that $Q(\sigma)= \hat{Q}(\sigma)$ for $\sigma$ a $k$-colored permutation of length $l < n$ and let $\pi$ be a $k$-colored permutation of length $n$.  Since the growth diagram for $\pi^{-1}$ is simply the reflection of the growth diagram for $\pi$ around the diagonal line $y=x$, then $\hat{P}(\pi)=\hat{Q}(\pi^{-1})$ and $\hat{Q}(\pi)=\hat{P}(\pi^{-1})$.  Let $\pi_{n-1}$ be the first $n-1$ elements in the $k$-colored permutation $\pi$.  Then $\hat{Q}(\pi_{n-1})$ is the shape of the $\hat{Q}(\pi)$ tableau at the $(n-1)$st stage.  

By reflecting across the diagonal, we have $\hat{Q}(\pi_{n-1})=\hat{P}_{n-1}(\pi^{-1})$ where $\hat{P}_{n-1}(\pi^{-1})$ is the tableau for the square diagram consisting of the first $n-1$ rows of the square diagram for $\pi^{-1}$.  Let $R$ represent the $k$-ribbon that $\hat{P}_{n-1}(\pi^{-1})$ and $\hat{P}(\pi^{-1})$ differ by.  Then $R$ also represents the $k$-ribbon that $\hat{Q}(\pi_{n-1})$ and $\hat{P}(\pi^{-1})$ differ by.  From the growth diagrams we know that the shape of $\hat{Q}(\pi)$ is equal to the shape of $\hat{P}(\pi)$ which is equal to the shape of $\hat{P}(\pi^{-1})$.  Thus $R$ represents the $k$-ribbon that $\hat{Q}(\pi)$ and $\hat{Q}(\pi_{n-1})$ differ by.  We must show that $R$ is also the $k$-ribbon that $Q(\pi)$ and $Q(\pi_{n-1})$ differ by, which means $R$ must be the $k$-ribbon that $P(\pi)$ and $P(\pi_{n-1})$ differ by, since $Q$ represents a recording tableau for the insertion tableau $P$.

Since $ev(P(\pi))=\hat{P}(\pi)$ then these tableau have the same shape and $ev(P(\pi))$ has the same shape as $P(\pi)$ by construction so $P(\pi)$ has the same shape as $\hat{P}(\pi)$.  By reflection, the shape of $\hat{P}(\pi)$ is the same as the shape of $\hat{P}(\pi^{-1})$.  Suppose $a^{j_n}$ is the last element in the $k$-colored permutation $\pi$ and let $\sigma$ be the $k$-colored permutation in $S_{n-1}$ obtained from $\pi$ by deleting the last element $a^{j_n}$ and replacing all elements $i^{j_i}$ with $i>a$ by $(i-1)^{j_i}$.  By the method of insertion, the shape of $P(\pi_{n-1}) = P(\sigma)$ and by Fomin's growth diagram we have that the shape of $\hat{Q}(\pi_{n-1})$ will be the same as the shape of $\hat{Q}(\sigma)$, since $\hat{Q}(\pi_{n-1})$ is the path tableau for the growth diagram of the first $n-1$ columns of the square diagram for $\pi$, i.e. for all but the last element $a^{j_n}$ of the square diagram for $\pi$.  

By definition, the shape of $P(\sigma)$ is equal to the shape of $Q(\sigma)$ and by induction $Q(\sigma) = \hat{Q}(\sigma)$.  Thus the $k$-ribbon that $\hat{Q}(\pi_{n-1})$ and $\hat{Q}(\pi)$ differ by is the same as the $k$-ribbon that $P(\pi_{n-1})$ and $P(\pi)$ differ by, which is the same as the $k$-ribbon that $Q(\pi_{n-1})$ and $Q(\pi)$ differ by since $Q$ is a recording tableau for $P$.  By induction, $Q(\pi_{n-1}) = \hat{Q}(\pi_{n-1})$ and since $Q(\pi_{n-1})$ differs from $Q(\pi)$ in the same $k$-ribbon that $\hat{Q}(\pi_{n-1})$ differs from $\hat{Q}(\pi)$ by, then $Q(\pi) = \hat{Q}(\pi)$.

\end{proof}

\begin{Theorem}
The evacuation procedure is a bijection between standard $k$-ribbon Fibonacci tableaux and $k$-ribbon Fibonacci path tableaux.  
\end{Theorem}
\begin{proof}
The evacuation algorithm is, by definition, an injection from standard $k$-ribbon Fibonacci tableaux to $k$-ribbon Fibonacci path tableaux.  The growth diagrams of Fomin show that $k^n n!$ equals the number of pairs $(\hat{P},\hat{Q})$ where $\hat{P}$ and $\hat{Q}$ are $k$-ribbon Fibonacci path tableaux of the same shape.  The insertion algorithm given in Section \ref{insertion} shows $k^n n!$ equals the number of pairs $(P, Q)$ where $P$ is a standard $k$-ribbon Fibonacci tableau and $Q$ is a $k$-ribbon Fibonacci path tableau.  Since $Q=\hat{Q}$ by Theorem \ref{Qtableau}, then the number of standard $k$-ribbon Fibonacci tableaux must equal the number of $k$-ribbon Fibonacci path tableaux.  Hence, the evacuation algorithm is a bijection.  
\end{proof}

\section{Color-to-Spin}
\label{colortospin}

For a pair $(P,Q)$ in which $P$ is a standard $k$-ribbon Fibonacci tableau and $Q$ is a $k$-ribbon Fibonacci path tableau, we define 
\[
vert(P,Q) = \sum_i (j_i-1)
\]
where the sum is over all $k$-ribbons in $P$ and $Q$ and $j_i$ is the height of each ribbon.

If $c$ is a column containing two $k$-ribbons, say one of height $j_1$ on top of one of height $k+1-j_1$, then the contribution of this column to $vert$ is $(j_1 - 1) + (k+1-j_1 - 1) = k-1$.  Thus the $vert$ statistic can be determined simply from the shape of $P$ as $(k-1)$ times the number of columns with two $k$-ribbons plus the sum of $(j_i - 1)$ for all columns containing a single $k$-ribbon of height $j_i$.  Since $P$ and $Q$ have the same shape then $vert(P,Q)$ is simply twice the value of the $vert$ statistic for $P$.

For each column $c$ of height 2 in $P$ containing a $k$-ribbon of height $j_1$ on top of a $k$-ribbon of height $k+1-j_1$, we define $split_c = k-j_1$.  We then define $split_P$ to be the sum of $split_c$ for all columns of height 2 in $P$.  We define $split_Q$ similarly for $Q$.  Then

\[
split(P,Q)=split_Q-split_P.
\]

We can determine the $split$ statistic from the shadow lines of the square diagram, since these lines determine the $P$ tableau for the pair $(P,Q)$.  Shadow lines with a single $X^{j_i}$ on them do not contribute to the $split$ statistic since they correspond to columns of height 1.  For a line with 2 $X$'s on it, the contribution to $split(P)$ is $k-j_2$ where $j_2$ is the color of the $X$ in the rightmost column of the shadow line.  This is a direct result of the manner in which the $P$ tableau can be constructed from the shadow lines.

The contribution of this line to $split(Q)$ is a bit trickier to determine.  When determining the $Q$ tableau in the insertion process, $k$-ribbons never move from being the bottom $k$-ribbon in a column to being the top $k$-ribbon in a column.  Thus in the recording tableau $Q$, the bottom $k$-ribbon of a column of height 2 is created first and its shape is recorded permanently in the $Q$ tableau.  In addition, once a column contains 2 $k$-ribbons, it will have 2 $k$-ribbons throughout the remainder of the insertion process.  The leftmost $X^{j_1}$ on a line with 2 $X$'s determines the shape of the bottom $k$-ribbon in the column in the $Q$-tableau.  Then the top ribbon in this column in $Q$ has height $k+1-j_2$.  Thus the contribution of this column to $split(Q)$ is $k-(k+1-j_2) = j_2-1$.

We now define
\[
spin(P,Q) = \frac{1}{2}vert(P,Q) + split(P,Q) = vert(P) + split(P,Q).
\]

In the previous example of $(P,Q)$

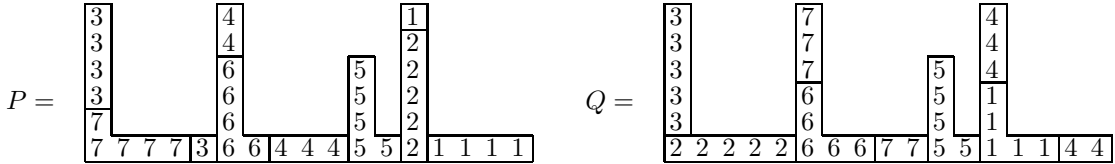
\setlength{\unitlength}{.7cm}
\begin{center}
\begin{picture}(18,3.5)(0,0)
\thinlines

\put(-1.5,1){$P=$}

\put(0,0){\line(1,0){8.5}}
\put(0,0){\line(0,1){3}}
\put(.5,.5){\line(1,0){2}}
\put(0,1){\line(1,0){.5}}
\put(0,3){\line(1,0){.5}}
\put(.5,.5){\line(0,1){2.5}}
\put(2,0){\line(0,1){.5}}
\put(2.5,0){\line(0,1){3}}
\put(2.5,2){\line(1,0){.5}}
\put(2.5,3){\line(1,0){.5}}
\put(3,.5){\line(0,1){2.5}}
\put(3,.5){\line(1,0){2}}
\put(3.5,0){\line(0,1){.5}}
\put(5,0){\line(0,1){2}}
\put(5,2){\line(1,0){.5}}
\put(5.5,.5){\line(0,1){1.5}}
\put(5.5,.5){\line(1,0){.5}}
\put(6,0){\line(0,1){3}}
\put(6,2.5){\line(1,0){.5}}
\put(6,3){\line(1,0){.5}}
\put(6.5,0){\line(0,1){3}}
\put(6.5,.5){\line(1,0){2}}
\put(8.5,0){\line(0,1){.5}}

\put(.1,.1){$7$}
\put(.1,.6){$7$}
\put(.1,1.1){$3$}
\put(.1,1.6){$3$}
\put(.1,2.1){$3$}
\put(.1,2.6){$3$}
\put(.6,.1){$7$}
\put(1.1,.1){$7$}
\put(1.6,.1){$7$}
\put(2.1,.1){$3$}
\put(2.6,.1){$6$}
\put(2.6,.6){$6$}
\put(2.6,1.1){$6$}
\put(2.6,1.6){$6$}
\put(2.6,2.1){$4$}
\put(2.6,2.6){$4$}
\put(3.1,.1){$6$}
\put(3.6,.1){$4$}
\put(4.1,.1){$4$}
\put(4.6,.1){$4$}
\put(5.1,.1){$5$}
\put(5.1,.6){$5$}
\put(5.1,1.1){$5$}
\put(5.1,1.6){$5$}
\put(5.6,.1){$5$}
\put(6.1,.1){$2$}
\put(6.1,.6){$2$}
\put(6.1,1.1){$2$}
\put(6.1,1.6){$2$}
\put(6.1,2.1){$2$}
\put(6.1,2.6){$1$}
\put(6.6,.1){$1$}
\put(7.1,.1){$1$}
\put(7.6,.1){$1$}
\put(8.1,.1){$1$}

\put(9.5,1){$Q=$}

\put(11,0){\line(1,0){8.5}}
\put(11,0){\line(0,1){3}}
\put(11.5,.5){\line(1,0){2}}
\put(11,.5){\line(1,0){.5}}
\put(11,3){\line(1,0){.5}}
\put(11.5,.5){\line(0,1){2.5}}

\put(13.5,0){\line(0,1){3}}
\put(13.5,1.5){\line(1,0){.5}}
\put(13.5,3){\line(1,0){.5}}
\put(14,.5){\line(0,1){2.5}}
\put(14,.5){\line(1,0){2}}
\put(15,0){\line(0,1){.5}}

\put(16,0){\line(0,1){2}}
\put(16,2){\line(1,0){.5}}
\put(16.5,.5){\line(0,1){1.5}}
\put(16.5,.5){\line(1,0){.5}}
\put(17,0){\line(0,1){3}}
\put(17,1.5){\line(1,0){.5}}
\put(17,3){\line(1,0){.5}}
\put(17.5,.5){\line(0,1){2.5}}
\put(17.5,.5){\line(1,0){2}}
\put(18.5,0){\line(0,1){.5}}
\put(19.5,0){\line(0,1){.5}}

\put(11.1,.1){$2$}
\put(11.1,.6){$3$}
\put(11.1,1.1){$3$}
\put(11.1,1.6){$3$}
\put(11.1,2.1){$3$}
\put(11.1,2.6){$3$}
\put(11.6,.1){$2$}
\put(12.1,.1){$2$}
\put(12.6,.1){$2$}
\put(13.1,.1){$2$}
\put(13.6,.1){$6$}
\put(13.6,.6){$6$}
\put(13.6,1.1){$6$}
\put(13.6,1.6){$7$}
\put(13.6,2.1){$7$}
\put(13.6,2.6){$7$}
\put(14.1,.1){$6$}
\put(14.6,.1){$6$}
\put(15.1,.1){$7$}
\put(15.6,.1){$7$}
\put(16.1,.1){$5$}
\put(16.1,.6){$5$}
\put(16.1,1.1){$5$}
\put(16.1,1.6){$5$}
\put(16.6,.1){$5$}
\put(17.1,.1){$1$}
\put(17.1,.6){$1$}
\put(17.1,1.1){$1$}
\put(17.1,1.6){$4$}
\put(17.1,2.1){$4$}
\put(17.1,2.6){$4$}
\put(17.6,.1){$1$}
\put(18.1,.1){$1$}
\put(18.6,.1){$4$}
\put(19.1,.1){$4$}

\end{picture}
\end{center}

we have $vert(P,Q) = 30$, $split_P = 8$, $split_Q=4$, $split(P,Q)=-4$ and $spin(P,Q) =  15-4=11$.  

For a colored permutation $\pi$ we define 
\[
color(\pi) = \sum_i (j_i-1)
\]
where $j_i$ is the color of each element $x_i$ in $\pi$.  

\begin{Theorem}
If $\pi$ is a $k$-colored permutation, colored by the colors $1$ through $k$ and $(P,Q)$ is the pair of tableaux obtained through the $k$-ribbon Fibonacci insertion algorithm, then 
\[
color(\pi) = spin(P,Q).
\]
\end{Theorem}

\begin{proof}

Since every $X^{j_i}$ in the square diagram lies on some shadow line, we will prove this result by showing that the contribution of each shadow line to $color(\pi)$ is the same as the contribution of that shadow line to $spin(P,Q)$.

Suppose the shadow line $L$ contains a single $X^{j_i}$ for some color $j_i$.  Then the contribution of this line to $color(\pi)$ is $j_i-1$.  The contribution of this line to $vert(P)$ is $j_i-1$ since this single $X^{j_i}$ corresponds to a column with a single $k$-ribbon of height $j_i$.  Since there is only one $X$ in this column, the contribution of this $X$ to $split(P,Q)$ is zero and so the contribution to $spin(P,Q)$ is $j_i-1$.  Thus the contribution of this line to $color(\pi)$ equals the contribution of the line to $spin(P,Q)$.

Suppose the shadow line $L$ contains two $X$'s, with the leftmost being $X^{j_1}$ and the rightmost being $X^{j_2}$.  The contribution of this line to $color(\pi)$ is $(j_1 - 1) + (j_2 - 1) = j_1 + j_2 - 2$.  The contribution of this line to $vert(P)$ is $k-1$ since this line corresponds to a column with two $k$-ribbons.  The contribution of this line to $split(P)$ is $k-j_2$ since the rightmost $X$ on the line, which determines the height of the top $k$-ribbon in the column, is an $X^{j_2}$.  The contribution of this line to $split(Q)$ is $k-(k+1-j_1) = j_1-1$.  The total contribution of this line to $spin(P,Q)$ is thus $(k-1) + [(j_1-1)-(k-j_2)] = (k-1) + (j_1 + j_2 - 1 - k) = j_1 + j_2 -2$ which is the same as the contribution to $color(\pi)$.

\end{proof}

\section{Knuth Relations for Fibonacci Tableaux}
\label{knuth}

In the setting of Young's lattice, the Schensted correspondence provides a nice relationship between permutations and pairs of chains in the lattice which can be interpreted as standard Young tableaux.  Furthermore, Knuth relations provide an interesting equivalence relation on $S_n$ by determining the set of all permutations that give the same $P$-tableau through the Schensted correspondence.  In this section, we explore the Fibonacci analogue of these Knuth relations, as generated by the insertion algorithm for $k$-ribbon Fibonacci tableaux described in the preceding sections.

We first define the notion of Knuth equivalence for permutations as it relates to Young tableaux.  Two permutations, $\pi$ and $\sigma$ in $S_n$ are said to be $P$-equivalent if they give the same $P$ tableau under the Schensted correspondence.  An alternate description of $P$-equivalence in terms of the permutations themselves was given by Knuth \cite{Knu}.  Two permutations $\pi$ and $\sigma$ in $S_n$ differ by a {\it {Knuth relation of the first kind}}, written $\pi \stackrel{1}{\approxeq} \sigma$, if for $x < y < z$ and
\[
\pi = x_1 \dots yxz \dots x_n \ \ {\text {and}}\ \  \sigma = x_1 \dots yzx \dots x_n\ \  {\text {or vice versa.}}
\]

Two permutations $\pi$ and $\sigma$ in $S_n$ differ by a {\it {Knuth relation of the second kind}}, written $\pi \stackrel{2}{\approxeq} \sigma$, if for $x < y < z$ and
\[
\pi = x_1 \dots xzy \dots x_n\ \  {\text {and}}\ \  \sigma = x_1 \dots zxy \dots x_n\ \  {\text {or vice versa.}}
\]

We say that two permutations $\pi$ and $\sigma$ are {\it {Knuth equivalent}}, written $\pi \stackrel{K}{\approxeq} \sigma$ if there is a sequence of permutations $\alpha_i$ such that
\[
\pi = \alpha_1 \stackrel{i_1}{\approxeq} \alpha_2 \stackrel{i_2}{\approxeq} \cdots \stackrel{i_l}{\approxeq} \alpha_{l+1}=\sigma
\]
where $i_j \in \{ 1, 2 \}$ for all $j=1, \dots, l$.
It is a fundamental theorem of Knuth \cite{Knu} that two permutations are $P$-equivalent if and only if they are Knuth equivalent.

In the context of the Fibonacci tableaux, we will say that two $k$-colored permutations are Fibonacci $P$-equivalent if they give the same $P$ tableau under the $k$-ribbon Fibonacci insertion algorithm.  However, the Knuth relations on the permutations do not preserve Fibonacci $P$ equivalence, i.e. two $k$-colored permutations may be Knuth equivalent but give different $P$ tableau under the $k$-ribbon Fibonacci insertion.  We will now describe the set of $k$-colored permutations that do give the same Fibonacci $P$ tableau under $k$-ribbon Fibonacci insertion.

In Section \ref{geometric}, we saw that given a $k$-colored permutation $\pi$, we can use shadow lines to directly determine the $P$ tableau created by the insertion algorithm.  With this tool in hand, we can now determine the set of $k$-colored permutations which are identified with a fixed tableau $P$ under the insertion correspondence, and hence describe a different equivalence relation on the set of $k$-colored permutations.

Let $L$ be a shadow line in the square diagram that contains two $X$'s, a leftmost $X^{j_1}$ in row $l$ and a rightmost $X^{j_2}$ in column $m$ and row $r$.  Then $X^{j_1}$ in row $l$ corresponds to the element $l^{j_1}$ in the permutation and the $X^{j_2}$ corresponds to the element $r^{j_2}$.  Since the shadow line $L$ would be the same if $X^{j_1}$ were in any of the columns 1 through $m-1$ as long as it is to the left of the column in which $X^{j_2}$ appears, then $l^{j_1}$ can appear in any of the positions 1 through $m-1$ in the permutation as long as it is to the left of the element $r$.  The set of all such equivalences for any shadow line in the square diagram with two $X$'s on them form the set of all permutations which give the same $P$ tableau under the $k$-ribbon Fibonacci insertion algorithm.  

For example, consider the permutation $\pi=4^3 5^1 2^1 1^4 3^2$ with the following square diagram.

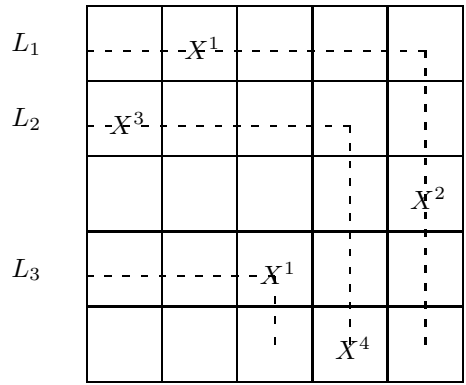
\setlength{\unitlength}{1cm}
\begin{center}
\begin{picture}(7,4.5)(1,1)
\thinlines

\put(1,0){\line(0,1){5}}
\put(2,0){\line(0,1){5}}
\put(3,0){\line(0,1){5}}
\put(4,0){\line(0,1){5}}
\put(5,0){\line(0,1){5}}
\put(6,0){\line(0,1){5}}

\put(1,0){\line(1,0){5}}
\put(1,1){\line(1,0){5}}
\put(1,2){\line(1,0){5}}
\put(1,3){\line(1,0){5}}
\put(1,4){\line(1,0){5}}
\put(1,5){\line(1,0){5}}

\put(1.3, 3.3){$X^3$}
\put(2.3, 4.3){$X^1$}
\put(3.3, 1.3){$X^1$}
\put(4.3, .3){$X^4$}
\put(5.3, 2.3){$X^2$}

\dashline{.1}(1,4.4)(5.5,4.4) 
\dashline{.1}(5.5,4.4)(5.5,0.5) 
\put(0,4.4){$L_1$}

\dashline{.1}(1,3.4)(4.5,3.4) 
\dashline{.1}(4.5,3.4)(4.5,0.5) 
\put(0,3.4){$L_2$}

\dashline{.1}(1,1.4)(3.5,1.4) 
\dashline{.1}(3.5,1.4)(3.5,0.5) 
\put(0,1.4){$L_3$}

\end{picture}
\end{center}
\vspace{.5in}

The set of all $5$-colored permutations that give the same $P$ tableau as $\pi=4^3 5^1 2^1 1^4 3^2$ under the $k$-ribbon insertion algorithm is:
\[
\begin{array}{ccc}
5^1 4^3 2^1 1^4 3^2& & 5^1 2^1 4^3 1^4 3^2\\
 4^3 5^1 2^1 1^4 3^2& & 2^1 5^1 4^3 1^4 3^2\\
 4^3 2^1 5^1 1^4 3^2& & 2^1 4^3 5^1 1^4 3^2\\
 4^3 2^1 1^4 5^1 3^2& & 2^1 4^3 1^4 5^1 3^2.
\end{array}
\]

Notice that this set is independent of the value of $k$, in the sense that it would remain unchanged provided $n=5$ remains fixed and $k$ is any value at least 4.

In the case of Young tableaux, Knuth was able to give a description of the Schensted algorithm for semistandard permutations which produces a pair of semistandard Young tableaux.  One question for further research is to determine the correct Fibonacci insertion algorithm for semistandard permutations and to find the proper definition of a semistandard Fibonacci tableau.

\end{document}